\documentclass[%
  12pt, 
  noamsfonts, 
]{amsart}
\usepackage[T1]{fontenc}
\usepackage[%
  varg, 
]{newtx}
\usepackage{cancel}
\usepackage{aliascnt}
\usepackage{mathtools}
\usepackage{physics2}
\usephysicsmodule{ab,ab.braket,xmat}
\usepackage{keytheorems}
\usepackage[%
  shortlabels,
  inline,
]{enumitem} 
\usepackage[%
  scr=boondoxupr, 
  cal=cm, 
  bb=dsfontserif, 
]{mathalpha}
\usepackage[Symbol]{upgreek}
\usepackage{xfrac}
\usepackage{accents}
\usepackage[dvipsnames]{xcolor} 
\usepackage{fixdif}
\usepackage{mathcommand} 
\usepackage{tikz} 
\usepackage[%
  pagebackref, 
]{hyperref}
\usepackage{zref-clever}
\usepackage{caption} 
\usepackage{aligned-overset} 
\numberwithin{equation}{section}
\mathtoolsset{showonlyrefs}
\setlist[enumerate,1]{
  label = \textup{(\alph*)},
  ref = \textup{(\alph*)}
}

\zcsetup{%
  cap,
  abbrev,
  nameinlink,
  counterresetby={
    theorem=section,
    proposition=section,
    lemma=section,
    corollary=section,
    conjecture=section,
    definition=section,
    remark=section,
    notation=section,
    example=section,
    warning=section,
    construction=section,
    convention=section,
    recollection=section
  },
}
\hypersetup{%
  bookmarksnumbered=true,
  colorlinks=true,
  linkcolor=blue,
  citecolor=cyan,
  urlcolor=cyan,
  pdfstartview={FitBH null},
}
\renewcommand*{\backrefalt}[4]{%
  \ifcase #1 %
  No citations.%
  \or
  $\uparrow$ #2%
  \else
  $\uparrow$ #2%
  \fi
}
\providecommand{\noopsort}[1]{} 

\newkeytheorem{mainthm}[%
  name=Theorem,%
  counter-format=\Alph{mainthm},%
]

\newkeytheorem{%
  theorem,%
  proposition,%
  lemma,%
  corollary,%
  conjecture,%
}[sharenumber=equation]

\newkeytheorem{%
  definition,%
  remark,%
  notation,%
  example,%
  warning,%
  construction,%
  convention,%
  recollection,%
}[sharenumber=equation,style=remark]

\makeatletter
\zref@newprop{listitem}{}
\NewDocumentCommand \ListInThm { m O{\textup} O{(\roman{enumi})} O{.} }
{
  \AtBeginEnvironment{#1}
  {
    \zcsetup{
      countertype={enumi=#1},
      counterresetby={enumi=#1}
    }
    \zref@setcurrent{listitem}{#2{#3}}
    \zref@localaddprop\ZREF@mainlist{listitem}
    \setlist[enumerate,1]
    {
      label={#2{#3}},
      ref={\csname the#1\endcsname#4{#3}}
    }%
  }
}
\makeatother
\NewDocumentCommand \itemref { m }
{\zcref[noname,ref=listitem]{#1}}

\ListInThm{mainthm}
\ListInThm{theorem}
\ListInThm{proposition}
\ListInThm{lemma}
\ListInThm{corollary}
\ListInThm{conjecture}

\usetikzlibrary{%
  cd,%
  arrows,%
  backgrounds,%
  calc,%
  decorations,%
  decorations.pathmorphing,%
  shapes,%
  tikzmark%
}
\tikzcdset{%
  scale cd/.style = %
  {%
    every label/.append style = {scale = #1},%
    cells = {nodes = {scale = #1}}%
  }%
}

\NewDocumentCommand \concept {m}{\textbf{#1}}
\NewDocumentCommand \nota {m}{{\color{orange}#1}}
\NewDocumentCommand \nocolor {m} {{\color{black}#1}}
\DeclareDocumentMathCommand\txand%
{ O{\quad} }{#1\text{and}#1}
\DeclareDocumentMathCommand\txforall%
{ O{\quad} }{#1\text{for all }}
\declaremathcommand{\placeholder}{\:\cdot\:}

\providecommand\given{}
\DeclareDocumentMathCommand{\SetSymbol}{O{}}{\nonscript\:#1\vert\allowbreak\nonscript\:\mathopen{}}

\DeclarePairedDelimiterX{\SetcAux}[1]{\lbrace}{\rbrace}{%
  \renewcommand\given{\SetSymbol[\delimsize]}#1%
}
\DeclarePairedDelimiterX{\SetaAux}[1]{\langle}{\rangle}{%
  \renewcommand\given{\SetSymbol[\delimsize]}#1%
}
\DeclarePairedDelimiterX{\SetbAux}[1]{\lbrack}{\rbrack}{%
  \renewcommand\given{\SetSymbol[\delimsize]}#1%
}
\DeclarePairedDelimiterX{\SetpAux}[1]{\lparen}{\rparen}{%
  \renewcommand\given{\SetSymbol[\delimsize]}#1%
}
\DeclarePairedDelimiterX{\SetbbAux}[1]{.}{.}{%
  \renewcommand\given{\SetSymbol[\delimsize]}%
  \mathopen{}\delimsize\lbrack\mkern-4mu\delimsize\lbrack
  #1%
  \delimsize\rbrack\mkern-4mu\delimsize\rbrack\mathclose{}%
}
\DeclarePairedDelimiterX{\SetccAux}[1]{.}{.}{%
  \renewcommand\given{\SetSymbol[\delimsize]}%
  \mathopen{}\delimsize\lbrace\mkern-4mu\delimsize\lbrace
  #1%
  \delimsize\rbrace\mkern-4mu\delimsize\rbrace\mathclose{}%
}
\DeclarePairedDelimiterX{\SetppAux}[1]{.}{.}{%
  \renewcommand\given{\SetSymbol[\delimsize]}%
  \mathopen{}\delimsize\lparen\mkern-4mu\delimsize\lparen
  #1%
  \delimsize\rparen\mkern-4mu\delimsize\rparen\mathclose{}%
}
\ExplSyntaxOn
\cs_new_protected:Npn \agv_set_dispatch:nnnn #1#2#3#4
{
  \group_begin:
  \tl_if_novalue:nTF {#1}
  { \cs_set_eq:NN \__agv_set_base \SetcAux }
  {
    \str_case:nnF {#1}
    {
      {a}  { \cs_set_eq:NN \__agv_set_base \SetaAux  }
      {b}  { \cs_set_eq:NN \__agv_set_base \SetbAux  }
      {c}  { \cs_set_eq:NN \__agv_set_base \SetcAux  }
      {p}  { \cs_set_eq:NN \__agv_set_base \SetpAux  }
      {bb} { \cs_set_eq:NN \__agv_set_base \SetbbAux }
      {cc} { \cs_set_eq:NN \__agv_set_base \SetccAux }
      {pp} { \cs_set_eq:NN \__agv_set_base \SetppAux }
    }
    {
      \PackageError{agv-perm-draft}
      {Unknown~\string\Set~delimiter~selector~`#1'}
      {Use~a,b,c,p,bb,cc,or~pp.}
    }
  }
  \IfBooleanTF{#2}
  { \__agv_set_base * {#4} }
  {
    \IfNoValueTF{#3}
    { \__agv_set_base {#4} }
    { \__agv_set_base [#3] {#4} }
  }
  \group_end:
}
\DeclareDocumentMathCommand{\AgvSetDispatch}{mmmm}
{ \agv_set_dispatch:nnnn {#1}{#2}{#3}{#4} }
\ExplSyntaxOff
\DeclareDocumentMathCommand{\Set}{e{_} s o m}
{\AgvSetDispatch{#1}{#2}{#3}{#4}}
\DeclareDocumentMathCommand{\dbrack}{m}{\Set_{bb}{#1}}
\DeclareDocumentMathCommand{\dbrace}{m}{\Set_{cc}{#1}}
\DeclareDocumentMathCommand{\dparen}{m}{\Set_{pp}{#1}}

\DeclareDocumentMathCommand \fun { m e{^_} d() }
{%
  \operatorname{#1}%
  \IfValueT{#2}{\sp{#2}}%
  \IfValueT{#3}{\sb{#3}}%
  \IfNoValueTF{#4}{}{\ab(#4)}%
}
\declaremathcommand{\N}{\mathbb{N}}
\declaremathcommand{\Z}{\mathbb{Z}}
\declaremathcommand{\Q}{\mathbb{Q}}
\declaremathcommand{\R}{\mathbb{R}}
\declaremathcommand{\C}{\mathbb{C}}
\declaremathcommand{\F}{\mathbb{F}}
\declaremathcommand{\k}{\mathbb{k}}
\DeclareFontFamily{U}{mathx}{}
\DeclareFontShape{U}{mathx}{m}{n}{<-> mathx10}{}
\DeclareFontSubstitution{U}{mathx}{m}{n}
\DeclareSymbolFont{agvlargebox}{U}{mathx}{m}{n}
\DeclareMathSymbol{\bigboxtimes}{\mathop}{agvlargebox}{"D2}

\DeclareMathOperator{\End}{End}

\DeclareMathOperator{\gr}{gr}

\DeclareMathOperator{\Hom}{Hom}

\DeclareMathOperator{\id}{id}

\DeclareMathOperator{\Ker}{Ker}

\DeclareMathOperator{\Res}{Res}

\DeclareMathOperator{\Spec}{Spec}

\DeclareMathOperator{\wt}{wt}
\renewmathcommand\le\leqslant
\renewmathcommand\ge\geqslant
\declaremathcommand\one{{\scriptstyle\mathfrak{I}}}
\declaremathcommand\unit{\mathbb{1}}
\declaremathcommand\iu{\mathtt{i}}
\declaremathcommand\e{\mathsf{e}}
\declaremathcommand\Rf{\operatorname{\mathbf{R}}}
\declaremathcommand\qbar{\mathscr{q}}
\declaremathcommand\SG{\mathfrak{S}}
\declaremathcommand\vac{\mathbf{1}}
\declaremathcommand\cch{\mathscr{c}}
\declaremathcommand\cfv{\upomega}
\declaremathcommand\hollowcolon{{}^{\circ}_{\circ}}
\declaremathcommand\pt{\mathsf{p}}
\declaremathcommand\tpt{\tilde{\mathsf{p}}}

\declaremathcommand{\ctensor}{\mathbin{\widehat{\otimes}}}
\DeclareDocumentMathCommand\odv{m}%
{\frac{\d}{\d{#1}}}%
\DeclareDocumentMathCommand\pdv{m}%
{\frac{\partial}{\partial{#1}}}
\DeclareDocumentMathCommand\dual{m}
{\ifblank{#1}{(\:\cdot\:)}{#1}^{\ast}}
\DeclareDocumentMathCommand\rldual{m}
{\ifblank{#1}{(\:\cdot\:)}{#1}^{\dagger}}
\DeclareDocumentMathCommand\opp{m}
{\ifblank{#1}{(\:\cdot\:)}{#1}^{\mathrm{op}}}
\DeclareDocumentMathCommand\invo{m}%
{\prescript{\theta}{}{#1}}%
\DeclareDocumentMathCommand\vect{m}%
{\boldsymbol{#1}}
\DeclareDocumentMathCommand\grp{m}%
{\fun{\mathsf{#1}}}%
\DeclareDocumentMathCommand\cat{m}%
{\operatorname{\mathsf{#1}}}%

\DeclareDocumentMathCommand\vo{mm}%
{#1_{(#2)}}%
\DeclareDocumentMathCommand\lo{mm}%
{#1_{[#2]}}%
\DeclareDocumentMathCommand\VL{m}%
{L_{#1}}%

\DeclareDocumentMathCommand\NL{O{\bullet}}%
{\mathsf{N}_{\mathsf{L}}^{#1}}
\DeclareDocumentMathCommand\NR{O{\bullet}}%
{\mathsf{N}_{\mathsf{R}}^{#1}}
\DeclareDocumentMathCommand\NLR{O{\bullet}}%
{\mathsf{N}^{#1}}
\DeclareDocumentMathCommand\cNL{O{\bullet}}%
{\prescript{\mathrm{c}}{}{\mathsf{N}}_{\mathsf{L}}^{#1}}
\DeclareDocumentMathCommand\cNR{O{\bullet}}%
{\prescript{\mathrm{c}}{}{\mathsf{N}}_{\mathsf{R}}^{#1}}
\DeclareDocumentMathCommand\cNLR{O{\bullet}}%
{\prescript{\mathrm{c}}{}{\mathsf{N}}^{#1}}

\DeclareDocumentMathCommand\normord{m}%
{\mathopen{\hollowcolon}\mathinner{#1}\mathclose{\hollowcolon}}%

\title{Twisted enveloping algebras for permutation orbifolds}
\author{Xu Gao}
\address{%
  Department of Mathematics, School of Sciences, Shanghai Maritime University,
  1550 Haigang Avenue, Pudong New Area,
  Shanghai 201306, China
}
\email{gaoxu@shmtu.edu.cn}
\author{Yiyi Zhu}
\address{%
  School of Mathematics and Statistics, Guangdong University of Technology,
  No. 161 Yinglong Road, Tianhe District, Guangzhou, 510520, China
}
\email{yzhu51@gdut.edu.cn}
\newcommand{\xu}[1]%
{\textcolor{violet}{\underline{XU:} #1 }}
\newcommand{\yiyi}[1]%
{\textcolor{purple}{\underline{YIYI:} #1 }}
\begin{document}

\begin{abstract}
  We formulate a topological formal calculus for continuous distributions with values in linearly topologized spaces.
  We apply it here to construct permutation-twisted fields directly from products and branch specializations.
  For a vertex operator algebra $V$, the key structural observation is that modes of tensors with at most one non-vacuum factor topologically generate the permutation-twisted enveloping algebra.
  Combining this observation with the calculus and $\Delta_N$, we obtain
  an explicit isomorphism $\mathscr{U}_{\sigma}(V^{\otimes N})\cong\mathscr{U}(V)$
  for an $N$-cycle $\sigma$, with grading and canonical filtrations rescaled by $N$.
  For arbitrary permutations, the cycle factors are rescaled before taking their degreewise completed tensor product.
\end{abstract}

\maketitle

\tableofcontents

\section*{Introduction}

Orbifold theory studies the representations of a fixed-point vertex operator algebra (VOA) $V^G$ through the VOA $V$ and its finite symmetry group $G$. Twisted sectors supply representations of $V^G$: for each $g\in G$, a $g$-twisted $V$-module restricts to a $V^G$-module. Permutation symmetries offer an explicit reduction: over the complex numbers, the twisted-module theory of $V^{\otimes N}$ for an $N$-cycle $\sigma\in\mathfrak{S}_N$ can be described using a single copy of $V$ \cite{BDM02}. We revisit this reduction through a topological formal calculus and the observation that modes of tensors with at most one non-vacuum factor topologically generate the twisted enveloping algebra.

\subsection*{Where the problem stands}

Permutation-twisted representations have been studied through lattice constructions \cite{BHL07,DXY16,DXY17,DXY18} and fusion products \cite{DLXY24,DXY25b,G26}. The explicit cyclic module reduction is due to Barron, Dong, and Mason \cite{BDM02}. Over the complex numbers, they use the coordinate-change operator $\Delta_N(z)$ to construct weak $\sigma$-twisted $V^{\otimes N}$-modules from weak $V$-modules and establish equivalences for weak, admissible and ordinary modules. For a finite-order automorphism $g$ acting diagonally, they obtain the corresponding equivalences with the $g^N$-twisted categories for $V$, together with the cyclic rationality equivalence. These cyclic constructions require neither rationality nor $C_2$-cofiniteness. Equivariance determines all tensor positions from one-position fields, and a local system of fields extends the formulas to arbitrary tensor states. For arbitrary permutations, \cite{BDM02} tensors the cycle constructions and, under rationality hypotheses on the factors, establishes rationality and classifies irreducible twisted modules.

Dong, Xu, and Yu \cite{DXY22} construct the cyclic Zhu algebra isomorphism $A_\sigma(V^{\otimes N})\cong A(V)$ and, for an arbitrary permutation with $k$ disjoint cycles, $A_\sigma(V^{\otimes N})\cong A(V)^{\otimes k}$, again over the complex numbers and without rationality or $C_2$-cofiniteness assumptions. Their proof uses the spanning classes of orbit sums of tensors with one non-vacuum factor and the operator $\Delta_N(1)$ to define maps in both directions and verify the quotient relations and multiplication. The module and Zhu algebra constructions thus single out the same controlling states. We study their universal modes in the completed twisted enveloping algebra.

\subsection*{What topology makes possible}

Formal calculus must accommodate infinite coefficient sums when coefficients converge to zero without eventually vanishing. Part~I treats this problem through continuous distributions with values in linearly topologized spaces. For discrete coefficient spaces, continuity recovers lower-truncated Laurent series; in several variables, the topology of the formal test-function space records the order of expansion. The formulation is independent of a choice of VOA or module, and the permutation construction below is one application. Related formal-distribution calculus in the non-Archimedean Banach setting is developed in \cite{NAVOA}; here the coefficient and test-function spaces carry linear topologies.

Universal modes of high index need not vanish. For the algebra-valued fields considered here, we establish uniform bounds ensuring that only finitely many terms contribute modulo each open left ideal. Completeness (degreewise completeness for sums in a fixed degree) then gives the coefficients of products, joint local expansions and branch specializations.

In \cite{BDM02}, one-position fields first generate a local system carrying a vertex algebra structure; a homomorphism from $V^{\otimes N}$ to this intermediate vertex algebra supplies the full twisted module structure. Here the field of an arbitrary tensor state is constructed directly by a product of coordinate-transformed fields specialized to distinct branches. Local expansions of these products give the full-state twisted Jacobi identities. The construction works both in degreewise complete graded enveloping algebras and in $\End(M)$ with its compact-open topology, treating universal coefficients and arbitrary weak modules by the same formal operations.

\subsection*{What controls the cycle reduction}

The structural observation, proved in \zcref{lem:tensor-generation}, is that modes of tensors with at most one non-vacuum factor topologically generate the permutation-twisted enveloping algebra. We call these states \emph{very simple tensors}, and write
\[
  \mathfrak{x}_a^i
  =\vac^{\otimes(i-1)}\otimes a\otimes\vac^{\otimes(N-i)}
  \qquad(a\in V,\ 1\le i\le N).
\]
Every tensor-state mode belongs to the closure of the subalgebra generated by their modes. The proof reconstructs arbitrary tensor-state fields by iterating normally ordered products and the twisted Jacobi identity, with coefficient sums converging in the canonical topology. For a cyclic twist, equivariance relates the tensor positions, so values on one position determine any continuous homomorphism on the whole algebra.

The formulas of \cite{BDM02} specify the comparison on these generators. Let $g$ have period $T$ and act diagonally, let $\sigma$ be an $N$-cycle, and put $\mathfrak{x}_a=\mathfrak{x}_a^1$. The two assignments are
\[
  \varphi\bigl(J_{n/N}(\mathfrak{x}_a)\bigr)
  =J_n(\Delta_N(1)a),
  \qquad
  \psi\bigl(J_n(a)\bigr)
  =J_{n/N}(\mathfrak{x}_{\Delta_N(1)^{-1}a})
  \quad(n\in\sfrac{1}{T}\Z).
\]
The full-state field constructions verify the defining relations; degree and neighborhood bounds give continuity, and closed-ideal descent yields homomorphisms on the completed quotients. Their composites fix the generators by cancellation of $\Delta_N(1)$ with its inverse, so topological generation proves that the maps are inverse.

The constructions work over a characteristic-zero field $\k$ containing the required roots of unity, without a rationality assumption on $V$. In the basic cyclic case $g=\id_V$, they give
\[
  \mathscr{U}_{\sigma}(V^{\otimes N})
  \cong\mathfrak{d}_N^{\ast}\mathscr{U}(V),
\]
where $\mathscr{U}(V)=\mathscr{U}_{\id_V}(V)$. The rescaling $\mathfrak{d}_N^{\ast}$ defined in \zcref{nota:grading-dilation} leaves the underlying algebra and topology unchanged and divides degrees and both left and right canonical neighborhood indices by $N$. More generally, \zcref{thm:main-iso} gives mutually inverse isomorphisms of almost-canonically seminormed rings
\[
  \varphi\colon\mathscr{U}_{\sigma g}(V^{\otimes N})
  \overset{\sim}{\longrightarrow}\mathfrak{d}_N^{\ast}\mathscr{U}_{g^N}(V),
  \qquad \psi=\varphi^{-1}.
\]

For a permutation $\sigma$ with disjoint cycle lengths $N_1,\ldots,N_k$, \zcref{thm:main-iso-general} combines the cyclic comparison and the tensor product theorem to give
\[
  \mathscr{U}_{\sigma g}(V^{\otimes N})
  \cong
  \mathfrak{d}_{N_1}^{\ast}\mathscr{U}_{g^{N_1}}(V)\ctensor\cdots\ctensor
  \mathfrak{d}_{N_k}^{\ast}\mathscr{U}_{g^{N_k}}(V).
\]
The factors are rescaled before taking the total grading and completing degreewise: original factor degrees $p_1,\ldots,p_k$ contribute total degree $\sum_i p_i/N_i$. Put $\tau=\sigma g$ and $h=g^{N_1}\otimes\cdots\otimes g^{N_k}$. For unequal cycle lengths, this completion need not coincide with $\mathscr{U}_h(V^{\otimes k})$, which uses the usual total grading.

\subsection*{What follows from the comparison}

The indexed filtration comparison identifies level algebras and the mode multiplication linking different degrees. For a single cycle, \zcref{cor:higher-zhu-isomorphism} gives
\[
  A_{\sigma g,n/N}(V^{\otimes N})\cong A_{g^N,n}(V)
  \qquad(n\in\sfrac{1}{T}\N).
\]
For arbitrary permutations, \zcref{cor:mode-transition-tensor} gives tensor product decompositions of the mode transition algebras, preserving their multiplication and the bigradings induced by the factorwise rescalings. Taking their $(0,0)$-components gives the twisted Zhu algebra decompositions of \zcref{cor:zhu-tensor}, including those of \cite{DXY22} when $g=\id_V$.

The transition-algebra comparison also transports strong identity elements and the strong identity condition, which controls reconstruction of exhaustive modules from level data. In particular, \zcref{cor:sic-general} shows that $\mathscr{U}_{\tau}(V^{\otimes N})$ satisfies this condition if and only if $\mathscr{U}_h(V^{\otimes k})$ does.

For exhaustive modules, the cyclic equivalences are pullbacks along the completed-algebra isomorphisms; the field constructions also apply to arbitrary weak modules. For a general permutation, \zcref{thm:module-equivalences-general} gives
\[
  \begin{aligned}
    \cat{WMod}_{\tau}(V^{\otimes N})
    &\simeq\bigboxtimes_{i=1}^k\cat{WMod}_{g^{N_i}}(V)
      \simeq\cat{WMod}_h(V^{\otimes k}),\\
    \cat{ExMod}_{\tau}(V^{\otimes N})
    &\simeq\bigboxtimes_{i=1}^k\cat{ExMod}_{g^{N_i}}(V)
      \simeq\cat{ExMod}_h(V^{\otimes k}),
  \end{aligned}
\]
where $\cat{WMod}$ and $\cat{ExMod}$ denote weak and exhaustive module categories, and $\boxtimes$ is the Kelly tensor product of $\k$-linear Grothendieck categories, encoding commuting actions of the cycle factors. The conformal zero mode on the corresponding $\tau$-twisted module is $\sum_i\VL{0}^{(i)}/N_i$ plus the scalar shift in \zcref{lem:virasoro-zero-mode}, where $\VL{0}^{(i)}$ denotes the Virasoro zero mode of the $i$-th factor. For unequal cycle lengths, admissibility uses rescaled degrees and the ordinary condition involves this weighted sum. When all cycle lengths agree, the usual admissible and ordinary conditions are recovered after uniform rescaling.

Finally, \zcref{cor:rationality-general} shows that $V^{\otimes N}$ is $\tau$-rational if and only if $V^{\otimes k}$ is $h$-rational, using the paper's convention of geometrically ordinary simple modules over $\k$.

\subsection*{Organization of the paper}
Part~I develops the topological formal calculus. Section~\ref{sec:topological-formal-calculus} introduces continuous distributions and their formal test-function spaces; \zcref{sec:formal-expansions-residues} develops ordered expansions and boundary residues, including fractional powers.

Part~II supplies the VOA and algebraic frameworks. Section~\ref{sec:voa-acs-framework} recalls VOAs and twisted modules, while \zcref{sec:acsr-framework} recalls almost-canonically seminormed rings, exhaustive modules, level algebras and mode transition algebras. Section~\ref{sec:voa-enveloping-algebras} constructs the twisted enveloping algebras and establishes a twisted Jacobi criterion for field assignments; \zcref{sec:enveloping-tensor-products} treats tensor products of VOAs, enveloping algebras and modules.

Part~III carries out the permutation construction. Section~\ref{sec:set-up} establishes topological generation and develops the tools for coordinate change, field products and branch specialization. Sections~\ref{sec:varphi} and~\ref{sec:psi} construct $\varphi$ and $\psi$, respectively. Section~\ref{sec:conclusion} proves the inverse identities, passes to arbitrary permutations and derives the applications.

\begin{center}
  \begin{tikzpicture}[
      section node/.style={font=\normalsize,inner sep=3pt},
      part frame/.style={draw,dashed,minimum height=3.8cm},
      >=latex
    ]
      \node[part frame,minimum width=2cm] (part-i) at (0,0) {};
      \node[part frame,minimum width=4.4cm] (part-ii) at (4.6,0) {};
      \node[part frame,minimum width=4.4cm] (part-iii) at (10.4,0) {};
      \node[font=\small\bfseries] at (0,1.45) {Part I};
      \node[font=\small\bfseries] at (4.6,1.45) {Part II};
      \node[font=\small\bfseries] at (10.4,1.45) {Part III};

      \node[section node] (s1) at (0,0.15) {\S\ref{sec:topological-formal-calculus}};
      \node[section node] (s2) at (0,-0.95) {\S\ref{sec:formal-expansions-residues}};
      \draw[->] (s1) -- (s2);

      \node[section node] (s3) at (3.5,0.55) {\S\ref{sec:voa-acs-framework}};
      \node[section node] (s4) at (5.7,0.55) {\S\ref{sec:acsr-framework}};
      \node[section node] (s5) at (4.6,-0.4) {\S\ref{sec:voa-enveloping-algebras}};
      \node[section node] (s6) at (4.6,-1.35) {\S\ref{sec:enveloping-tensor-products}};
      \draw[->] (s3) -- (s5);
      \draw[->] (s4) -- (s5);
      \draw[->] (s5) -- (s6);

      \node[section node] (s7) at (10.4,0.55) {\S\ref{sec:set-up}};
      \node[section node] (s8) at (9.2,-0.4) {\S\ref{sec:varphi}};
      \node[section node] (s9) at (11.6,-0.4) {\S\ref{sec:psi}};
      \node[section node] (s10) at (10.4,-1.35) {\S\ref{sec:conclusion}};
      \draw[->] (s7) -- (s8);
      \draw[->] (s7) -- (s9);
      \draw[->] (s8) -- (s10);
      \draw[->] (s9) -- (s10);

      \node at ($(part-i.east)!0.5!(part-ii.west)$) {$\Longrightarrow$};
      \node at ($(part-ii.east)!0.5!(part-iii.west)$) {$\Longrightarrow$};
    \end{tikzpicture}
\end{center}

\subsection*{Conventions}
Throughout this paper, $\nota{\k}$ is a field of characteristic zero. 
Unless otherwise stated, all vector spaces and tensor products are over $\k$, and all linear maps are $\k$-linear.
For simplicity, we will assume that $\k$ contains roots of unity that are necessary for the constructions.
Over a smaller field, all constructions and statements are understood after the corresponding base change.
We will always discard the (possible) topology on $\k$ and treat it as a discrete field.

\part{Topological formal calculus}\label{part:preliminaries}

Formal calculus packages infinitely many coefficient identities into a single generating-series identity. Here we formulate this calculus in terms of continuous distributions on topologized spaces of formal test functions. For a discrete coefficient space, continuity recovers lower-truncated Laurent series; for a linearly topologized coefficient space, the corresponding coefficients are required to converge to zero. In several variables, the topology of the test-function space records the order of expansion and yields uniform coefficient bounds. The distinction between different expansion orders is thus part of the formulation itself.

Under the completeness hypotheses stated below, the resulting coefficient sums can be recovered from finite sums in discrete quotients. These bounds justify the ordered products, diagonal expansions and fractional boundary residues developed in this part; for algebra-valued products, neighborhoods given by open left ideals supply the required control. The formulation does not require a choice of VOA or module. We later apply it to the universal coefficient fields in \zcref{sec:voa-enveloping-algebras,sec:set-up} and to the module fields of \zcref{rem:weak-module-distributions}. Our treatment is inspired by the Banach setting of \cite{NAVOA}; the coefficient and test-function spaces here are handled through linear topologies (cf.~\zcref{rec:topological-vector-spaces}).

\section{Formal series and continuous distributions}
\label{sec:topological-formal-calculus}

We first recall formal series and residues, then introduce continuous distributions.

\subsection{Formal series}
The basic objects are \concept{(unbounded) formal power series}, which are formal infinite sums of the form
\[
  F(z_1,\cdots,z_n) = \sum_{i_1,\cdots,i_n \in \Z} F_{i_1,\cdots,i_n} z_1^{-i_1-1} \cdots z_n^{-i_n-1},
\]
where the coefficients $F_{i_1,\cdots,i_n}$ are taken from the given vector space $U$.
The space of all such formal power series is denoted by $\nota{U\dbrack{z_1^{\pm 1}, \cdots, z_n^{\pm 1}}}$.

Recall that the \concept{residue} of a formal power series in $U\dbrack{z_1^{\pm 1}, \cdots, z_n^{\pm 1}}$ with respect to the variable $z_k$ is defined by linearly extending the rule
\[
  \nota{\Res_{z_k}} z_1^{-i_1-1} \cdots z_n^{-i_n-1} = \delta_{i_k,0}z_1^{-i_1-1} \cdots \xcancel{z_k^{-i_k-1}}\cdots z_n^{-i_n-1}.
\]
And the multiple residue is defined by iterating the single-variable one:
\[
  \Res_{z_1,\cdots,z_n} := \Res_{z_1} \circ \cdots \circ \Res_{z_n}.
\]

Geometrically, a formal power series $F(z_1,\cdots,z_n) \in U\dbrack{z_1^{\pm 1}, \cdots, z_n^{\pm 1}}$ can be viewed as an \emph{$U$-valued distribution} on the poly-punctured affine line $(\mathbb{A}^1\setminus\Set{0})^n = \Spec\k[z_1^{\pm 1}, \cdots, z_n^{\pm 1}]$: 
it acts on a Laurent polynomial $\phi(z_1,\cdots,z_n)$ via taking the residue with respect to $z_1,\cdots,z_n$:
\[
  \Res_{z_1,\cdots,z_n} \phi(z_1,\cdots,z_n)F(z_1,\cdots,z_n),
\]
where the product $\phi(z_1,\cdots,z_n)F(z_1,\cdots,z_n)$ is defined via the scale multiplication of $U$.
In this way, we obtain a linear isomorphism
\begin{equation}\label{eq:distribution}
  U\dbrack{z_1^{\pm 1}, \cdots, z_n^{\pm 1}} \cong
  \Hom_{\k}(\k[z_1^{\pm 1}, \cdots, z_n^{\pm 1}], U).
\end{equation}

\begin{notation}
  Although we view $F(z_1,\cdots,z_n)$ as a distribution, we should avoid using the notation $F(\phi)$ for its action on $\phi$ as it may be confused with the evaluation of the power series. 
  In what follows, we will employ braket notations: for instance,
  \[
    \nota{\braket < \nocolor{F} | \nocolor{\phi}  > }
    \quad\text{or}\quad
    \nota{\braket < \nocolor{F(z_1,\cdots,z_n)} | \nocolor{\phi(z_1,\cdots,z_n)}  > }
  \]
  will denote the action of $F$ on $\phi$, defined as above.
\end{notation}

\subsection{Topological vector space}

We use linear topologies to interpret the infinite sums of modes, for example, in \zcref{def:enveloping-algebra}.

\begin{recollection}\label{rec:topological-vector-spaces}
  Recall that a \concept{topological vector space} is a vector space $U$ equipped with a topology so that addition and scalar multiplication are
  continuous.  
  In \cite[\S 3.6.1]{BD}, Beilinson and Drinfeld consider topological vector spaces that are
  \begin{enumerate}
    \item \concept{linearly topologized}: open subspaces form a basis of
        neighborhoods of $0$;
    \item \concept{separated}: the intersection of all open subspaces is $0$;
        and
    \item \concept{complete}: every Cauchy net converges.
  \end{enumerate}
  Further discussion of these assumptions can be found in \cite{B08}; see also \cite{P24}.
  We also refer to \cite[Chap. II.6]{LefschetzAT} and \cite[\href{https://stacks.math.columbia.edu/tag/07E7}{Tag 07E7}]{stacks-project} for treatments of topological vector spaces.
\end{recollection}
\begin{convention}\label{conv:graded-topological-spaces}
  We will encounter graded topological vector spaces. 
  For $U=\bigoplus_{p\in\Gamma}U_p$, the linear topology is understood in the category of $\Gamma$-graded vector spaces: the neighborhoods of zero have a basis of graded subspaces.  
  The same neighborhoods topologize the underlying vector space.
  Unless otherwise specified, convergence refers to this topology. Such convergence implies degreewise convergence, but not conversely in general. 
  Hence, completeness after forgetting the grading does not follow from the degreewise completeness.
\end{convention}

\begin{convention}\label{conv:sep-comp-assumptions}
  Unless otherwise stated, all (graded) topological vector spaces are linearly topologized and separated.
\end{convention}

\begin{recollection}\label{rec:topological-associative-algebras}
  By a \concept{(graded) topological associative algebra} we mean a linearly
  topologized, separated, (graded) complete (graded) associative algebra such that the neighborhoods of zero have a basis of open left (graded) ideals and the right multiplication is continuous.
  Equivalently, it is an associative algebra object in the normally ordered monoidal category of \cite[\S 3.6.1]{BD}.
\end{recollection}
\begin{remark}\label{rem:laurent-field-not-bd}
  The field $\k\dparen{t}$ with its usual valuation topology is not a topological associative algebra in the above sense, since its basic lattices $t^n\k\dbrack{t}$ are not ideals of $\k\dparen{t}$; see \cite[\S 3.6.1, footnote 33]{BD}.
\end{remark}

\begin{example}
  Any vector space can be equipped with the \concept{discrete} topology and this topology verifies all the assumptions (since our field $\k$ is equipped with the discrete topology\footnote{Conversely, if $\k$ is not discrete, then the scalar multiplication on a nontrivial discrete vector space will never be continuous. This is why we must insist on discarding (possible) topologies on $\k$ and view it as an abstract field.}).
  \emph{The usual formal calculus can thus be embedded into the topological formal calculus.}
\end{example}
\begin{recollection}\label{rec:discrete-quotients}
  Any complete topological vector space can be realized as a projective limit of discrete spaces.\footnote{For a degreewise complete graded space, this identification holds degreewise.}
  In this sense, we may think them as pro-objects of discrete vector spaces. 
  Using this perspective, constructions and arguments can often be reduced to the discrete case by passing to the projective limit.
\end{recollection}

\begin{definition}\label{def:summable}
  Let $U$ be a topological vector space, and let $(v_i)_{i \in I}$ be a family of elements of $U$ indexed by an arbitrary set $I$.
  We say that the family $(v_i)_{i \in I}$ is \concept{summable} if the net of finite partial sums
  \[
    \sum_{i \in J} v_i,\qquad J\subset I \text{ finite, }
  \]
  converges in $U$.
  In this case, we denote the limit by $\nota{\sum_{i \in I} v_i}$.
\end{definition}
\begin{remark}
  A family in a graded topological vector space is called \concept{degreewise summable} if its projection to each homogeneous component is summable. Summability implies degreewise summability, but not conversely in general.
\end{remark}
The following criterion for summability is evident from \zcref{rec:discrete-quotients}.
\begin{lemma}\label{lem:summability}
  Let $U$ be a topological vector space and $(v_i)_{i\in I}$ a family contained in a complete subspace of $U$. Then the family is summable if and only if its image in every discrete quotient of $U$ has finite support.
\end{lemma}
\begin{corollary}
  If a sequence $(v_n)_{n\in\N}$ in a complete topological vector space converges to $0$, then it is summable\footnote{This is very different from standard real analysis.}.
  The same holds in a degreewise complete graded topological vector space when all $v_n$ belong to one fixed homogeneous component.
\end{corollary}

Another important family of topological vector spaces is the compact ones.

\begin{recollection}
  A topological vector space is said to be \concept{(linearly) compact} if every family of its closed affine subspaces with the finite intersection property has a non-empty intersection.
  A topological vector space that is both compact and discrete is finite-dimensional.
  Hence, one can characterize compact vector spaces as \emph{projective limits} of finite-dimensional vector spaces. 
  On the other hand, one can characterize discrete vector spaces as \emph{inductive limits} of finite-dimensional vector spaces.
\end{recollection}

\begin{example}\label{eg:continuous-hom}
  The space $\Hom_{\k}^{\mathtt{cont}}(U,V)$ of continuous linear maps between two topological vector spaces $U$ and $V$ is topologized with the \emph{(linearly) compact-open topology}: a basis of neighborhoods of $0$ is given by the subspaces of continuous linear maps from compact subspaces of $U$ to open subspaces of $V$.
\end{example}

\begin{remark}
  When $U$ is discrete, every linear map from it is automatically continuous. 
  Hence, $\Hom_{\k}^{\mathtt{cont}}(U,V)$ is simply $\Hom_{\k}(U,V)$. 
  However, even when $V$ is also discrete, the topology on $\Hom_{\k}(U,V)$ is not discrete in general.
  Indeed, a sequence of linear maps $f_n \in \Hom_{\k}(U,V)$ converges if and only if for every $u \in U$, the sequence $(f_n(u))_{n\in\N}$ converges in $V$ (when $V$ is discrete, this means the sequence is eventually constant).
\end{remark}

\subsection{Continuous distributions}
Imposing continuity on the distributions of \eqref{eq:distribution} gives a topological enhancement of the formal calculus.
\begin{notation}\label{nota:continuous-distributions}
  Equip $\k[z^{\pm1}]$ with the linear topology induced by the valuation topology of $\k\dparen{z}$; its basic open subspaces are $z^N\k[z]$, for $N\in\Z$.  This linear topology is separated but not complete.

  For a vector space $U$ equipped with a linear topology, we write
  \[
    \nota{\mathcal{D}_{U}(z)}
    := \Hom_{\k}^{\mathtt{cont}}(\k[z^{\pm1}],U).
  \]
  Under the residue identification, this is the subspace
  \begin{equation}\label{eq:continuous-distribution}
    \mathcal{D}_{U}(z)
    = \Set*{F(z)=\sum_{i\in\Z}F_iz^{-i-1}
      \given \lim_{i\to\infty}F_i=0}
    \subset U\dbrack{z^{\pm1}}.
  \end{equation}
  Indeed, continuity at $0$ says precisely that for every open subspace $W\subseteq U$, one has $F_i\in W$ for all sufficiently large $i$. In particular, $\mathcal{D}_{U}(z)=U\dparen{z}$ when $U$ is discrete.

  When $U$ is complete, restriction from the completed test-function space gives a canonical isomorphism
  \[
    \Hom_{\k}^{\mathtt{cont}}(\k\dparen{z},U)
    \xrightarrow{\ \sim\ }
    \mathcal{D}_U(z),
  \]
  and \zcref{rec:discrete-quotients} gives
  \[
    \mathcal{D}_{U}(z) = \varprojlim_{W} \mathcal{D}_{U/W}(z) = \varprojlim_{W} (U/W)\dparen{z},
  \]
  where $W$ ranges over open subspaces of $U$.
\end{notation}

For two variables, we first recall the topology on $\k\dparen{z}\dparen{w}$.
\begin{recollection}\label{rec:multi-valuation-topology}
  The standard \emph{multi-valuation topology} on $\k\dparen{z}\dparen{w}$ is the topology induced by the normally ordered tensor product; see \cite[\S 3.6.1]{BD}, or \cite[\S 4.1]{Morrow12} for the equal-characteristic higher local field topology. 
  Explicitly, a basis of open subspaces is given by
  \[
    \mathcal{V}(m_\bullet)
    :=\Set*{\sum_{j\in\Z}\phi_j(z)w^j\in\k\dparen{z}\dparen{w}
      \given \phi_j(z)\in z^{-m_j}\k\dbrack{z}\text{ for every }j},
  \]
  where $m_\bullet$ ranges over all sequences satisfying $m_j\in\Z\cup\{\infty\}$ and $m_j=\infty$ for $j\gg0$, and we adopt the convention $z^{-\infty}\k\dbrack{z}:=\k\dparen{z}$.
  Equivalently, one may replace $z^{-m_j}\k\dbrack{z}$ by an arbitrary open subspace $L_j\subseteq\k\dparen{z}$, subject to the condition $L_j=\k\dparen{z}$ for $j\gg0$.
\end{recollection}

\begin{notation}\label{nota:ordered-continuous-distributions}
  Give $\k[z^{\pm1},w^{\pm1}]$ the linear topology induced from the multi-valuation topology on $\k\dparen{z}\dparen{w}$.  
  For a vector space $U$ equipped with a linear topology, we write
  \[
    \nota{\mathcal{D}_{U}(z,w)}
    :=\Hom_{\k}^{\mathtt{cont}}(\k[z^{\pm1},w^{\pm1}],U).
  \]
  Thus the order of the variables is essential: $\mathcal{D}_{U}(z,w)$ and $\mathcal{D}_{U}(w,z)$ use different topologies on their test-function spaces.
  For complete $U$, restriction from the completed test-function space $\k\dparen{z}\dparen{w}$ induces a canonical isomorphism
  \[
    \Hom_{\k}^{\mathtt{cont}}(\k\dparen{z}\dparen{w},U)
    \xrightarrow{\ \sim\ }
    \mathcal{D}_U(z,w).
  \]
\end{notation}

\begin{proposition}\label{prop:ordered-continuous-distributions}
  Under the residue identification \eqref{eq:distribution}, $\mathcal{D}_{U}(z,w)$ is identified with the subspace of $U\dbrack{z^{\pm 1}, w^{\pm 1}}$ whose members are the double series
  \[
    F(z,w)=\sum_{i,j\in\Z}F_{i,j}z^{-i-1}w^{-j-1}
  \]
  satisfying the following conditions:
  \begin{enumerate}
    \item\label{prop:ordered-continuous-distributions:pointwise} 
        for every fixed $j$, one has $\lim\limits_{i\to\infty}F_{i,j}=0$ in $U$; and
    \item\label{prop:ordered-continuous-distributions:uniform} 
        for every open subspace $W\subseteq U$, there is $j_0\in\Z$ such that
        \[
          F_{i,j}\in W
          \qquad\text{for every }i\in\Z\text{ and every }j\ge j_0.
        \]
  \end{enumerate}
\end{proposition}
\begin{proof}
  Let $\Lambda\colon\k[z^{\pm1},w^{\pm1}]\to U$ be continuous and set $F_{i,j}:= \braket <\Lambda|z^iw^j>$. 
  For each fixed $j$, the coefficient-line map
  \[
    \k[z^{\pm1}]\longrightarrow\k[z^{\pm1},w^{\pm1}],\qquad
    \phi(z)\longmapsto \phi(z)w^j,
  \]
  is continuous. 
  Its composite with $\Lambda$ belongs to $\mathcal{D}_{U}(z)$, which gives \itemref{prop:ordered-continuous-distributions:pointwise} for $F_{i,j}$. 
  Given an open subspace $W\subseteq U$, choose a basic multi-valuation neighborhood $\mathcal{V}(m_\bullet)$ such that $\mathcal{V}(m_\bullet)\cap\k[z^{\pm1},w^{\pm1}]\subseteq\Lambda^{-1}(W)$. 
  Since $m_j=\infty$ for $j\gg0$, the whole polynomial coefficient-line $\k[z^{\pm1}]w^j$ lies in $\mathcal{V}(m_\bullet)$ for all sufficiently large $j$. 
  This gives \itemref{prop:ordered-continuous-distributions:uniform}, uniformly in $i$.

  Conversely, suppose we have a double series $F(z,w)$ satisfying \itemref{prop:ordered-continuous-distributions:pointwise,prop:ordered-continuous-distributions:uniform}.  
  By \itemref{prop:ordered-continuous-distributions:pointwise}, for each $j$ there is a continuous map
  \[
    \Lambda_j\colon\k[z^{\pm1}]\longrightarrow U,
    \qquad z^i\longmapsto F_{i,j}.
  \]
  For a Laurent polynomial $\phi=\sum_j\phi_j(z)w^j$, define
  \[
    \Lambda(\phi):=\sum_j\braket <\Lambda_j|\phi_j>,
  \]
  where the sum is finite.
  To see that $\Lambda$ is continuous, fix an open $W\subseteq U$ and choose $J$ such that $\Lambda_j(\k[z^{\pm1}])\subseteq W$ for $j\ge J$. 
  For each $j<J$, choose an open subspace $L_j\subseteq\k\dparen{z}$ whose intersection with $\k[z^{\pm1}]$ is mapped into $W$ by $\Lambda_j$, and put $L_j=\k\dparen{z}$ for $j\ge J$. 
  The intersection of the resulting basic multi-valuation neighborhood $\sum_jL_jw^j$ with $\k[z^{\pm1},w^{\pm1}]$ is mapped into $W$. 
  Thus $\Lambda$ is continuous and has the prescribed coefficients.
\end{proof}

\begin{remark}\label{rem:ordered-continuous-distributions}
  If $U$ is discrete, \zcref{prop:ordered-continuous-distributions} reduces to
  \[
    \mathcal{D}_{U}(z,w)=U\dparen{z}\dparen{w}.
  \]
  For complete $U$, \zcref{rec:discrete-quotients} also gives
  \[
    \mathcal{D}_U(z,w)
    = \varprojlim_W \mathcal{D}_{U/W}(z,w)
    = \varprojlim_W (U/W)\dparen{z}\dparen{w},
  \]
  where $W$ ranges over open subspaces of $U$.
\end{remark}

\begin{remark}
  For graded $U$, the distribution spaces above use the underlying topology of \zcref{conv:graded-topological-spaces}. We do not consider distribution spaces in the graded scenario since the test-function spaces are NOT graded topological vector spaces.
\end{remark}

\begin{lemma}\label{lem:products-of-continuous-distributions}
  Let $U$ be a topological associative algebra in the sense of \ref{rec:topological-associative-algebras}. If
  \[
    A(z)\in\mathcal{D}_U(z)
    \txand
    B(w)\in\mathcal{D}_U(w),
  \]
  then we have
  \[
    A(z)B(w)\in\mathcal{D}_U(z,w)
    \txand
    B(w)A(z)\in\mathcal{D}_U(w,z).
  \]
\end{lemma}
\begin{proof}
  Write $A(z)=\sum_iA_iz^{-i-1}$ and $B(w)=\sum_jB_jw^{-j-1}$. 
  For each fixed $j$, continuity of right multiplication by $B_j$ gives $A_iB_j\to0$ as $i\to\infty$. 
  Now let $W\subseteq U$ be an open left ideal.
  Since $B_j\to0$, one has $B_j\in W$ for $j\gg0$, and then
  \[
    A_iB_j\in UW\subseteq W
  \]
  for every $i$.
  Since open left ideals form a neighborhood basis, \zcref{prop:ordered-continuous-distributions} gives the desired containments.
\end{proof}
\begin{remark}
  Joint continuity of multiplication alone does not imply the lemma. 
  For example, take $U=\k\dparen{t}$ with its valuation topology, put $A_{-n}=t^{-n}$ for $n\ge1$ and $A_i=0$ otherwise, and put $B_j=t^j$ for $j\ge0$ and $B_j=0$ otherwise. 
  Then $A_i\to0$ and $B_j\to0$, but $A_{-(j+1)}B_j=t^{-1}\notin\k\dbrack{t}$ for every $j\ge0$. 
  Hence $A(z)B(w)\notin\mathcal{D}_U(z,w)$, in agreement with \zcref{rem:laurent-field-not-bd}.
\end{remark}

\section{Formal expansions and residues}
\label{sec:formal-expansions-residues}

We define ordered expansions and boundary residues in the two-points scenario, and then extend these constructions to fractional powers.

\subsection{Ordered expansions and boundary residues}
Consider the configuration schemes 
\[
  \nota{\mathfrak{X}_{2+n}}
  :=\Set*{
    (\pt_0,\pt_1,\ldots,\pt_n,\pt_\infty)
    \in(\mathbb{P}^{1})^{n+2}
    \given
    \pt_0=0,\ \pt_\infty=\infty,\
    \pt_i\neq\pt_j\text{ for }i\neq j
  }.
\]
For $n=2$, we use $w$ and $z$ as the coordinates of $\pt_1$ and $\pt_2$, respectively.  
Then, the function ring of $\mathfrak{X}_{2+2}$ is
\[
  \nota{\mathcal{R}} := \k[z^{\pm1},w^{\pm1},(z-w)^{-1}].
\]
\begin{notation}\label{nota:relative-boundary}
  We view $\mathfrak{X}_{2+2}$ as a relative affine scheme over $\mathfrak{X}_{2+1}=\Spec\k[w^{\pm1}]$.  
  Its relative compactification is
  \[
    \overline{\mathfrak{X}}_{2+2}^{\mathrm{rel}}
    :=\mathbb{P}^{1}_{\mathfrak{X}_{2+1}},
  \]
  whose relative boundary divisor consists of three components:
  \[
    z=0,\qquad
    z=\infty,\qquad\text{and}\qquad
    z=w.
  \]
  Corresponding to them, we have three \emph{ordered expansions}:
  \[
  \begin{aligned}
    \nota{\iota_{w,z}} &\colon \mathcal{R} \to \k\dparen{w}\dparen{z},\\
    \nota{\iota_{z,w}} &\colon \mathcal{R} \to \k\dparen{z}\dparen{w},\\
    \nota{\iota_{w,z-w}} &\colon \mathcal{R} \to \k\dparen{w}\dparen{z-w}.
  \end{aligned}
  \]
  Explicitely, $\iota_{w,z}$ and $\iota_{z,w}$ are expansions at $z=0$ and $z=\infty$ respectively: they keep powers of $z,w$, and send $(z-w)^{\Box}$ ($\Box \in \Z$) to
  \[
    \sum_{k=0}^\infty \binom{\Box}{k} (-w)^{\Box-k} z^k \in \k\dparen{w}\dparen{z}
    \qquad\text{and}\qquad
    \sum_{k=0}^\infty \binom{\Box}{k} z^{\Box-k} (-w)^k \in \k\dparen{z}\dparen{w}
  \]
  respectively. 
  At the diagonal $z=w$, $\iota_{w,z-w}$ keeps powers of $w,(z-w)$, and sends $z^{\Box}$ ($\Box \in \Z$) to
  \[
    \sum_{k=0}^\infty \binom{\Box}{k} w^{\Box-k} (z-w)^k \in \k\dparen{w}\dparen{z-w}.
  \]
\end{notation}

We extend the diagonal expansion from rational functions to distributions by applying the same binomial formula termwise:
\begin{lemma}\label{lem:diagonal-expansion}
  Let $U$ be a topological vector space, and let
  \[
    F(z,w)\in\mathcal{D}_U(z,w)\cap\mathcal{D}_U(w,z).
  \]
  Assume that $U$ is complete. Then termwise binomial expansion defines
  \[
    \iota_{w,z-w}F\in\mathcal{D}_U(w,z-w).
  \]
  For graded $U$, degreewise completeness is sufficient if coefficients of monomials of the same total degree are homogeneous of the same degree.
\end{lemma}
\begin{proof}
  Only monomials of total degree $p+k$ contribute to the coefficient of $w^p(z-w)^k$ in the diagonal expansion. 
  By \zcref{prop:ordered-continuous-distributions}, lower bounds for the powers of $z$ and $w$ leave only finitely many contributions modulo every open subspace. Their sum exists by \zcref{lem:summability} and completeness.
  This shows the well-definedness of the diagonal expansion $\iota_{w,z-w}$.

  For fixed $k$, these bounds also give a lower bound for the powers of $w$ modulo every open subspace. Since no negative powers of $z-w$ occur in the expansion, by \zcref{prop:ordered-continuous-distributions}, we have  $\iota_{w,z-w}F\in\mathcal{D}_U(w,z-w)$.
\end{proof}

\begin{remark}\label{rem:diagonal-poles}
  To expand a distribution with a diagonal pole, we first multiply it by a power of $z-w$ so that the lemma applies, and then divide the resulting expansion by the same power. The binomial identity shows that the result is independent of the chosen power.
\end{remark}

Recall that we have used $\Res_z$ for the extraction of the coefficient of $z^{-1}$ in a formal power series in $z$. 
We write $\oint_{z=\bullet}\cdots\d{z}$ for the boundary residues defined below. \emph{This notation is just a formal symbol in this paper.} For local expansions of complex meromorphic differentials, these operations agree with the usual residues, given by contour integrals normalized by $\frac{1}{2\pi\iu}$.
\begin{definition}\label{def:boundary-residue}
  Let $U$ be a complete topological vector space and fix $\phi(z,w)\in\mathcal{R}$.
  Using the canonical extensions in \zcref[noname]{nota:continuous-distributions,nota:ordered-continuous-distributions}, for
  \[
    F_0\in\mathcal{D}_U(w,z),\qquad
    F_\infty\in\mathcal{D}_U(z,w),\qquad
    F_\Delta\in\mathcal{D}_U(w,z-w),
  \]
  define elements of $\mathcal{D}_U(w)$ by their action on $\psi(w)\in\k\dparen{w}$:
  \begin{align*}
    \braket < \nota{\oint_{z=0}\phi F_0\d{z}} | \psi >
      &:=\braket < F_0 | \iota_{w,z}(\phi)\psi(w) >,\\
    \braket < \nota{\oint_{z=\infty}\phi F_\infty\d{z}} | \psi >
      &:=-\braket < F_\infty | \iota_{z,w}(\phi)\psi(w) >,\\
    \braket < \nota{\oint_{z=w}\phi F_\Delta\d{z}} | \psi >
      &:=\braket < F_\Delta | \iota_{w,z-w}(\phi)\psi(w) >.
  \end{align*}
  We call $\nota{\oint_{z=\bullet}\cdots\d{z}}$ the \concept{boundary residue} along $z=\bullet$.
  In the double-series notation, we can write them as follows:
  \begin{enumerate}
    \item 
        $\displaystyle \oint_{z=0}{\phi(z,w)F_0(z,w)}\d{z} = \Res_{z} \iota_{w,z}(\phi(z,w)) F_0(z,w)$;
    \item 
        $\displaystyle \oint_{z=\infty}{\phi(z,w)F_\infty(z,w)}\d{z} = -\Res_{z} \iota_{z,w}(\phi(z,w)) F_\infty(z,w)$;
    \item 
        $\displaystyle \oint_{z=w}{\phi(z,w)F_\Delta(z,w)}\d{z} = \Res_{z-w} \iota_{w,z-w}(\phi(z,w)) F_\Delta(z,w)$.
  \end{enumerate}
\end{definition}
\begin{remark}\label{rem:graded-boundary-residues}
  More generally, the coefficient formulas define boundary residues whenever the terms contributing to each coefficient form a summable family. For degreewise complete graded $U$, the bounds in \zcref{prop:ordered-continuous-distributions} ensure summability if, for each coefficient, all contributing terms are homogeneous of the same degree. The same applies to the fractional formulas in \zcref{def:boundary-residue-g}.
\end{remark}

Extending the orindary residue theorem for projective line linearly, we have:
\begin{proposition}[Residue theorem]\label{prop:residue-theorem}
  Let $F \in U\otimes_\k\mathcal{R}$, and write
  \[
    F_0:=\iota_{w,z}F,\qquad
    F_\infty:=\iota_{z,w}F,\qquad
    F_\Delta:=\iota_{w,z-w}F.
  \]
  Then, for every $\phi(z,w)\in\mathcal{R}$,
  \[
    \oint_{z=0}\phi F_0\d{z}
    {}+\oint_{z=\infty}\phi F_\infty\d{z}
    {}+\oint_{z=w}\phi F_\Delta\d{z}
    =0.
  \]
\end{proposition}

\subsection{Topological formal calculus for fractional powers}

Orbifold theory requires formal series with fractional powers. 
In what follows, we fix $T \in \N$ and a primitive $T$-th root of unity $\nota{\zeta_T}$.
\begin{convention}\label{conv:frac-power}
  Following the usual convention in orbifold theory, for a series in $U\dbrack{z^{\pm\sfrac{1}{T}}}$, we write $F(z)$ rather than $F(z^{\sfrac{1}{T}})$. 
  This convention suppresses two pieces of information which have to be made explicit:
  \begin{enumerate}
    \item 
        The actual formal coordinate is $z^{\sfrac{1}{T}}$, and it has $T$ branches. 
        The deck transformation
        \[
          \nota{\varsigma_T}\colon
          z^{\sfrac{1}{T}}\longmapsto
          \zeta_Tz^{\sfrac{1}{T}}
        \]
        fixes $z$ while permutes the branch. 
        We write $F(\varsigma_T(z))$ for the series obtained by applying this coordinate change.
    \item 
        A further substitution $z\mapsto z^{\sfrac{1}{N}}$ changes the denominator of the exponents. 
        We insist that $F(z^{\sfrac{1}{N}})$ means the series in $U\dbrack{z^{\pm\sfrac{1}{NT}}}$ obtained by making this substitution to $F(z) \in U\dbrack{z^{\pm\sfrac{1}{T}}}$, not a series in $U\dbrack{z^{\pm\sfrac{1}{N}}}$.
  \end{enumerate}
\end{convention}

\begin{recollection}[Geometric picture]
  Let $\pi_T \colon \mathbb{P}^{\sfrac{1}{T}} \rightarrow \mathbb{P}^1$ be the degree-$T$ cyclic cover, ramified at $0$ and $\infty$.  
  For $n\ge0$, consider the configuration space
  \[
    \nota{\mathfrak{X}_{2^{\sfrac{1}{T}}+n}}
    :=\Set*{
      (\pt_0,\pt_1,\ldots,\pt_n,\pt_\infty)
      \in(\mathbb{P}^{\sfrac{1}{T}})^{n+2}
      \given
      \pt_0=0,\ \pt_\infty=\infty,\
      \pi_T(\pt_i)\neq\pi_T(\pt_j)\text{ for }i\neq j
    }.
  \]
  For $n=2$, use $w^{\sfrac{1}{T}}$ and $z^{\sfrac{1}{T}}$ as the coordinates of $\pt_1$ and $\pt_2$, respectively.  
  Then, the function ring of $\mathfrak{X}_{2^{\sfrac{1}{T}}+2}$ is
  \[
    \nota{\mathcal{R}_T}
    :=\k[z^{\pm\sfrac{1}{T}},w^{\pm\sfrac{1}{T}},(z-w)^{-1}].
  \]
  We regard $\mathfrak{X}_{2^{\sfrac{1}{T}}+2}$ as a relative affine scheme over $\mathfrak{X}_{2^{\sfrac{1}{T}}+1}=\Spec\k[w^{\pm\sfrac{1}{T}}]$. 
  Its relative compactification and the corresponding orbifold are
  \[
    \nota{\overline{\mathfrak{X}}_{2^{\sfrac{1}{T}}+2}^{\mathrm{rel}}}
    :=\mathbb{P}^{\sfrac{1}{T}}_{\mathfrak{X}_{2^{\sfrac{1}{T}}+1}}
    \qquad\text{and}\qquad
    \nota{\overline{\mathcal{X}}_{2^{\sfrac{1}{T}}+2}^{\mathrm{rel}}}
    :=\ab[
      \overline{\mathfrak{X}}_{2^{\sfrac{1}{T}}+2}^{\mathrm{rel}}
      /\mu_T^{(z)}
    ],
  \]
  where the relative deck group $\mu_T^{(z)}$ acts on $z^{\sfrac{1}{T}}$ and fixes $w^{\sfrac{1}{T}}$.  
  The relative boundary consists of
  \[
    z^{\sfrac{1}{T}}=0,\qquad
    z^{\sfrac{1}{T}}=\infty,\qquad\text{and}\qquad
    z^{\sfrac{1}{T}}=\zeta_T^{-r}w^{\sfrac{1}{T}}
    \quad(0\le r<T).
  \]
  Their images in $\overline{\mathcal{X}}_{2^{\sfrac{1}{T}}+2}^{\mathrm{rel}}$ are $z=0$, $z=\infty$, and $z=w$ respectively.
\end{recollection}

\begin{notation}\label{nota:relative-boundary-g}
  Corresponding to the boundary components, we have the ordered expansions
  \[
  \begin{aligned}
    \nota{\iota_{w,z}}&\colon
      \mathcal{R}_T\longrightarrow
      \k\dparen{w^{\sfrac{1}{T}}}\dparen{z^{\sfrac{1}{T}}},\\
    \nota{\iota_{z,w}}&\colon
      \mathcal{R}_T\longrightarrow
      \k\dparen{z^{\sfrac{1}{T}}}\dparen{w^{\sfrac{1}{T}}},\\
    \nota{\iota_{w,z-w;r}}&\colon
      \mathcal{R}_T\longrightarrow
      \k\dparen{w^{\sfrac{1}{T}}}\dparen{z-w}
      \qquad(0\le r<T).
  \end{aligned}
  \]
  The first two keep fractional powers of $z$ and $w$ and expand $(z-w)^\Box$, for $\Box\in\Z$, as
  \[
    \sum_{k\ge0}\binom{\Box}{k}(-w)^{\Box-k}z^k
    \qquad\text{and}\qquad
    \sum_{k\ge0}\binom{\Box}{k}z^{\Box-k}(-w)^k
  \]
  respectively. 
  At the $r$-th diagonal branch, powers of $w$ and $z-w$ are kept, while for $\Box \in \sfrac{1}{T}\Z$,
  \[
    \iota_{w,z-w;r}(z^\Box)
    =\zeta_T^{-rT\Box}\sum_{k\ge0}
      \binom{\Box}{k}w^{\Box-k}(z-w)^k.
  \]
  Equivalently,
  \[
    \iota_{w,z-w;r}\phi(z,w)
    =\iota_{w,z-w;0}\phi(\varsigma_T^{-r}(z),w)
    \qquad(\phi\in\mathcal{R}_T).
  \]
\end{notation}

With the branch formulas above, \zcref{lem:diagonal-expansion} also holds for fractional powers, with the target $\mathcal{D}_U(w^{\sfrac{1}{T}},z-w)$ instead.

\begin{definition}\label{def:boundary-residue-g}
  Let $U$ be as in \zcref{def:boundary-residue} and fix $\phi(z,w)\in\mathcal{R}_T$. For
  \[
    F_0\in\mathcal{D}_U(w^{\sfrac{1}{T}},z^{\sfrac{1}{T}}),\qquad
    F_\infty\in\mathcal{D}_U(z^{\sfrac{1}{T}},w^{\sfrac{1}{T}}),\qquad
    F_\Delta=(F_{\Delta,r})_{0\le r<T}
    \in\prod_{r=0}^{T-1}
      \mathcal{D}_U(w^{\sfrac{1}{T}},z-w),
  \]
  we first abbreviate the residue on the $r$-th branch of the cyclic cover by
  \[
    \nota{\oint_{z=w}^{[r]}\phi F_{\Delta,r}\d{z}}
    :=\oint_{z^{\sfrac{1}{T}}=\zeta_T^{-r}w^{\sfrac{1}{T}}}
      \phi(z,w)F_{\Delta,r}(z,w)\d{z}
    =\Res_{z-w}\iota_{w,z-w;r}(\phi(z,w))F_{\Delta,r}(z,w).
  \]
  Then, the three boundary residues on $\overline{\mathcal{X}}_{2^{\sfrac{1}{T}}+2}^{\mathrm{rel}}$ are defined by
  \begin{align*}
    \nota{\oint_{z=0}\phi F_0\d{z}}
      &:=\frac{1}{T}\oint_{z^{\sfrac{1}{T}}=0}
          \phi(z,w)F_0(z,w)\d{z} 
        =\Res_z\iota_{w,z}(\phi(z,w))F_0(z,w),\\
    \nota{\oint_{z=\infty}\phi F_\infty\d{z}}
      &:=\frac{1}{T}\oint_{z^{\sfrac{1}{T}}=\infty}
          \phi(z,w)F_\infty(z,w)\d{z}
        =-\Res_z\iota_{z,w}(\phi(z,w))F_\infty(z,w),\\
    \nota{\oint_{z=w}\phi F_\Delta\d{z}}
      &:=\frac{1}{T}\sum_{r=0}^{T-1}
        \oint_{z=w}^{[r]}\phi F_{\Delta,r}\d{z}
        =\Res_{z-w}\frac{1}{T}\sum_{r=0}^{T-1}\iota_{w,z-w;r}(\phi(z,w))F_{\Delta,r}(z,w).
  \end{align*}
  The factor $\tfrac{1}{T}$ at $0$ and $\infty$ accounts for the ramification of the degree-$T$ cover. 
  At a diagonal branch, $z-w$ is a local coordinate, and the cover is unramified.  
  The average over its $T$ lifts gives the root-stack residue along $z=w$.
\end{definition}

\begin{proposition}[Fractional residue theorem]
  \label{prop:fractional-residue-theorem}
  Let $F\in U\otimes_\k\mathcal{R}_T$, and write
  \[
    F_0:=\iota_{w,z}F,\qquad
    F_\infty:=\iota_{z,w}F,\qquad
    F_{\Delta,r}:=\iota_{w,z-w;r}F
    \quad(0\le r<T).
  \]
  Put $F_\Delta:=(F_{\Delta,r})_{0\le r<T}$.
  Then, for every $\phi(z,w)\in\mathcal{R}_T$,
  \[
    \oint_{z=0}\phi F_0\d{z}
    +\oint_{z=\infty}\phi F_\infty\d{z}
    +\oint_{z=w}\phi F_\Delta\d{z}
    =0.
  \]
\end{proposition}
\begin{proof}
  Apply the ordinary residue theorem on the degree-$T$ cover
  $\overline{\mathfrak{X}}_{2^{\sfrac{1}{T}}+2}^{\mathrm{rel}}$.  
\end{proof}

\subsection{Eigenspace decomposition}

\begin{notation}\label{nota:eigenspace-decomposition}
  Let $U$ be a (topological) vector space with a (continuous) automorphism $g$, and fix a \emph{period} $T$ of $g$, i.e. a positive integer satisfying $g^T=\id_U$.
  We can decompose $U$ into eigenspaces of $g$:
  \[
    \nota{U^{[r]}} := \Set*{ u \in U \given g(u) = \zeta_{T}^{r}u }.
  \]
  where $\zeta_{T}$ is a fixed primitive $T$-th root of unity, and $[r]$ stands for the congruence class of $r$ modulo $T$.
  More concretely, for any $u \in U$, let
  \begin{equation}\label{eq:proj-r}
    \nota{u^{(r)}} :=
    \frac{1}{T}\sum_{k=0}^{T-1} \zeta_{T}^{-rk}g^{k}(u)
    \in U^{[r]},
  \end{equation}
  i.e. the projection to the $[r]$-component. 
  Then, we have the following decomposition
  \begin{equation}\label{eq:decomposition-vect}
    u  =  u^{(n_1)} + \cdots + u^{(n_T)},
  \end{equation}
  whenever $n_1,\cdots,n_T$ form a complete set of representatives modulo $T$.

  For a formal power series $F(z) \in U\dbrack{z^{\pm1}}$, we let $g$ act on $F(z)$ by acting on the coefficients.
\end{notation}
\begin{remark}
  Enlarging the period $T$ only adds zero eigenspaces: $U^{[r]}=0$ unless $[r]$ is annihilated by the order of $g$.
\end{remark}

\begin{lemma}\label{lem:equivariant-diagonal-residue}
  Let $U$ be a topological vector space and $g$ a continuous automorphism of it with a period $T$. Set
  \[
    F(t) = \sum_{n\in\Z}F_nt^{-n-1}\in\mathcal{D}_U(t).
  \]
  For $m\in\sfrac{1}{T}\Z$, one has
  \[
    \frac{1}{T}\sum_{r=0}^{T-1}
      \oint_{z=w}^{[r]}z^m(g^rF(z-w))\d{z}
    = \sum_{k\ge0}
      \binom{m}{k}w^{m-k}F_k^{(Tm)}
      \in\mathcal{D}_U(w^{\sfrac{1}{T}}).
  \]
\end{lemma}
\begin{proof}
  By \zcref{nota:relative-boundary-g},
  \[
    \oint_{z=w}^{[r]}z^m(g^rF(z-w))\d{z}
    =\zeta_T^{-rTm}\sum_{k\ge0}
      \binom{m}{k}w^{m-k}g^rF_k.
  \]
  Averaging over $r$ and applying \eqref{eq:proj-r} to each $F_k$ gives the formula.
\end{proof}



\part{Twisted enveloping algebras}
\label{part:twisted-enveloping-algebras}

After recalling basics of vertex operator algebras, twisted modules and almost-canonically seminormed rings, we introduce twisted enveloping algebras and study their properties in our topological framework.

\section{Vertex operator algebras and their twisted modules}
\label{sec:voa-acs-framework}

\begin{definition}\label{def:voa}
  A \concept{vertex operator algebra} (\emph{VOA} for short) is a quadruple $(V, Y(\cdot,z), \vac, \cfv)$, throughout simply denoted by $V$, where
  \begin{itemize}
    \item $\nota{V} = \bigoplus_{k\in\Z} V_k$ is a graded vector space with $\dim{V_k} < \infty$ and $V_k=0$ for all sufficiently negative $k$;
    \item $\nota{Y( \nocolor{\placeholder} , z )}\colon V \to (\End{V})\dbrack{z^{\pm1}}$ is a linear map assigning to each $a \in V$ the \concept{vertex operator}
      \[
        Y(a,z) = \sum_{n\in\Z} \nota{\vo{a}{n}} z^{-n-1},
      \]
      which is a formal power series with coefficients in $\End{V}$;
    \item $\nota{\vac} \in V_0$ and $\nota{\cfv} \in V_2$ are two distinguished vectors, called the \concept{vacuum vector} and the \concept{conformal vector} respectively.
  \end{itemize}
  These data must satisfy the following axioms:
  \begin{enumerate}[label={\textbf{V}\arabic*}]
    \item \concept{\small(Truncation property)}
      for all $a,b \in V$: $\vo{a}{n}b = 0$ if $n$ is sufficiently large;
    \item \concept{\small(Vacuum property)}
      $Y( \vac, z ) = \id_V$;
    \item \concept{\small(Creation property)}
      for all $a \in V$:
      $Y( a, z ){\vac} \in V\dbrack{z}$ and $\lim\limits_{z\to 0}Y( a, z ){\vac} = a$;
    \item \concept{\small(Jacobi identity)}
      for all $a, b \in V$: the following identity of residues\footnote{See \zcref{def:boundary-residue} for the notations.}
      \begin{align*}
        \MoveEqLeft
        \oint_{z=w}
        \phi(z, w) Y(Y(a, z-w) b, w)  \d{z} \\
        &=
        -\oint_{z=\infty}
        \phi(z, w) Y(a, z) Y(b, w) \d{z}
        -\oint_{z=0}
        \phi(z, w) Y(b, w) Y(a, z) \d{z}
      \end{align*}
      holds for every rational function $\phi(z, w) \in \k[z^{\pm1}, w^{\pm1}, (z-w)^{\pm1}]$;
    \item \concept{\small(Virasoro relations)}
      defining $\nota{\VL{n}}$ ($n \in \Z$) as the coefficients of $Y( \cfv, z )$:
      \[
        Y( \cfv, z ) =
        \sum_{n \in \mathbb{Z}}
        \vo{\cfv}{n}z^{-n-1} =
        \sum_{n \in \mathbb{Z}}
        \VL{n}z^{-n-2},
      \]
      then there is a number $\nota{\cch}\in\k$, called the \concept{central charge}, such that
      \[
        [\VL{m},\VL{n}] =
        (m-n)\VL{m+n} +
        \delta_{m+n,0}\frac{(m-1)m(m+1)}{12}\cch,
      \]
      where $m,n \in \Z$;
    \item \concept{\small($\VL{0}$-spectrum)}
      $\VL{0}|_{V_{n}} = n$ whenever $n \in \Z$;
    \item \concept{\small($\VL{-1}$-derivative property)}
      for all $a \in V$: $Y( \VL{-1}a, z ) = \pdv{z} Y( a, z )$.
  \end{enumerate}
\end{definition}
\begin{remark}
  The \emph{truncation property} implies that the vertex operator $Y(a, z)$ lies in the subspace $\mathcal{D}_{\End{V}}(z)$ (cf. \zcref{nota:continuous-distributions}) of $\End{V}\dbrack{z^{\pm1}}$.
\end{remark}
\begin{remark}\label{rem:topological-voa}
  If we replace the vector space $V$ by topological vector space and replace the endomorphism algebra $\End{V}$ by the continuous endomorphism algebra, then we need to also replace the \emph{truncation property} by the following \emph{convergence property}:
  \begin{enumerate}
    \item[\textbf{V}1']
      for all $a,b \in V$: the sequence $\vo{a}{n}b$ converges to zero as $n$ goes to infinity.
  \end{enumerate}
\end{remark}

\begin{definition}
  For any homogeneous $a \in V$, we call the eigenvalue of $a$ under the operator $\VL{0}$ the \concept{weight} of $a$, and denote it by $\nota{\wt a}$.
\end{definition}

\begin{definition}\label{def:voa-automorphism}
  An automorphism of a vertex operator algebra $V$ is a linear transformation $\nota{g} \in \grp{GL}(V)$ preserving $\vac$ and $\cfv$ and satisfying
  \[
    g Y(a, z) g^{-1} = Y( g(a) , z )
    \txforall a \in V.
  \]
\end{definition}

The notion of weak twisted modules is defined as follows.
\begin{definition}\label{def:twisted-module}
  Let $V$ be a vertex operator algebra and let $g$ be an automorphism of $V$ with a period $T$.
  A \concept{weak $g$-twisted $V$-module} is a vector space $M$ equipped with a linear map
  $Y_{M}(\placeholder,z) \colon V \to (\End{M})\dbrack{z^{\pm\sfrac{1}{T}}}$ assigning to each $a \in V$ a \concept{$g$-twisted vertex operator}
  \[
    \nota{Y_{M}(\nocolor{a},z)} = \sum_{n\in\sfrac{1}{T}\Z} \nota{J^{M}_{n}( \nocolor{z^{-n-\VL{0}}a} )},
  \]
  satisfying the following axioms: 
  \begin{enumerate}[label={\textbf{M}\arabic*}]
    \item \concept{\small(Truncation property)}
      for all $a \in V$, $u \in M$: $J^{M}_{n}(a)u = 0$ if $n$ is sufficiently large;
    \item \concept{\small(Vacuum property)}
      $Y_{M}( \vac , z ) = \id_{M}$;
    \item \concept{\small($g$-equivariance property)}
      for all $a\in V$:
      \[
        Y_{M}( g(a) , \varsigma_{T}(z) ) = Y_{M}( a , z ),
      \]
      where $\varsigma_{T}$ is the coordinate change $z^{\sfrac{1}{T}}\mapsto \zeta_{T} z^{\sfrac{1}{T}}$ (cf. \zcref{conv:frac-power});

    \item \concept{\small(Twisted Jacobi identity)}
      for all $a, b \in V$, the following identity of residues
      \begin{align*}
        \MoveEqLeft
        \frac{1}{T}\sum_{r=0}^{T-1}\oint_{z=w}^{[r]}
        \phi(z, w) Y_{M}(Y(g^r.a, z-w) b, w)
        \d{z} \\
        &=
        -\oint_{z=\infty}
        \phi(z, w) Y_{M}(a, z) Y_{M}(b, w)
        \d{z}
        -\oint_{z=0}
        \phi(z, w) Y_{M}(b, w) Y_{M}(a, z)
        \d{z}
      \end{align*}
      holds for every rational function $\phi(z,w)\in\mathcal{R}_T$.
  \end{enumerate}
  When $g$ is the identity automorphism, we call $M$ a \concept{weak (untwisted) $V$-module}.
\end{definition}

We write $\nota{\cat{WMod}_g(V)}$ for the category of weak $g$-twisted $V$-modules.

\begin{remark}\label{rem:weak-module-distributions}
  With the compact-open topology in \zcref{eg:continuous-hom}, $\End(M)$ is a complete topological associative algebra, and the truncation property is equivalent to
  \[
    Y_M(a,z)\in\mathcal D_{\End(M)}(z^{\sfrac{1}{T}})
    \qquad(a\in V).
  \]
\end{remark}

\begin{definition}
  For any $a\in V$ and any $n\in\sfrac{1}{T}\Z$, the linear transformation $J^{M}_{n}(a) \in \End{M}$ is called a \concept{mode} of \concept{degree} \nota{$-n$} of the vertex operator $Y_{M}$.
  For simplicity, we often omit the superscript $M$ in $J^{M}_{n}$ when there is no confusion.
\end{definition}

The categories of modules that we are interested in are the following.

\begin{definition}
  \label{def:ex-module}
  A weak $g$-twisted $V$-module is called an \concept{exhaustive $g$-twisted $V$-module} if it satisfies the following \emph{weak truncation property}:
  \begin{quote}
    for any $u \in M$, we require $J_{n_1}(a^1) \cdots J_{n_k}(a^k)u = 0$ whenever $n_1+\cdots+n_k$ is sufficiently large.
  \end{quote}
\end{definition}
We write $\nota{\cat{ExMod}_g(V)}$ for the category of exhaustive $g$-twisted $V$-modules.

\begin{definition}\label{def:admissible}
  A weak $g$-twisted $V$-module $M$ is called an \concept{admissible\footnotemark $g$-twisted $V$-module}%
  \footnotetext{The wording of \cite[Definition 3.3]{DLM-TR} alone does not specify whether the grading is part of the data. Here we treat the category of admissible modules as a full subcategory of $\cat{WMod}_g(V)$. The disadvantage is that the category of admissible modules is not abelian in general. However, under $g$-rationality, it coincides with the semisimple Grothendieck category $\cat{ExMod}_g(V)$.}
  if there exists an $\sfrac{1}{T}\N$-grading $M=\bigoplus_{n\in \sfrac{1}{T}\N}M(n)$ such that\footnote{Here and what follows, $M(n)$ is set to be $0$ if $n<0$.}
  \begin{equation*}
    J_{m}(a)M(n)\subset M(n-m)
  \end{equation*}
  for any $a\in V$, $m\in \sfrac{1}{T}\Z$, and any $n\in \sfrac{1}{T}\N$.
\end{definition}

Any admissible $g$-twisted $V$-module $M$ must be exhaustive. In fact, any $u \in M(n)$ can be annihilated by modes of total degree greater than $n$.

\begin{definition}\label{def:ordinary}
  A weak $g$-twisted $V$-module $M$ is called an \concept{ordinary $g$-twisted $V$-module} if it admits a decomposition $M=\bigoplus_{\lambda\in\k}M_{\lambda}$ into eigenspaces of $\VL{0}$ such that $M_{\lambda + n} = 0$ for sufficiently negative $n$ and each eigenspace $M_{\lambda}$ is finite-dimensional.

  A weak $g$-twisted $V$-module $M$ is called \concept{geometrically ordinary} if it becomes ordinary after a base change of $\k$.
\end{definition}

Any ordinary $g$-twisted $V$-module is geometrically ordinary, and any geometrically ordinary $g$-twisted $V$-module is admissible. Conversely, a simple $\sfrac{1}{T}\N$-graded $g$-twisted $V$-module is geometrically ordinary if and only if its homogeneous components are finite-dimensional.

\begin{definition}\label{def:g-rationality}
  A vertex operator algebra $V$ is called \concept{$g$-rational} if every admissible $g$-twisted $V$-module is a direct sum of geometrically ordinary simple $g$-twisted modules.
\end{definition}
\begin{remark}
  This definition follows the spirit of \cite{Z96}.
  When $\k =$ the field of complex numbers, it is equivalent to Dong--Li-Mason's definition by \cite[Theorem 8.1(c)]{DLM-TR}.
\end{remark}

\section{Almost-canonically seminormed rings and exhaustive modules}
\label{sec:acsr-framework}

The concept of almost-canonically seminormed rings was introduced in the appendix of \cite{DGK23}, tracing back to \cite{MNT10}, as a unified framework for various versions of enveloping algebras of vertex operator algebras. Mode transition algebras and the strong identity condition were also introduced there and further developed in \cite{DGK24,Cai26,BFHY26,SIC}.

\begin{definition}\label{def:almost-canonical-seminormed-ring}
  An \emph{almost-canonical (left) seminorm} on a graded algebra $U = \bigoplus_{n \in \sfrac{1}{T}\Z} U_{n}$ is a system of neighborhoods $\nota{\NL{U}}$ of $0$ in $U$ verifying the following conditions:
  \begin{enumerate}[label=(\texttt{A}\arabic*)]
    \item\label{A1} Each neighborhood $\NL[n]U$ is a graded subspace of $U$.
    \item\label{A2} The \emph{canonical seminorm} $\nota{\cNL[\bullet]U}:=UU_{\le-\bullet}$ is contained in and degreewise dense in $\NL[\bullet]U$.
    \item\label{A3} $(\NL[n]U_{p})U_{q}\subset\NL[n-q]U_{p+q}$ and $U_{p}(\NL[n]U_{q})\subset\NL[n]U_{p+q}$ for all $n,p,q\in\sfrac{1}{T}\Z$.
  \end{enumerate}
  Note that, by \zcref[noname]{A3}, each $\NL[n]U$ is a left graded ideal of $U$. Hence, we also call such a system $\NL{U}$ a system of \emph{left neighborhoods}.
\end{definition}

\begin{remark}
  \cite{SIC} uses $\Z$-grading. Multiplying all degrees and neighborhood indices by $T$ converts our $\sfrac{1}{T}\Z$-grading to that convention, without changing the homogeneous subspaces, algebra or topology.
\end{remark}

\begin{remark}
  If $U$ itself is already a graded topological algebra, then we require that each $\NL[n]U$ is a graded closed subspace of $U$.
\end{remark}

\begin{notation}
  The system of \emph{right neighborhoods} associated to $\NL$ is defined, for $n,p\in\sfrac{1}{T}\Z$, by \[\nota{\NR[n]U_p}:=\NL[n-p]U_p.\]
  From this definition, one can verify that $\NR$ satisfies properties analogous to those of $\NL$ with a changing of sides.
  For instance, each $\NR[n]U$ is a graded right ideal of $U$.

  Note that, the left and right neighborhoods $\NL$ and $\NR$ coincide on the subalgebra $U_{0}$; that is, $\NL[\bullet]U_{0} = \NR[\bullet]U_{0}$. So we omit the subscripts and simply denote them by $\nota{\NLR[\bullet]U_{0}}$.
\end{notation}

\begin{notation}\label{nota:grading-dilation}
  Let $U$ be a $\sfrac{1}{T}\Z$-graded almost-canonically seminormed ring and $N$ a positive integer.
  We write $\nota{\mathfrak{d}_N^{\ast}}U$ for the pullback along the dilation $\mathfrak{d}_N\colon\sfrac1{NT}\Z\longrightarrow\sfrac{1}{T}\Z$, $p\longmapsto Np$.
  It has the same algebra and topology as $U$, with degrees and neighborhood indices divided by $N$:
  \[
    (\mathfrak{d}_N^{\ast}U)_p=U_{Np},\qquad
    \NL[n]\ab(\mathfrak{d}_N^{\ast}U)=\NL[Nn]U,\qquad
    \NR[n]\ab(\mathfrak{d}_N^{\ast}U)=\NR[Nn]U
    \quad\ab(p,n\in\sfrac1{NT}\Z).
  \]
  Do not confuse it with enlarging the index lattice from $\sfrac{1}{T}\Z$ to $\sfrac1{NT}\Z$, which keeps the existing degrees and assigns zero homogeneous components to the added indices.
\end{notation}

\begin{definition}
  For $n\in\sfrac{1}{T}\N$ and $p\in\sfrac{1}{T}\Z$, define
  \begin{equation}\label{eq:level-quotients}
    \nota{\mathfrak{L}^{n}_{p}}:=(U/\NL[n+\sfrac{1}{T}]U)_{p},\quad
    \nota{\mathfrak{R}^{n}_{p}}:=(U/\NR[n+\sfrac{1}{T}]U)_{p},\txand
    \nota{\mathsf{A}_{n}}:= \mathfrak{L}^{n}_{0} = \mathfrak{R}^{n}_{0} = U_0/\NLR[n+\sfrac{1}{T}]U_0,
  \end{equation}
  The algebra $\mathsf{A}_{n}$ is called the \concept{level-$n$ algebra} of $U$.
  Note that $\mathfrak{L}^{n}$ (resp. $\mathfrak{R}^{n}$) is a discrete, graded $(U,\mathsf{A}_{n})$-bimodule (resp. $(\mathsf{A}_{n},U)$-bimodule).

  We also introduce the following auxiliary notations:
  \[
    \nota{\mathfrak{L}^{\gr n}_{p}} := \Ker(\mathfrak{L}^{n}_{p} \to \mathfrak{L}^{n-\sfrac{1}{T}}_{p})
    \txand
    \nota{\mathfrak{R}^{\gr n}_{p}} := \Ker(\mathfrak{R}^{n}_{p} \to \mathfrak{R}^{n-\sfrac{1}{T}}_{p}).
  \]
  At the bottom step we put $\mathfrak{L}^{-\sfrac{1}{T}}=\mathfrak{R}^{-\sfrac{1}{T}}=0$, consistently with $\NL[0]U=\NR[0]U=U$.
  Note that they are precisely the associated graded pieces of the filtrations $\NL[\bullet]U_{p}$ and $\NR[\bullet]U_{p}$.
\end{definition}

\begin{remark}\label{rem:zhu-comparison}
  It is shown in \cite{Han20,Han25} that, for the $g$-twisted enveloping algebra $\mathscr{U}_{g}(V)$, the algebras $\mathsf{A}_{n}$ are isomorphic to the algebras $A_{g,n}(V)$ in \cite{DLM-TR}, and the $(\mathsf{A}_{m},\mathsf{A}_{n})$-bimodules $\mathfrak{L}^{n}_{m-n}$ (= $\mathfrak{R}^{m}_{m-n}$) are isomorphic to the $(A_{g,m}(V),A_{g,n}(V))$-bimodules $A_{g,m,n}(V)$ in \cite{DJtams}.
\end{remark}

\begin{definition}\label{def:mode-transition-algebra}
  The \concept{mode transition algebra} of $U$ is the bigraded $U$-bimodule
  \[
    \nota{\mathfrak{A}_{\bullet,-\circ}}:=\mathfrak{L}^{0}_{\bullet}\otimes_{\mathsf{A}_{0}}\mathfrak{R}^{0}_{-\circ}.
  \]
  Its multiplication $\nota{\star}$ is induced by the isomorphism
  \begin{align*}
    \mathfrak{R}^{0}\otimes_{U}\mathfrak{L}^{0}&\longrightarrow\mathsf{A}_{0}\\
    [\alpha]^{\mathsf{R}}_{0}\otimes[\beta]^{\mathsf{L}}_{0}&\longmapsto
    \begin{dcases*}
      [\alpha\beta]_{0} & if $\deg(\alpha)+\deg(\beta)=0$,\\
      0 & otherwise.
    \end{dcases*}
  \end{align*}
\end{definition}

\begin{definition}\label{def:sic}
  We say that the \concept{strong identity condition} holds (for $U$) if there are \emph{strong identity elements} $\nota{\one_{n}}\in\mathfrak{A}_{n,-n}$ ($n\in\sfrac{1}{T}\N$) satisfying the identities
  \begin{equation}\label{eq:sic}
    \one_{n}\star\mathfrak{a}=\mathfrak{a}=\mathfrak{a}\star\one_{m}
    \txforall\mathfrak{a}\in\mathfrak{A}_{n,-m}
    \quad(n,m\in\sfrac{1}{T}\N).
  \end{equation}
\end{definition}

\begin{definition}\label{def:modules}
  A left $U$-module $M$ is said to be
  \begin{itemize}
    \item
      \concept{exhaustive} if it is the union of the subspaces
      \[
        \nota{\Omega_{n}(M)} :=
        \Set*{ u \in M \given \NL[n+\sfrac{1}{T}]U \cdot u = 0 };
      \]
    \item
      \concept{gradable} if it is exhaustive and admits a lower-bounded grading $M = \bigoplus_{n \in \sfrac{1}{T}\Z} M(n)$ satisfying
      \[
        U_{p}M(n) \subset M(p+n) \txforall p,n\in\sfrac{1}{T}\Z.
      \]
  \end{itemize}
  For right $U$-modules, we define exhaustiveness and gradability similarly.
\end{definition}

We write $\nota{\cat{Mod}(A)}$ and $\nota{\cat{ExMod}(U)}$ for the categories of left $A$-modules and exhaustive left $U$-modules, respectively. 
It is well-known that $\cat{Mod}(A)$ is a $\k$-linear Grothendieck category. 
By \cite[Proposition 2.3]{SIC}, so is $\cat{ExMod}(U)$.

\begin{remark}\label{rem:rescaling-modules}
  Since $\mathfrak{d}_N^{\ast}U$ has the same underlying algebra and topology as $U$, rescaling gives $\cat{ExMod}(\mathfrak{d}_N^{\ast}U)=\cat{ExMod}(U)$.
  A lower-bounded graded $U$-module becomes a lower-bounded graded
  $\mathfrak{d}_N^{\ast}U$-module after dividing its degrees by $N$.
\end{remark}

The following is a consequence of \cite[Theorem A]{SIC}:
\begin{theorem}
  Suppose the strong identity condition holds for $U$. Then,
  \begin{enumerate}
    \item $\cat{ExMod}(U)$ coincides with the category of gradable $U$-modules;
    \item $\Omega_0\colon\cat{ExMod}(U)\to\cat{Mod}(\mathsf A_0)$ is an equivalence.
  \end{enumerate}
\end{theorem}

It is shown in \cite[Remark 3.4.6]{DGK23} that the enveloping algebra of a rational VOA satisfies the strong identity condition.
More general is ture:

\begin{definition}[{\cite[Definition 9.1]{SIC}}]\label{def:k-rationality}
  An almost-canonically seminormed ring $U$ is called \concept{$\k$-rational} if every lower-bounded graded left $U$-module is a graded direct sum of simple graded modules with finite-dimensional homogeneous components.
\end{definition}


\begin{theorem}[{\cite[Theorem C]{SIC}}]\label{thm:rational-sic}
  If $U$ is $\k$-rational, then $\cat{ExMod}(U)$ is semisimple and $U$ satisfies the strong identity condition.
\end{theorem}

\begin{remark}\label{rem:rationality-formulations}
  \zcref{def:k-rationality} is equivalent to the following condition: gradable $U$-modules are semisimple and simple gradable $U$-modules have finite-dimensional homogeneous components (in any lower-bounded $U$-module grading).
  A $\k$-rational $U$ has only finitely many simple exhaustive modules up to isomorphism.
\end{remark}

\begin{recollection}[{\cite{LGS17}}]\label{recol:kelly-module-categories}
  For $\k$-linear Grothendieck categories $\mathcal{C}$ and $\mathcal{D}$, their \concept{Deligne--Kelly tensor product} is a $\k$-linear Grothendieck category $\mathcal{C}\boxtimes\mathcal{D}$ equipped with a $\k$-bilinear functor
  \[
    \nota{\boxtimes}\colon\mathcal{C}\times\mathcal{D}
    \longrightarrow\mathcal{C}\boxtimes\mathcal{D}
  \]
  preserving colimits in each variable and satisfying the following universal property: for every $\k$-linear Grothendieck category $\mathcal{E}$, precomposition with $\boxtimes$ induces an equivalence
  \[
    \cat{Cocont}(\mathcal{C}\boxtimes\mathcal{D},\mathcal{E})
    \simeq
    \cat{Cocont}(\mathcal{C},\mathcal{D};\mathcal{E}).
  \]
  Here $\cat{Cocont}$ denotes categories of functors that are $\k$-linear and preserve colimits separately in each input variable.

  Let $A^{\tt I}$ and $A^{\tt II}$ be $\k$-algebras, and let
  \[
    \mathcal{C}^{\tt I}\subset\cat{Mod}(A^{\tt I}),
    \qquad
    \mathcal{C}^{\tt II}\subset\cat{Mod}(A^{\tt II})
  \]
  be full subcategories closed under submodules, quotients and arbitrary direct sums.
  Then both are Grothendieck categories, and tensoring modules identifies their Kelly product with the full subcategory of $\cat{Mod}(A^{\tt I}\otimes A^{\tt II})$ whose restrictions belong to $\mathcal{C}^{\tt I}$ and $\mathcal{C}^{\tt II}$, respectively.
\end{recollection}

\begin{recollection}[{\cite[\S 10]{SIC}}]\label{recol:acsr-tensor-products}
  Let $U^{\tt I}$ and $U^{\tt II}$ be almost-canonically seminormed rings.
  With the total grading, $U^{\tt I}\otimes U^{\tt II}$ is an almost-canonically seminormed ring whose neighborhoods are
  \[
    \NL[n](U^{\tt I} \otimes U^{\tt II})_{p} = \sum_{\substack{i+j=n\\a+b=p}}\NL[i]U^{\tt I}_{a}\otimes \NL[j]U^{\tt II}_{b}\txand
    \NR[n](U^{\tt I} \otimes U^{\tt II})_{p} = \sum_{\substack{i+j=n\\a+b=p}}\NR[i]U^{\tt I}_{a}\otimes \NR[j]U^{\tt II}_{b}.
  \]
  We write $\nota{U^{\mathtt{I}} \ctensor U^{\mathtt{II}}}$ for its degreewise completion. 
  Note that, seeing \cite[Remark 10.1]{SIC}, degreewise completion preserves the discrete quotients and hence the mode transition algebra, the strong identity condition, and any application there.

  Tensoring modules induces an equivalence
  \[
    \cat{ExMod}(U^{\tt I})\boxtimes\cat{ExMod}(U^{\tt II})
    \simeq
    \cat{ExMod}(U^{\tt I}\ctensor U^{\tt II}).
  \]
  If both factors satisfy the strong identity condition, then so does $U^{\tt I}\ctensor U^{\tt II}$, with strong identity elements given by
  \[
    \one_n = \sum_{p+q=n}\one^{\tt I}_p\otimes\one^{\tt II}_q.
  \]
  The converse holds if all factors are nonzero. Indeed, the proof of \cite[Theorem 10.10]{SIC} applies in our context without the stated finiteness assumptions: projecting the tensor-product strong identity relations to bidegree $(0,0)$ in one factor and applying a linear functional sending its level-zero identity to $1$ gives strong identities for the other factor.
\end{recollection}

\begin{remark}\label{rem:rescaled-tensor-products}
  For a positive integer $N$, uniform rescaling preserves the
  underlying algebra and topology of the completed tensor product:
  \[
    \mathfrak{d}_N^{\ast}\ab(U^{\mathtt{I}}\ctensor U^{\mathtt{II}})
    \cong\mathfrak{d}_N^{\ast}U^{\mathtt{I}}
      \ctensor\mathfrak{d}_N^{\ast}U^{\mathtt{II}}.
  \]
  Rescaling the factors by different amounts, saying $N_1$ and $N_2$, may change the completed algebra and its topology.
  Nevertheless, \zcref{rem:rescaling-modules}
  identifies the exhaustive module categories of the rescaled factors
  with those of the original factors.
  Thus \zcref{recol:acsr-tensor-products} gives
  \[
    \cat{ExMod}(\mathfrak{d}_{N_1}^{\ast}U^{\mathtt{I}}
      \ctensor\mathfrak{d}_{N_2}^{\ast}U^{\mathtt{II}})
    \simeq
      \cat{ExMod}(U^{\mathtt{I}})
      \boxtimes\cat{ExMod}(U^{\mathtt{II}})
    \simeq
      \cat{ExMod}(U^{\mathtt{I}}\ctensor U^{\mathtt{II}}).
  \]
  However, such an equivalence does not preserve the gradings of the modules.
\end{remark}

\begin{lemma}\label{lem:tensor-k-rationality}
  Let $U^{\tt I}$ and $U^{\tt II}$ be nonzero almost-canonically seminormed rings, and put $U=U^{\tt I}\ctensor U^{\tt II}$ with the total grading.
  Then $U$ is $\k$-rational if and only if both $U^{\tt I}$ and $U^{\tt II}$ are $\k$-rational.
\end{lemma}
\begin{proof}
  Suppose both factors are $\k$-rational.
  By \cite[Theorem C]{SIC}, they satisfy the strong identity condition, and their canonical modules $\mathfrak{L}^{0}(U^{\tt I})$ and $\mathfrak{L}^{0}(U^{\tt II})$ have finite-dimensional homogeneous components.
  In particular, their level-zero algebras $\mathsf{A}_0(U^{\tt I})$ and $\mathsf{A}_0(U^{\tt II})$ are finite-dimensional semisimple algebras.
  By \zcref{recol:acsr-tensor-products}, $U$ satisfies the strong identity condition and $\mathsf{A}_0(U)\cong\mathsf{A}_0(U^{\tt I})\otimes\mathsf{A}_0(U^{\tt II})$, which is also semisimple since $\k$ is of characteristic zero.
  Moreover, the decomposition
  \begin{equation}
    \label{eq:weighted-tensor-module-grading}
    \tag{\textasteriskcentered}
    \mathfrak{L}^{0}_{s}(U)\cong
    \bigoplus_{p+q=s}
      \mathfrak{L}^{0}_{p}(U^{\tt I})\otimes\mathfrak{L}^{0}_{q}(U^{\tt II})
  \end{equation}
  shows that $\mathfrak{L}^{0}(U)$ has finite-dimensional components.
  The reconstruction equivalence in \cite[Theorem A]{SIC}, given by $\mathfrak{L}^{0}(U)\otimes_{\mathsf{A}_0(U)}- \colon \cat{Mod}(\mathsf{A}_0(U)) \to \cat{ExMod}(U)$, now gives $\k$-rationality following \zcref{rem:rationality-formulations}.

  Conversely, suppose $U$ is $\k$-rational.
  Its canonical modules $\mathfrak{L}^{0}(U)$ and $\mathfrak{R}^{0}(U)$ have finite-dimensional homogeneous components, so the finiteness hypotheses of \cite[Theorem D]{SIC} hold.
  Since $U$ satisfies the strong identity condition, that theorem gives the same condition for both factors.
  The tensor formula $\mathsf{A}_0(U) \cong \mathsf{A}_0(U^{\tt I})\otimes\mathsf{A}_0(U^{\tt II})$ shows that $\mathsf{A}_0(U^{\tt I})$ and $\mathsf{A}_0(U^{\tt II})$ are nonzero finite-dimensional semisimple algebras just as $\mathsf{A}_0(U)$.
  On the other hand, \eqref{eq:weighted-tensor-module-grading} shows that $\mathfrak{L}^{0}(U^{\tt I})$ and $\mathfrak{L}^{0}(U^{\tt II})$ have finite-dimensional homogeneous components.
  Then the same reconstruction argument proves $\k$-rationality of both factors.
\end{proof}

\section{Twisted enveloping algebras}
\label{sec:voa-enveloping-algebras}

In this section, $V$ is a vertex operator algebra and $g$ is an automorphism of $V$ with a fixed period $T$. We construct the $g$-twisted enveloping algebra $\mathscr{U}_{g}(V)$ and study its properties.

\subsection{Definition of the twisted enveloping algebra}

Starting from the free unital associative algebra on the modes $J_n(a)$ ($a\in V$, $n\in\sfrac{1}{T}\Z$), we take the degreewise completion with respect to the canonical seminorm of \zcref{def:almost-canonical-seminormed-ring} and then divide by the closure of the ideal generated by the vacuum, $g$-equivariance, and twisted Jacobi relations. These relations are the mode-level analogues of the axioms of a weak $g$-twisted $V$-module (cf. \zcref{def:twisted-module}), with truncation replaced by convergence of the mode coefficients $J_n(a)$ to zero as $n\to\infty$ in the canonical topology. 

\begin{construction}
  Let $\nota{\mathbb{U}_{\sfrac{1}{T}}(V)}$ be the linear space spanned by the (possibly empty) words
  \[
    J_{n_1}(a^1) \cdots J_{n_k}(a^k),
  \]
  where $a^1,\cdots,a^k\in V$ and $n_1,\cdots,n_k\in\sfrac{1}{T}\Z$. We further assume that $J_{n_1}(a^1) \cdots J_{n_k}(a^k)$ is linear in each $a^i$ for $1 \le i \le k$.
  It is a unital associative algebra with the multiplication given by \emph{concatenation} and the identity $1_{\mathbb{U}_{\sfrac{1}{T}}(V)}$ given by the empty word.
  This algebra is $\sfrac{1}{T}\Z$-graded with the grading given by
  \[
    \deg J_{n_1}(a^1)\cdots J_{n_k}(a^k) = -n_1-\cdots-n_k.
  \]

  Let $\widehat{\mathbb{U}}_{\sfrac{1}{T}}(V)_{\bullet}$ be the degreewise completions of $\mathbb{U}_{\sfrac{1}{T}}(V)_{\bullet}$ with respect to its \emph{canonical seminorm} (cf. \zcref{def:almost-canonical-seminormed-ring}), that is:
  \[
    \nota{\widehat{\mathbb{U}}_{\sfrac{1}{T}}(V)_{\bullet}}
    :=
    \varprojlim_{n}
    \frac{\mathbb{U}_{\sfrac{1}{T}}(V)_{\bullet}}{\NL[n]\mathbb{U}_{\sfrac{1}{T}}(V)_{\bullet}}
    =
    \varprojlim_{n'}
    \frac{\mathbb{U}_{\sfrac{1}{T}}(V)_{\bullet}}{\NR[n']\mathbb{U}_{\sfrac{1}{T}}(V)_{\bullet}}.
  \]
  Elements in $\widehat{\mathbb{U}}_{\sfrac{1}{T}}(V)$ are called ($\sfrac{1}{T}$-fractional) \concept{modes} and among them, $J_{n}(a)$ ($a \in V, n \in \sfrac{1}{T}\Z$) are called \concept{simple modes}.
\end{construction}

\begin{definition}
  \label{def:enveloping-algebra}
  The \concept{$g$-twisted enveloping algebra} $\nota{\mathscr{U}_{g}(V)}$ is the quotient graded algebra of $\widehat{\mathbb{U}}_{\sfrac{1}{T}}(V)$ by the following homogeneous relations:
  \begin{description}[labelindent=\parindent]
    \item[Vacuum property]
      \[
        J_{n}(\vac) \sim \delta_{n,0}1_{\widehat{\mathbb{U}}_{\sfrac{1}{T}}(V)}
        \qquad (n \in \sfrac{1}{T}\Z);
      \]
    \item[$g$-equivariance property]
      \[
        \zeta_{T}^{-nT}J_{n}(g(a)) \sim J_{n}(a)
        \qquad (a \in V, n \in \sfrac{1}{T}\Z);
      \]
    \item[Twisted Jacobi identity (component form)]\footnote{Recall \zcref{nota:eigenspace-decomposition} for the notation $a^{(mT)}$.}
      \begin{align*}
        \MoveEqLeft
        \sum_{i\ge 0}\binom{\wt a+m-1}{i}
        J_{m+n+k}( \vo{a^{(mT)}}{k+i}b )
        \\
        &\sim
        \sum_{i\ge 0}(-1)^{i}\binom{k}{i}(
          J_{m+k-i}(a)J_{n+i}(b) -
          (-1)^{k}J_{n+k-i}(b)J_{m+i}(a)
        ),
      \end{align*}
      where $a,b\in V$, $m,n\in\sfrac{1}{T}\Z$, and $k\in\Z$.
  \end{description}
  For convenience, we write $\nota{\widehat{\mathbb{I}}_{g}}$ for the closure of the homogeneous two-sided ideal of $\widehat{\mathbb{U}}_{\sfrac{1}{T}}(V)$ generated by the above relations. Thus,
  \[
    \mathscr{U}_{g}(V)=\widehat{\mathbb{U}}_{\sfrac{1}{T}}(V)/\widehat{\mathbb{I}}_{g}.
  \]
  We give this quotient its degreewise quotient topology.
  The closures of the canonical left neighborhoods, and their images in the quotient, form bases of graded open left ideals. Hence $\widehat{\mathbb{U}}_{\sfrac{1}{T}}(V)$ and $\mathscr{U}_{g}(V)$ are graded topological associative algebras in the sense of \zcref{rec:topological-associative-algebras}.
  For simplicity, when $g$ is the identity automorphism, we omit the subscripts $\sfrac{1}{1}$ and $\id$.
\end{definition}
\begin{remark}\label{rem:exhaustive-module-correspondence}
  The mode action gives an equivalence
  \[
    \cat{ExMod}_g(V)
    \simeq\cat{ExMod}(\mathscr{U}_g(V)).
  \]
  It also identifies compatible lower-bounded gradings and their homogeneous components, so admissible $g$-twisted $V$-modules correspond to gradable $\mathscr{U}_g(V)$-modules.
  By \zcref{rem:rationality-formulations}, $V$ is $g$-rational if and only if $\mathscr{U}_g(V)$ is $\k$-rational.
\end{remark}
\begin{remark}\label{rem:twisted-jacobi-eigenstates}
  The $g$-equivariance property implies that $J_{n}(a) = 0$ whenever $a \in V^{[r]}$ and $n\not\in \tfrac{r}{T}+\Z$. Under this observation, the twisted Jacobi identity can be rewritten as
  \begin{equation}
    \begin{aligned}
      \MoveEqLeft
      \sum_{i\ge 0}\binom{\wt a+m-1}{i}
      J_{m+n+k}( \vo{a}{k+i}b )
      \\
      &\sim
      \sum_{i\ge 0}(-1)^{i}\binom{k}{i}(
        J_{m+k-i}(a)J_{n+i}(b) -
        (-1)^{k}J_{n+k-i}(b)J_{m+i}(a)
      )
    \end{aligned}
  \end{equation}
  for all $a\in V^{[r]}, b\in V$, $m\in \tfrac{r}{T}+\Z, n\in\sfrac{1}{T}\Z$, and $k\in\Z$.
\end{remark}
\begin{remark}\label{rem:enveloping-choice-of-T}
  The graded algebra $\mathscr{U}_g(V)$, with its canonical topology, is independent of the chosen period. Indeed, under period enlargement with compatible choices of roots of unity, the additional modes vanish by \zcref{rem:twisted-jacobi-eigenstates}, and the existing modes retain their degrees.
\end{remark}

\begin{convention}\label{conv:formal-vertex-operator}
  It is convenient to introduce the following \emph{formal vertex operator}
  \[
    a \in V \longmapsto
    \nota{Y_{g}[\nocolor{a},z]} :=
    \sum_{n\in\sfrac{1}{T}\Z}
    J_{n}( z^{-n-\VL{0}}a )
    \in \widehat{\mathbb{U}}_{\sfrac{1}{T}}(V)\dbrack{z^{\pm \sfrac{1}{T}}}.
  \]
  With boundary residues interpreted coefficientwise as in \zcref{rem:graded-boundary-residues}, the above homogeneous relations take the following form:
  \begin{description}[labelindent=\parindent]
    \item[Vacuum property]
      \[
        Y_{g}[\vac,z] \sim 1_{\widehat{\mathbb{U}}_{\sfrac{1}{T}}(V)};
      \]
    \item[$g$-equivariance property]
      \[
        Y_{g}[ g(a) , \varsigma_{T}(z) ] \sim Y_{g}[ a , z ]
        \txforall a \in V;
      \]
    \item[Twisted Jacobi identity (operator form)]
      \begin{align*}
        \MoveEqLeft
        \frac{1}{T}\sum_{r=0}^{T-1}\oint_{z=w}^{[r]}
        \phi(z, w) Y_{g}[Y(g^r.a, z-w) b, w]
        \d{z} \\
        &\sim
        -\oint_{z=\infty}
        \phi(z, w) Y_{g}[a, z] Y_{g}[b, w]
        \d{z}
        -\oint_{z=0}
        \phi(z, w) Y_{g}[b, w] Y_{g}[a, z]
        \d{z}
      \end{align*}
      where $a,b\in V$ and $\phi(z,w)\in\mathcal{R}_T$.
  \end{description}
  In other words, viewing each $J_{n}(a)$ as a linear transformation on $\widehat{\mathbb{U}}_{\sfrac{1}{T}}(V)$, the above relations say that $\mathscr{U}_{g}(V)$ behaves like a weak $g$-twisted $V$-module with vertex operator $Y_{g}[\placeholder,z]$, except that it does not satisfy the truncation property.
\end{convention}
\begin{remark}
  For any $a\in V$, the canonical seminorm gives the genuine convergence
  \[
    \lim_{n\to\infty}J_n(a)=0
  \]
  in the induced topology on the underlying vector space of $\widehat{\mathbb{U}}_{\sfrac{1}{T}}(V)$.
  Indeed, if $n\ge N$, then $J_n(a)$ has degree $-n\le-N$, and hence
  \[
    J_n(a)\in \widehat{\mathbb{U}}_{\sfrac{1}{T}}(V)
      \widehat{\mathbb{U}}_{\sfrac{1}{T}}(V)_{\le-N}
      \subseteq\NL[N]\widehat{\mathbb{U}}_{\sfrac{1}{T}}(V).
  \]
  If $a$ is homogeneous, then
  \[
    Y_g[a,z]=\sum_{n\in\sfrac{1}{T}\Z}J_n(a)z^{-n-\wt a}.
  \]
  Hence, by \zcref{nota:continuous-distributions} and linearity,
  for every $a\in V$ we have
  \[
    Y_g[a,z]\in
    \mathcal{D}_{\widehat{\mathbb{U}}_{\sfrac{1}{T}}(V)}
      (z^{\sfrac{1}{T}}).
  \]
\end{remark}
\begin{notation}\label{nota:field-relation-ideal}
  Let $V$ be a VOA and $R$ a (graded) topological associative algebra in the sense of \zcref{rec:topological-associative-algebras}. 
  A linear map $F[-,z]\colon V\to\mathcal D_R(z^{\sfrac{1}{T}})$ is called a \concept{field assignment}. 
  For homogeneous $a\in V$, write
  \[
    F[a,z]=\sum_{n\in\sfrac{1}{T}\Z}\nota{F_n(a)}z^{-n-\wt a}.
  \]
  It is called \concept{graded} if $R$ is graded and $F_n(a)\in R_{-n}$.

  For a (graded) field assignment $F$, let $\nota{\mathbb{I}^g_F}$ denote the closed two-sided ideal of $R$ generated by the coefficients of its vacuum, $g$-equivariance and twisted Jacobi relations. This ideal is homogeneous in the graded case.


  When $R$ is graded, all open left ideals below are understood to be graded.
\end{notation}
\begin{remark}
  The two forms of the twisted Jacobi identity are equivalent.  
  Apply \zcref{lem:equivariant-diagonal-residue} to the single-variable distribution
  \[
    t^kY(\placeholder,t)b
    \in\mathcal{D}_{\Hom_\k(V,V)}(t),
  \]
  with $g$ acting by precomposition on the placeholder.  
  VOA truncation gives this membership and makes the following coefficient sum finite.
  Evaluating at $a$ and applying $Y_g[-,w]$ gives
  \begin{align*}
    \MoveEqLeft
    \frac{1}{T}\sum_{r=0}^{T-1}\oint_{z=w}^{[r]}
      z^{m+\wt a-1}(z-w)^kY_g[Y(g^ra,z-w)b,w]\d{z} \\
    &=\sum_{i\ge0}\binom{m+\wt a-1}{i}w^{m+\wt a-1-i}
      Y_g[\vo{a^{(T(m+\wt a-1))}}{k+i}b,w].
  \end{align*}
  Taking $\phi(z,w)=z^{m+\wt a-1}w^{n+\wt b-1}(z-w)^k$ ($m,n\in\sfrac{1}{T}\Z$, $k\in\Z$), the projector above is $a^{(T(m+\wt a-1))}=a^{(mT)}$, and we obtain the homogeneous relation
  \begin{align*}
    \MoveEqLeft
    \sum_{p\in\sfrac{1}{T}\Z}\sum_{i\ge 0}\binom{m+\wt a-1}{i}J_{m+n+k-p-1}(\vo{a^{(mT)}}{k+i}b)w^{p}\\
    &\sim
    \sum_{p\in\sfrac{1}{T}\Z}\sum_{i\ge 0}(-1)^i\binom{k}{i}(J_{m+k-i}(a)J_{n+i-p-1}(b)-(-1)^{k}J_{n+k-i-p-1}(b)J_{m+i}(a))w^{p}.
  \end{align*}
  in $\mathcal{D}_{\widehat{\mathbb{U}}_{\sfrac{1}{T}}(V)}(w^{\sfrac{1}{T}})$. 
  Taking the coefficients of $w^{-1}$ gives the component form.
  Conversely, by varying $n\in\sfrac{1}{T}\Z$, we recover the above homogeneous relation.
  Since $z^{m+\wt a-1}w^{n+\wt b-1}(z-w)^k$ ($m,n\in\sfrac{1}{T}\Z$, $k\in\Z$) span $\mathcal{R}_T$, we get the operator form.
\end{remark}

\subsection{Properties of the twisted enveloping algebra}

We first give a criterion for the twisted Jacobi identity that will apply both to various field assignments used in this paper.

Let $R$ be a (graded) topological associative algebra in the sense of
\zcref{rec:topological-associative-algebras}, and let
\[
  F[-,z]\colon V\longrightarrow\mathcal{D}_R(z^{\sfrac{1}{T}})
\]
be a (graded) field assignment satisfying the \emph{$g$-equivariance property}:
\[
  F[g(-),\varsigma_T(z)]=F[-,z].
\]

\begin{definition}[Locality]\label{def:locality}
  We say that $F$ satisfies \concept{locality}, or is \concept{local}, if, for every $a,b\in V$, there is $K\in\N$ such that
  \[
    (z-w)^K F[a,z]F[b,w]
    =
    (z-w)^K F[b,w]F[a,z]
    \in R\dbrack{z^{\pm\sfrac{1}{T}},w^{\pm\sfrac{1}{T}}}.
  \]
\end{definition}

\begin{lemma}\label{lem:local-product-residues}
  Suppose that $F$ is local. For homogeneous $a,b\in V$, choose $K\in\N$ such that
  \[
    (z-w)^K F[a,z]F[b,w]=(z-w)^K F[b,w]F[a,z].
  \]
  The coefficients of this common distribution tend uniformly to $0$ as either exponent tends to $-\infty$. The diagonal expansions,
  \[
    \iota_{w,z-w;r}\ab(F[a,z]F[b,w])
    \in\mathcal{D}_R(w^{\sfrac{1}{T}},z-w)
    \qquad(0\le r<T)
  \] 
  interpreted as in \zcref{rem:diagonal-poles}, are well-defined and independent of the choice of $K$.
  For every $\phi(z,w)\in\mathcal{R}_T$, we have
  \[
    \begin{aligned}
      &\oint_{z=w}\phi(z,w)\,
        \ab(\iota_{w,z-w;r}\ab(F[a,z]F[b,w]))_{0\le r<T}\d{z}\\
      &\qquad=
        -\oint_{z=\infty}\phi(z,w) F[a,z]F[b,w]\d{z}
        -\oint_{z=0}\phi(z,w) F[b,w]F[a,z]\d{z}
    \end{aligned}
  \]
  coefficientwise in $w^{\sfrac{1}{T}}$.
\end{lemma}
\begin{proof}
  Modulo an open left ideal $W$, the two expressions for the common distribution give uniform lower bounds in $w$ and $z$, respectively.
  For each coefficient in its diagonal expansion, these bounds leave only finitely many contributing terms modulo $W$.
  Thus \zcref{lem:diagonal-expansion} applies.
  Division by $(z-w)^K$ gives the asserted expansions and their independence of $K$ by \zcref{rem:diagonal-poles}.

  By linearity, fix a test monomial $\phi\in\mathcal{R}_T$ and a $w^{\sfrac{1}{T}}$-coefficient of the boundary residues. The contributing terms in the common distribution have fixed total degree in $z,w$, so the lower bounds leave only finitely many modulo $W$ at each boundary.
  The \nameref{prop:fractional-residue-theorem} therefore gives the residue identity modulo $W$. Then separatedness gives the identity in $R$.
\end{proof}

\begin{definition}[Branchwise weak associativity]
  \label{def:branchwise-weak-associativity}
  A local field assignment $F$ satisfies \concept{branchwise weak associativity} if, for every open subspace $W\subset R$, $a\in V^{[q]}$, and $b\in V$, the following coefficientwise identity holds for all sufficiently large $\ell\in\N$:
  \begin{align*}
    z^{\ell+\sfrac{q}{T}}F[a,z]F[b,w]\Big|_{z=(z-w)+w}
    &\equiv
    \iota_{w,z-w;0}(z^{\ell+\sfrac{q}{T}})F[Y(a,z-w)b,w]
    \pmod W.
  \end{align*}
  On the left, substitution expands in nonnegative powers of $w$.
\end{definition}

\begin{proposition}
  \label{prop:jacobi-criterion}
  Suppose that $F$ has the vacuum property. Then its twisted Jacobi identity holds if and only if $F$ is local and satisfies branchwise weak associativity.
\end{proposition}
\begin{proof}
  The twisted Jacobi identity implies locality by testing Laurent polynomials multiplied by a sufficiently large power of $z-w$.

  Now assume locality. By linearity, take $a\in V^{[q]}$ and $b\in V$ homogeneous, and choose $K$ as in \zcref{lem:local-product-residues}.
  We claim that the twisted Jacobi identity is equivalent to the following identity for every fixed $\ell\in\N$:
  \begin{equation}
    \label{eq:jacobi-diagonal-expansion}
    \iota_{w,z-w;0}\ab(z^{\ell+\sfrac{q}{T}}F[a,z]F[b,w]) = 
    \iota_{w,z-w;0}(z^{\ell+\sfrac{q}{T}})F[Y(a,z-w)b,w].
  \end{equation}
  Assume the twisted Jacobi identity. Comparing it with the residue identity of \zcref{lem:local-product-residues}, using the test functions $z^{\ell+\sfrac{q}{T}}(z-w)^sw^n$ with $s\in\Z$ and $n\in\sfrac{1}{T}\Z$, shows that the product's and the iterate's diagonal boundary residues agree.
  In these tests, each diagonal boundary residue equals the branch-zero residue: the product has this property because $z^{\sfrac{q}{T}}F[a,z]$ has only integral powers in $z$; for the iterate, the roots of unity from the test function (which gives $\zeta_T^{-rq}$) and from $g^ra$ (which gives $\zeta_T^{rq}$) cancel with each other, as in \zcref{lem:equivariant-diagonal-residue}.
  Multiplication by $(z-w)^s$ followed by $\Res_{z-w}$ extracts the $(z-w)^{-s-1}$-coefficient. Thus varying $s$ and comparing $w^{\sfrac{1}{T}}$-coefficients gives \eqref{eq:jacobi-diagonal-expansion}.
  Conversely, $\iota_{w,z-w;0}(z^{\sfrac{q}{T}})$ is invertible, with constant term $w^{\sfrac{q}{T}}$.
  Cancelling it in the case $\ell=0$ of \eqref{eq:jacobi-diagonal-expansion} identifies the product's branch-zero diagonal expansion with $F[Y(a,z-w)b,w]$. Equivariance then gives
  \[
    \iota_{w,z-w;r}\ab(F[a,z]F[b,w])
    =F[Y(g^ra,z-w)b,w]
    \qquad(0\le r<T).
  \]
  Substituting these identities into the residue identity of \zcref{lem:local-product-residues} gives the twisted Jacobi identity.

  To compare \eqref{eq:jacobi-diagonal-expansion} with branchwise weak
  associativity, we first show that, for each open subspace $W\subset R$
  and all sufficiently large $\ell\in\N$,
  \[
    \iota_{w,z-w;0}\ab(z^{\ell+\sfrac{q}{T}}F[a,z]F[b,w])
    \equiv
    z^{\ell+\sfrac{q}{T}}F[a,z]F[b,w]\Big|_{z=(z-w)+w}
    \pmod W.
  \]
  Here powers of $z$ are expanded binomially in nonnegative powers of
  $z-w$ and $w$, respectively.
  Fix an open left ideal $W\subset R$. By \zcref{lem:local-product-residues}, the product multiplied by $(z-w)^K$ has lower bounds in both variables modulo $W$. Equivariance makes its $z$-powers integral after multiplication by $z^{\sfrac{q}{T}}$.
  Thus, for any sufficiently large $\ell$, $z^{\ell+\sfrac{q}{T}}(z-w)^K F[a,z]F[b,w]$ has only nonnegative integral $z$-powers modulo $W$. 
  For such $z$-powers, both binomial expansions are the same polynomial in $z-w$ and $w$. 
  The lower bounds leave only finitely many contributing monomials modulo $W$ in each coefficient, so the two expansions agree coefficientwise modulo $W$.
  Dividing by $(z-w)^K$ proves the displayed congruence.

  If \eqref{eq:jacobi-diagonal-expansion} holds, the congruence just proved yields branchwise weak associativity.
  Conversely, assume branchwise weak associativity. For each open left ideal $W$, choosing $\ell$ large enough for weak associativity and the displayed congruence gives \eqref{eq:jacobi-diagonal-expansion} modulo $W$. 
  The choice of $\ell$ may depend on $W$. Cancelling the common invertible scalar factor $\iota_{w,z-w;0}(z^{\ell+\sfrac{q}{T}})$ from both sides of \eqref{eq:jacobi-diagonal-expansion} gives
  \[
    \iota_{w,z-w;0}\ab(F[a,z]F[b,w])
    \equiv F[Y(a,z-w)b,w]\pmod W.
  \]
  Both sides are independent of $\ell$ and $W$. Since they agree coefficientwise modulo every open left ideal, separatedness gives equality in $R$. The first part of the proof then gives the twisted Jacobi identity.
\end{proof}

\begin{corollary}\label{cor:enveloping-identities}
  In $\mathscr{U}_{g}(V)$, we have the following identity:
  \begin{enumerate}
    \item\label{cor:enveloping-identities:derivative} $\VL{-1}$-derivative property
      \[
        Y_{g}[\VL{-1}a, z] = \pdv{z}Y_{g}[a, z]
        \txforall a \in V;
      \]
    \item\label{cor:enveloping-identities:product} Given $a, b\in V$, suppose $\vo{(a^{(r)})}{i}b=0$ for all $r$ and $i\ge 0$, then we have
      \[
        Y_{g}[\vo{a}{-1}b, z] = Y_{g}[a,z]Y_{g}[b,z].
      \]
  \end{enumerate}
\end{corollary}
\begin{proof}
  The first identity follows from the vacuum property and the component form of the twisted Jacobi identity (cf. \zcref{def:enveloping-algebra}) with $b=\vac$ and $k=-2$.

  For the second, the eigenspace decomposition reduces the proof to
  $a\in V^{[q]}$, with $b$ arbitrary. The hypothesis implies that every
  $Y(g^ra,z-w)b$ has no negative $(z-w)$-powers. For Laurent monomial tests
  in $z^{\sfrac{1}{T}},w^{\sfrac{1}{T}}$, all diagonal residues vanish, so the
  twisted Jacobi identity gives
  \[
    Y_g[a,z]Y_g[b,w]=Y_g[b,w]Y_g[a,z].
  \]
  Thus we may take $K=0$ in \zcref{lem:local-product-residues}.
  By \eqref{eq:jacobi-diagonal-expansion},
  \[
    \iota_{w,z-w;0}(Y_g[a,z]Y_g[b,w])
    =Y_g[Y(a,z-w)b,w].
  \]
  The $(z-w)^0$-coefficient on the left is the well-defined product
  $Y_g[a,w]Y_g[b,w]$, while on the right it is
  $Y_g[\vo{a}{-1}b,w]$. This proves the claim.
\end{proof}


The following lemma says that the enveloping algebra is invariant under conjugation of the automorphism.
\begin{lemma}\label{lem:conjugation-invariance}
  Let $h$ be another automorphism of the VOA $V$.
  Then, there exists an isomorphism of graded topological algebras
  \[
    \Upsilon \colon \mathscr{U}_{g}(V) \overset{\cong}{\longrightarrow} \mathscr{U}_{hgh^{-1}}(V)
    \colon
    J_{n}(a) \mapsto J_{n}(h(a)).
  \]
\end{lemma}
\begin{proof}
  The assignment $J_{n}(a) \mapsto J_{n}(h(a))$ preserves degree and is therefore continuous, so it extends to an automorphism of $\widehat{\mathbb{U}}_{\sfrac{1}{T}}(V)$. Since $h$ preserves the vacuum vector, it preserves the vacuum relation.
  For the equivariance property, we have
  \[
    \Upsilon(\zeta_{T}^{-nT}J_{n}(g(a))) = \zeta_{T}^{-nT}J_{n}(hg(a)) = \zeta_{T}^{-nT}J_{n}(hgh^{-1}h(a)) \sim J_{n}(h(a)) = \Upsilon(J_{n}(a)).
  \]
  For the twisted Jacobi identity, note that
  \[
    h(\vo{a}{k}b) = \vo{h(a)}{k}h(b).
  \]
  Moreover, $h$ preserves weights and intertwines the eigenspace projections for $g$ and $hgh^{-1}$. Hence the twisted Jacobi identities correspond under $\Upsilon$.
  Thus $\Upsilon$ descends to a continuous homomorphism $\mathscr{U}_g(V)\to\mathscr{U}_{hgh^{-1}}(V)$; the same construction with $h^{-1}$ gives its continuous inverse.
\end{proof}

\section{Tensor products of twisted enveloping algebras}
\label{sec:enveloping-tensor-products}

Given two VOAs
\[
  (V^{\tt I},Y^{\tt I}(\placeholder,z),\vac^{\tt I},\cfv^{\tt I})
  \txand
  (V^{\tt II},Y^{\tt II}(\placeholder,z),\vac^{\tt II},\cfv^{\tt II}),
\]
there is a VOA structure on the tensor product space $V^{\tt I}\otimes V^{\tt II}$:
\begin{itemize}
  \item The vertex operator $Y(\placeholder,z)$ is given by
    \[
      Y(a^{\tt I}\otimes a^{\tt II},z) := Y^{\tt I}(a^{\tt I},z)\otimes Y^{\tt II}(a^{\tt II},z),
    \]
  \item The vacuum vector is given by
    \[
      \vac := \vac^{\tt I}\otimes \vac^{\tt II},
    \]
  \item The conformal vector is given by
    \[
      \cfv := \cfv^{\tt I}\otimes \vac^{\tt II} + \vac^{\tt I}\otimes \cfv^{\tt II}.
    \]
\end{itemize}

The following lemma extends the tensor-product Jacobi identity argument in \cite{FHL} to the present setting.

\begin{lemma}\label{lem:tensor-fields}
  Let $g$ and $h$ be automorphisms of $V^{\tt I}$ and $V^{\tt II}$, respectively, with a common period $T$.
  Let $R$ be a (graded) topological associative algebra in the sense of \zcref{rec:topological-associative-algebras}.
  Let $F^{\tt I}[-,z]$ and $F^{\tt II}[-,z]$ be (graded) field assignments from $V^{\tt I}$ and $V^{\tt II}$, respectively, to $\mathcal{D}_R(z^{\sfrac{1}{T}})$.
  Suppose that they commute coefficientwise at independent variables:
  \[
    \ab[F^{\tt I}[a_1,z],F^{\tt II}[a_2,w]]=0
    \qquad(a_1\in V^{\tt I},\ a_2\in V^{\tt II}).
  \]
  Then the (graded) field assignment
  \[
    F[-\otimes-,z] := F^{\tt I}[-, z] F^{\tt II}[-, z] = F^{\tt II}[-, z] F^{\tt I}[-, z]
  \]
  is well-defined and the relation ideals satisfy
  \[
    \mathbb{I}^{g\otimes h}_{F}
    \subset
    \overline{\mathbb{I}^{g}_{F^{\tt I}}+\mathbb{I}^{h}_{F^{\tt II}}}.
  \]
  See \zcref{nota:field-relation-ideal} for the notation.
\end{lemma}

\begin{proof}
  We first show that $F$ is well-defined. For homogeneous
  $a_1\in V^{\tt I}$, $a_2\in V^{\tt II}$ and $n\in\sfrac{1}{T}\Z$,
  consider the sums
  \[
    \sum_{r\in\sfrac{1}{T}\Z}F^{\tt I}_{r}(a_1)F^{\tt II}_{n-r}(a_2) = \sum_{r\in\sfrac{1}{T}\Z}F^{\tt II}_{r}(a_2)F^{\tt I}_{n-r}(a_1).
  \]
  Fix an open left ideal $W$. For sufficiently negative $r$, the right factor $F^{\tt II}_{n-r}(a_2)$ lies in $W$; for sufficiently positive $r$, $F^{\tt I}_r(a_1)$ lies in $W$ and commutativity moves it to the right.
  Thus only finitely many summands are nonzero modulo $W$. Completeness gives summability; in the graded case, all summands have the same degree, so degreewise completeness suffices. The same argument shows that the sum lies in $W$ for sufficiently large $n$.

  Put $J=\overline{\mathbb{I}^g_{F^{\tt I}}+\mathbb{I}^h_{F^{\tt II}}}$.
  It remains to verify the defining relations modulo $J$.
  By definition of $J$, the factors satisfy their vacuum, equivariance and twisted Jacobi relations modulo $J$. Consequently, $F$ has the vacuum property and is $(g\otimes h)$-equivariant modulo $J$.

  Now we show that $F$ is local modulo $J$.
  Take homogeneous states $a_1,b_1\in V^{\tt I}$ and $a_2,b_2\in V^{\tt II}$.
  The locality argument in the proof of \zcref{prop:jacobi-criterion} shows that both factor assignments are local modulo $J$. Hence, there exist $K^{\tt I},K^{\tt II}\in\N$ such that
  \begin{align*}
    (z-w)^{K^{\tt I}}
    \ab [ F^{\tt I}[a_1,z], F^{\tt I}[b_1,w] ]
    &\equiv 0 \pmod J,\\
    (z-w)^{K^{\tt II}}
    \ab [ F^{\tt II}[a_2,z], F^{\tt II}[b_2,w] ]
    &\equiv 0 \pmod J.
  \end{align*}
  The lower bounds for the multiplied ordered products and commutativity
  between the two assignments justify the following coefficientwise
  calculation modulo $J$:
  \begin{align*}
    \MoveEqLeft
    (z-w)^{K^{\tt I}+K^{\tt II}}
      \ab[F[a_1\otimes a_2,z],F[b_1\otimes b_2,w]]\\
    &=(z-w)^{K^{\tt I}}
      \ab[F^{\tt I}[a_1,z],F^{\tt I}[b_1,w]]
      \,(z-w)^{K^{\tt II}}F^{\tt II}[a_2,z]F^{\tt II}[b_2,w]\\
    &\quad +(z-w)^{K^{\tt I}}F^{\tt I}[b_1,w]F^{\tt I}[a_1,z]
      \,(z-w)^{K^{\tt II}}
      \ab[F^{\tt II}[a_2,z],F^{\tt II}[b_2,w]]\\
    &\equiv0\pmod J.
  \end{align*}

  For each factor, applying \eqref{eq:jacobi-diagonal-expansion} on eigenstates and cancelling its scalar factor gives, by linearity, the branch-zero diagonal identity modulo $J$. Their diagonal poles can be cleared using $(z-w)^{K^{\tt I}}$ and $(z-w)^{K^{\tt II}}$, respectively.
  After clearing the poles, the same coefficientwise argument identifies the product of the two diagonal expansions with the diagonal expansion of the product modulo $J$.
  By the tensor-product vertex-operator formula and the definition of $F$, the product of the two iterates is the tensor iterate.
  Dividing by $(z-w)^{K^{\tt I}+K^{\tt II}}$ gives
  \[
    \iota_{w,z-w;0}\ab(F[a_1\otimes a_2,z]F[b_1\otimes b_2,w])
    \equiv F[Y(a_1\otimes a_2,z-w)(b_1\otimes b_2),w]\pmod J.
  \]
  By the first part of the proof of \zcref{prop:jacobi-criterion}, this diagonal identity and locality give the tensor twisted Jacobi relations modulo $J$. Together with the vacuum and equivariance relations, this proves $\mathbb{I}^{g\otimes h}_F\subset J$.
\end{proof}

\begin{proposition}\label{prop:enveloping-tensor-product}
  For finite-order automorphisms $g$ of $V^{\tt I}$ and $h$ of $V^{\tt II}$, respectively, there exists an isomorphism of graded topological rings
  \[
    \mathscr{U}_{g\otimes h}(V^{\tt I} \otimes V^{\tt II})
    \cong
    \mathscr{U}_{g}(V^{\tt I}) \ctensor \mathscr{U}_{h}(V^{\tt II}).
  \]
\end{proposition}
\begin{proof}
  Choose a common period $T$ of $g$ and $h$. First, we have two graded ring homomorphisms
  \[
    \iota^{\tt I} \colon \mathbb{U}_{\sfrac{1}{T}}(V^{\tt I}) \longrightarrow  \mathbb{U}_{\sfrac{1}{T}}(V^{\tt I}\otimes V^{\tt II})
    \txand
    \iota^{\tt II} \colon \mathbb{U}_{\sfrac{1}{T}}(V^{\tt II}) \longrightarrow \mathbb{U}_{\sfrac{1}{T}}(V^{\tt I}\otimes V^{\tt II})
  \]
  given by sending $J_n(a)$ to $J_n(a\otimes\vac)$ and sending $J_n(b)$ to $J_n(\vac\otimes b)$, respectively.
  Both maps preserve degree and thus extend continuously to the degreewise completions.
  For these extended maps, it is straightforward to check that
  \[
    \iota^{\tt I}(\widehat{\mathbb{I}}_g)
      \subset\widehat{\mathbb{I}}_{g\otimes h},
    \qquad
    \iota^{\tt II}(\widehat{\mathbb{I}}_h)
      \subset\widehat{\mathbb{I}}_{g\otimes h}.
  \]
  They therefore induce continuous ring homomorphisms
  \[
    \iota^{\tt I} \colon \mathscr{U}_{g}(V^{\tt I}) \longrightarrow \mathscr{U}_{g\otimes h}(V^{\tt I}\otimes V^{\tt II})
    \txand
    \iota^{\tt II} \colon \mathscr{U}_{h}(V^{\tt II}) \longrightarrow \mathscr{U}_{g\otimes h}(V^{\tt I}\otimes V^{\tt II}).
  \]
  Their images commute. Indeed, for eigenstates $a,b$, every product $\vo{(a\otimes\vac)}{n}(\vac\otimes b)$ with nonnegative $n$ vanishes.
  Taking $k=0$ in the component twisted Jacobi identity (\zcref{def:enveloping-algebra}) therefore shows that these modes commute. 

  The product $\iota^{\tt I} \otimes \iota^{\tt II} \colon \mathscr{U}_{g}(V^{\tt I})\otimes\mathscr{U}_{h}(V^{\tt II}) \longrightarrow \mathscr{U}_{g\otimes h}(V^{\tt I}\otimes V^{\tt II})$ is continuous. 
  Indeed, for any $x, u \in \mathscr{U}_{g}(V^{\tt I}), y, v \in \mathscr{U}_{h}(V^{\tt II})$ with $\deg u\le-i$ and $\deg v\le-j$, the aforementioned commutativity rewrites its image as  $\iota^{\tt I}(x)\iota^{\tt II}(y)\iota^{\tt I}(u)\iota^{\tt II}(v)$, where the last two factors have total degree at most $-i-j$. 
  The continuity of $\iota^{\tt I} \otimes \iota^{\tt II}$ allows us to extend it to a continuous ring homomorphism on the degreewise completed tensor product:
  \[
    \iota \colon \mathscr{U}_{g}(V^{\tt I}) \ctensor \mathscr{U}_{h}(V^{\tt II})
    \longrightarrow
    \mathscr{U}_{g\otimes h}(V^{\tt I}\otimes V^{\tt II}).
  \]

  Conversely, set $R=\mathscr{U}_{g}(V^{\tt I})\ctensor\mathscr{U}_{h}(V^{\tt II})$.
  The canonical inclusions from the two tensor factors to $R$ are continuous by \zcref{recol:acsr-tensor-products}, so
  \[
    F^{\tt I}[-,z]=Y_g[-,z]\otimes1,
    \qquad F^{\tt II}[-,z]=1\otimes Y_h[-,z]
  \]
  are graded field assignments for $V^{\tt I}$ and $V^{\tt II}$ respectively. 
  We have seen that they commute, and it is clear that  $\mathbb{I}^g_{F^{\tt I}}=\mathbb{I}^h_{F^{\tt II}}=0$.
  By \zcref{lem:tensor-fields}, the field $Y_g[a^{\tt I},z]\otimes Y_h[a^{\tt II},z]$ is well-defined.
  Taking its mode coefficients gives the assignment
  \[
    J_n(a^{\tt I}\otimes a^{\tt II})\longmapsto
    \sum_{k\in\sfrac{1}{T}\Z}J_k(a^{\tt I})\otimes J_{n-k}(a^{\tt II})
  \]
  It is bilinear in $a^{\tt I},a^{\tt II}$ and preserves degree, so it determines a graded ring homomorphism 
  \[
    \wp\colon\mathbb{U}_{\sfrac{1}{T}}(V^{\tt I}\otimes V^{\tt II})\to R.
  \]
  This map preserves the canonical neighborhoods and hence extends continuously to the degreewise completion.
  The relation-ideal containment in \zcref{lem:tensor-fields} shows that $\wp$ kills the defining relations.
  Hence, $\wp$ induces a continuous ring homomorphism
  \[
    \wp \colon \mathscr{U}_{g\otimes h}(V^{\tt I}\otimes V^{\tt II})
    \longrightarrow
    \mathscr{U}_{g}(V^{\tt I}) \ctensor \mathscr{U}_{h}(V^{\tt II}).
  \]

  It remains to check that $\iota$ and $\wp$ are inverse to each other.
  For $m, n\in \sfrac{1}{T}\Z$, we have
  \begin{align*}
    \wp \iota(J_n(a^{\tt I})\otimes J_m(a^{\tt II}))
    &=\wp(J_n(a^{\tt I}\otimes \vac)J_m(\vac\otimes a^{\tt II}))\\
    &=(J_n(a^{\tt I})\otimes 1) (1 \otimes J_m(a^{\tt II}))\\
    &=J_n(a^{\tt I})\otimes J_m(a^{\tt II}).
  \end{align*}
  Conversely, using the commutativity established above and
  \zcref{cor:enveloping-identities:product}, we obtain
  \begin{align*}
    \iota \wp(J_n(a^{\tt I}\otimes a^{\tt II})) 
    &= \iota\ab\Big(\sum_{k\in\sfrac{1}{T}\Z} J_k(a^{\tt I}) \otimes J_{n-k}(a^{\tt II}))\\
    &= \sum_{k\in\sfrac{1}{T}\Z} J_k(a^{\tt I}\otimes \vac) J_{n-k}(\vac\otimes a^{\tt II})\\
    &= J_n((a^{\tt I}\otimes \vac)_{-1}(\vac\otimes a^{\tt II}))\\
    &= J_n(a^{\tt I}\otimes a^{\tt II}).
  \end{align*}
  The displayed modes generate degreewise dense subalgebras on both sides. Hence $\iota$ and $\wp$ are inverse to each other.
\end{proof}

\begin{proposition}\label{prop:module-tensor-products}
  For finite-order automorphisms $g$ of $V^{\tt I}$ and $h$ of $V^{\tt II}$, tensoring modules induces equivalences
  \begin{align*}
    \cat{WMod}_g(V^{\tt I})\boxtimes
    \cat{WMod}_h(V^{\tt II})
    &\simeq
    \cat{WMod}_{g\otimes h}(V^{\tt I}\otimes V^{\tt II}),\\
    \cat{ExMod}_g(V^{\tt I})\boxtimes
    \cat{ExMod}_h(V^{\tt II})
    &\simeq
    \cat{ExMod}_{g\otimes h}(V^{\tt I}\otimes V^{\tt II}),
  \end{align*}
  compatible with the inclusions of exhaustive modules into weak modules.
\end{proposition}
\begin{proof}
  The weak equivalence amounts to identifying vector spaces carrying commuting weak actions of $V^{\tt I}$ and $V^{\tt II}$ with weak modules of the tensor VOA $V^{\tt I}\otimes V^{\tt II}$, by \zcref{recol:kelly-module-categories}.
  Restriction to the factors gives commuting weak actions by the twisted Jacobi identity, as in the proof of \zcref{prop:enveloping-tensor-product}.
  Conversely, a vector space with commuting weak actions, applying \zcref{lem:tensor-fields} with $R=\End(M)$, gives a weak action of $V^{\tt I}\otimes V^{\tt II}$ by multiplication of the factor fields.
  The constructions are inverse to each other due to the vacuum property and the following identity (which follows from \zcref{cor:enveloping-identities:product}):
  \[
    Y_M(a\otimes b,z)
    =Y_M(a\otimes\vac,z)Y_M(\vac\otimes b,z)
    \qquad(a\in V^{\tt I},\ b\in V^{\tt II}).
  \]

  The exhaustive equivalence follows from \zcref{rem:exhaustive-module-correspondence,prop:enveloping-tensor-product,recol:acsr-tensor-products}. The compatibility with the inclusions of exhaustive modules into weak modules is clear from the construction.
\end{proof}

\part{Permutation orbifold theory}

In this part, we dive into the permutation orbifold theory. We start from a single cycle and then generalize to arbitrary permutation.

\section{Set-up for permutation orbifold theory}
\label{sec:set-up}

We start with the following data:
\begin{itemize}
  \item $(\nota{V},\nota{Y(\placeholder,z)},\nota{\vac},\nota{\cfv})$, a vertex operator algebra;
  \item $\nota{V^{\otimes N}} := V\otimes\cdots\otimes V$ the tensor product of $N$ copies of $V$;
  \item $\nota{\sigma}$, a cyclic permutation on $V^{\otimes N}$, without loss of generality, we may assume \[\sigma = (12\cdots N);\]
  \item $\nota{g}$, an automorphism of $V$, also acting on $V^{\otimes N}$ diagonally; and
  \item $\nota{T}$, a period of $g$, namely $g^T = \id_V$.
\end{itemize}
Then, $V^{\otimes N}$ is a vertex operator algebra with the vertex operator
\begin{equation}
  \label{eq:vertex-operator-tensor}
  A=a_1\otimes\cdots\otimes a_N \in V^{\otimes N} \longmapsto Y(A,z) := Y(a_1,z) \otimes \cdots \otimes Y(a_N,z).
\end{equation}
The permutation $\sigma$ and the diagonal action of $g$ give rise to two commuting automorphisms of this vertex operator algebra. 
In what follows, put $\nota{\tau}=\sigma g$ and $\nota{h}=g^N$, with periods $NT$ and $T$, respectively.
The next two sections construct homomorphisms in both directions between the twisted enveloping algebras
\[
  \mathscr{U}_{\tau}(V^{\otimes N})
  \qquad\text{and}\qquad
  \mathscr{U}_{h}(V).
\]
To state the equivariance identities for the fields $Y_h$ and $Y_\tau$ simultaneously, we fix $\nota{\zeta}$ to be a primitive $NT$-th root of unity with $\zeta^N=\zeta_T$.
The coordinate changes $\varsigma_{T}$ and $\varsigma_{NT}$ are:
\[
  \varsigma_{T}(z) \colon z^{\sfrac{1}{T}}
  \longmapsto
  \zeta^{N}z^{\sfrac{1}{T}}
  \txand
  \varsigma_{NT}(z) \colon z^{\sfrac{1}{NT}}
  \longmapsto
  \zeta z^{\sfrac{1}{NT}}.
\]

\begin{remark}
  We recommend that the reader first treat the case $g = \id_V$, $T = 1$, to become familiar with the strategy.
\end{remark}

\subsection{Very simple tensors}

We begin with a special class of simple tensors in $V^{\otimes N}$, whose modes serve as topological generators of the twisted enveloping algebra $\mathscr{U}_\tau(V^{\otimes N})$.

\begin{notation}
  Simple tensors of the form $a\otimes\vac\otimes\cdots\otimes \vac$, where $a\in V$, will play a critical role in this paper. So we give it a specific notation:
  \begin{equation}
    \nota{\mathfrak{x}_{a}} := a\otimes\vac\otimes\cdots\otimes \vac \in V^{\otimes N}.
  \end{equation}
  More generally, we use $\nota{\mathfrak{x}_{a}^{i}}$ to denote the simple tensor whose $i$-th component is $a$ and all other components are $\vac$. Note that, by our convention,
  \begin{equation}
    \mathfrak{x}_{a}^{i} = \sigma^{i-1}(a\otimes\vac\otimes\cdots\otimes \vac)\in V^{\otimes N}.
  \end{equation}
  So we take the above as the definition of $\mathfrak{x}_{a}^{i}$ if $i$ is out of the range $1,\cdots,N$. 
  These simple tensors are called \concept{very simple tensors} in this paper.
\end{notation}

It turns out that all modes can be built from those of very simple tensors. 
\begin{proposition}
  \label{lem:tensor-generation}
  The modes of very simple tensors topologically generate $\mathscr{U}_\tau(V^{\otimes N})$. 
  More generally, a (graded) field assignment with values in a (graded) topological associative algebra and satisfying the twisted Jacobi identity is uniquely determined by its fields on very simple tensors.
\end{proposition}
\begin{proof}
  Let $F$ be such a field assignment on $V^{\otimes N}$.
  For a homogeneous $\tau$-eigenstate $A$, let $\alpha$ be the natural number below $NT$ such that $\zeta^{\alpha}$ is its eigenvalue, and let $B\in V^{\otimes N}$ be homogeneous.
  The $\tau$-twisted Jacobi identity then gives
  \begin{equation*}
    \label{eq:tensor-generation-product}
    \tag{\textasteriskcentered}
    F[\vo{A}{-1}B,z] = 
    \normord{F[A,z]F[B,z]} - 
    \sum_{j \ge 0}\binom{\frac{\alpha}{(NT)}}{j+1}z^{-j-1}F[\vo{A}{j}B,z],
  \end{equation*}
  where normal ordering places $F_n(A)$ to the left when $n+\wt A<1$ and to the right when $n+\wt A\ge1$.  
  In the part $n+\wt A<1$, the rightmost mode index tends to $+\infty$ as $n\to-\infty$; in the other part, it does so as $n\to+\infty$.
  Thus only finitely many terms survive modulo any open left ideal.
  Completeness gives summability; in the graded case the terms have the same degree, so degreewise completeness suffices.
  On the other hand, the correction is finite by truncation.  

  By linearity, it suffices to prove that the fields of any simple tensor $A=a_1\otimes\cdots\otimes a_N$ with each $a_i$ homogeneous are determined by those of very simple tensors. 
  We prove this by induction on the number $k$ of non-vacuum factors in $A$.
  The base case $k \le 1$ is immediate, so we suppose $k \ge 2$. 
  Choose $1 \le i \le N$ with $a_i\ne\vac$, put $a=a_i$, and let $B$ be obtained from $A$ by replacing its $i$-th factor with $\vac$.
  Then $B$ has less non-vacuum factors and $A = \vo{(\mathfrak{x}_a^i)}{-1}B$.
  Applying \eqref{eq:tensor-generation-product} to the $\tau$-eigencomponents of $\mathfrak{x}_a^i$, with $B$ as the second state, and summing gives
  \[
    F[A,z] = F[\vo{(\mathfrak{x}_a^i)}{-1}B,z] = \normord{F[\mathfrak{x}_a^i,z]F[B,z]} - \sum_{\alpha = 0}^{NT-1}\sum_{j \ge 0}\binom{\frac{\alpha}{(NT)}}{j+1}z^{-j-1}F[\vo{(\mathfrak{x}_a^i)^{(\alpha)}}{j}B,z].
  \]
  By the identity $\tau\mathfrak{x}_a^i = \mathfrak{x}_{ga}^{i+1}$, the $\tau$-eigen components of $\mathfrak{x}_a^i$ are linear combinations of very simple tensors. 
  For any very simple tensor, say $\mathfrak{x}_c^k$, its mode $\vo{(\mathfrak{x}_c^k)}{j}$ with $j \ge 0$ acts on $B$ by changing only the $k$-th factor and vanishes when that factor is $\vac$, since $\vo{c}{j}\vac=0$.
  Hence, by the induction hypothesis, all the fields on the right-hand side are determined by those of very simple tensors. Thus, so is the field $F[A,z]$.
  For $F=Y_\tau$, the same induction places every simple mode in the closed subalgebra generated by the modes of very simple tensors, proving the first assertion.
\end{proof}

\subsection{The $\Delta_N$-operator}

In \cite{BDM02}, the operator $\Delta_N(z)$ is defined and plays a crucial role in connecting $\sigma$-twisted $V^{\otimes N}$-modules and untwisted $V$-modules.
We use it to relate the twisted enveloping algebras $\mathscr{U}_{\sigma g}(V^{\otimes N})$ and $\mathscr{U}_{g^N}(V)$.

\begin{definition}
  The $\Delta_N$-operator is defined as follows:
  \[
    \nota{\Delta_N(z)} :=
    \fun{exp}(
      \sum_{j=1}^{\infty}c_{j}z^{-\sfrac{j}{N}}\VL{j}
    )N^{-\VL{0}}z^{(\sfrac{1}{N}-1)\VL{0}}
    \in (\End{V})\dbrack{z^{\pm\sfrac{1}{N}}},
  \]
  where the coefficients $c_j$ are uniquely determined by the coordinate change:
  \[
    x \longmapsto
    \fun{exp}(
      -\sum_{j=1}^{\infty}c_{j}x^{j+1}\pdv{x}
    )x
    = \tfrac{1}{N}(1+x)^{N} - \tfrac{1}{N}.
  \]
  Since $\VL{j}$ lowers weight by $j$ and $V$ is lower-bounded, the exponential is locally finite when acting on $V$. In particular, $\Delta_N(1)$ is a well-defined invertible operator on $V$.
\end{definition}

The following properties of $\Delta_N(z)$ will be used in this paper. For self-containedness, we provide the proofs here.
\begin{lemma}
  \label{lem:delta-vacuum}
  $\Delta_N(z)\vac = \vac$.
\end{lemma}
\begin{proof}
  This follows from $\VL{j}\vac = 0$ for all $j\ge 0$.
\end{proof}
\begin{lemma}
  \label{lem:delta-scaling}
  We have the following identity in $(\End{V})\dbrack{z^{\pm\sfrac{1}{N}}}$:
  \[
    \Delta_N(z) = z^{\sfrac{1}{N}\VL{0}}\Delta_N(1)z^{-\VL{0}}
  \]
\end{lemma}
\begin{proof}
  Since $\wt(\VL{j}a)=\wt a-j$ for $j\in \Z$, we have
  \begin{equation}
    z^{\sfrac{1}{N}\VL{0}}(\sum_{j=1}^{\infty}c_{j}\VL{j})=(\sum_{j=1}^{\infty}c_{j}z^{-\sfrac{j}{N}}\VL{j})z^{\sfrac{1}{N}\VL{0}}.
  \end{equation}
  Thus,
  \begin{equation}
    z^{\sfrac{1}{N}\VL{0}}\fun{exp}(\sum_{j=1}^{\infty}c_{j}\VL{j})=\fun{exp}(\sum_{j=1}^{\infty}c_{j}z^{-\sfrac{j}{N}}\VL{j})z^{\sfrac{1}{N}\VL{0}}.
  \end{equation}
  The desired identity follows immediately.
\end{proof}

\begin{lemma}\label{lem:delta-covariance}
  On the diagonal branch $z^{\sfrac{1}{N}}=w^{\sfrac{1}{N}}$, use $x=z^{\sfrac{1}{N}}-w^{\sfrac{1}{N}}$ as a local coordinate. Then the following identity holds as an $\End(V)\dparen{w^{-\sfrac{1}{N}}}$-valued formal distribution in $x$:
  \[
    \Delta_N(w)Y(a,z-w)\Delta_N(w)^{-1}
    =Y(\Delta_N(z)a,z^{\sfrac{1}{N}}-w^{\sfrac{1}{N}})
    \txforall a\in V.
  \]
\end{lemma}
\begin{proof}
  Put $\rho_w(u)=(w^{\sfrac{1}{N}}+u)^N-w$, and let $\mathcal{R}_{\rho_w}$ be its coordinate-change operator with the conventions of \cite[\S 6.3.1]{FBZ04}.
  The coordinate $\rho_w$ is obtained by successively applying $u\mapsto w^{-\sfrac{1}{N}}u$, $u\mapsto ((1+u)^N-1)/N$, and $u\mapsto Nwu$.
  Hence, the defining relation for $c_j$ gives
  \[
    \mathcal{R}_{\rho_w}
    =w^{\sfrac{1}{N}\VL{0}}\exp\ab\Big(\sum_{j\ge1}c_j\VL{j})(Nw)^{-\VL{0}}
    =\Delta_N(w).
  \]
  For $x=z^{\sfrac{1}{N}}-w^{\sfrac{1}{N}}$, we have $\rho_w(x)=z-w$ and
  \[
    \rho_w(x+u)-\rho_w(x)=(z^{\sfrac{1}{N}}+u)^N-z=\rho_z(u).
  \]
  Thus the coordinate-change formula \cite[Lemma 6.5.6]{FBZ04} reads
  \[
    Y(a,x)=\mathcal{R}_{\rho_w}Y(\mathcal{R}_{\rho_z}^{-1}a,z-w)\mathcal{R}_{\rho_w}^{-1}.
  \]
  Note that its proof applies over $\k\dparen{w^{-\sfrac{1}{N}}}$ since the calculation uses only formal identities and rational coefficients.
  Replacing $a$ by $\Delta_N(z)a$ in the above equation gives the desired identity.
\end{proof}

\subsection{Products of fields and branch specialization}

In the next two sections, we form products of fields and specialize their variables to distinct branches. The following lemma ensures that these products and specializations are well-defined.

\begin{lemma}
\label{lem:field-products}
  Let $R$ be a (graded) topological associative algebra in the sense of \zcref{rec:topological-associative-algebras}, and let $F[-,u]\colon V\to\mathcal{D}_R(u^{\sfrac{1}{T}})$ be a (graded) field assignment satisfying $h$-equivariance and the $h$-twisted Jacobi identity.
  For homogeneous states $a_1,\ldots,a_m\in V$ and independent variables $\vect{u}=(u_1,\ldots,u_m)$, choose $K_{ij}\in\N$ such that
  \[
    (u_i-u_j)^{K_{ij}}\ab[F[a_i,u_i],F[a_j,u_j]]=0
    \qquad(i<j),
  \]
  and put
  \[
    P(\vect{u}) := \prod_{i<j}(u_i-u_j)^{K_{ij}}.
  \]
  Then:
  \begin{enumerate}
    \item\label{lem:field-products:a} 
        The product distribution $P(\vect{u})F[a_1,u_1]\cdots F[a_m,u_m]$ is independent of the ordering, and its coefficients tend uniformly to zero as any exponent tends to $-\infty$. 
    \item\label{lem:field-products:b} 
        The local expansion of the ordered product $F[a_1,u_1]\cdots F[a_m,u_m]$ at the diagonal branch $u_i^{\sfrac{1}{T}}=\zeta_T^{-r}u_j^{\sfrac{1}{T}}$, in the sense of \zcref{rem:diagonal-poles}, is obtained by placing the selected pair adjacent and replacing it by
        \[
          F[Y(h^ra_i,u_i-u_j)a_j,u_j],
        \]
        with the other fields retained.
    \item\label{lem:field-products:c} 
        If the selected pairs share no indices and collisions occur only within them, applying \itemref{lem:field-products:b} to each pair gives a joint expansion independent of their order. The corresponding differences $u_i-u_j$ may be set equal.
    \item\label{lem:field-products:d} 
        For a formal parameter $x$ and $\lambda_i\in\k^\times$ with pairwise distinct $\lambda_i^T$, specializing the distribution in \itemref{lem:field-products:a} by
        \[
          u_i=\lambda_i^T x,\qquad
          u_i^{\sfrac{1}{T}}=\lambda_i x^{\sfrac{1}{T}},
        \]
        we obtain a well-defined distribution. 
        Moreover, its quotient by $P(\lambda_1^T x,\ldots,\lambda_m^T x)$ is an $R$-valued field independent of the ordering and the $K_{ij}$.
  \end{enumerate}
  These assertions extend to arbitrary states by multilinearity. They also hold when each $a_i$ has finite Laurent-polynomial dependence on its own variable, interpreted termwise.
\end{lemma}
For simplicity, we will denote the resulting field in \itemref{lem:field-products:d} by
\[
  \nota{F[a_1,u_1]\cdots F[a_m,u_m]\Big|_{u_i^{\sfrac{1}{T}}=\lambda_i x^{\sfrac{1}{T}}}}.
\]
\begin{proof}
  For \itemref{lem:field-products:a}, fix an open left ideal $W\subset R$.
  To obtain a uniform lower bound on the powers of $u_i$ modulo $W$, place its field rightmost.
  Multiplication by the polynomial $P(\vect{u})$ cannot decrease this lower bound, and locality identifies these distributions for all orderings, so these bounds hold simultaneously.

  For \itemref{lem:field-products:b}, first expand the product multiplied by $P$ on the indicated branch, using \zcref{lem:diagonal-expansion} and \zcref{nota:relative-boundary-g}. With exponents of the other variables fixed, \itemref{lem:field-products:a} gives uniform lower bounds in the selected two variables. Completeness gives the coefficient sums; in the graded case, homogeneity puts their terms in one degree, so degreewise completeness suffices.

  Place the selected pair adjacent and fix the other modes. 
  Then expand the pair using \eqref{eq:jacobi-diagonal-expansion} and $h$-equivariance. Multiplication by the fixed modes is continuous, and the remaining factors of $P$ produce only finitely many exponent shifts. Hence the multiplied product has the expansion obtained by replacing the pair as stated and retaining the expansion of $P$ and the other fields.

  Finally, multiply by the expansion of $P^{-1}$ on the same branch to recover the original product's local expansion: the pair factor contributes $(u_i-u_j)^{-K_{ij}}$, while the other factors are units where no additional collisions occur.

  For \itemref{lem:field-products:c}, consider the joint expansion of the product multiplied by $P$ in the selected differences $u_i-u_j$.
  The monomials contributing to each coefficient have a fixed sum of exponents, so by \itemref{lem:field-products:a} only finitely many terms contribute modulo every $W$. Their sum therefore exists in $R$ and does not depend on the order of expansion. 
  By \itemref{lem:field-products:b}, division by the expansion of $P$ gives the product with each selected pair replaced by the corresponding iterate and the factors for the other pairs are units at the collision locus.

  Each series $Y(h^ra_i,u_i-u_j)a_j$ is lower truncated in its own difference.
  Thus, after setting the differences equal, only finitely many exponent tuples contribute to any fixed power. This shows the last assertion of \itemref{lem:field-products:c}.

  Now, we prove \itemref{lem:field-products:d}. 
  Under that substitution, only monomials whose exponents have a fixed sum contribute to a given coefficient.
  By \itemref{lem:field-products:a}, there are only finitely many such terms modulo every $W$. Hence the sum exists in $R$.
  The resulting powers of $x$ are bounded below modulo $W$ by the bounds in \itemref{lem:field-products:a}.

  Since $P$ specializes to the invertible monomial
  \[
    \prod_{i<j}(\lambda_i^T-\lambda_j^T)^{K_{ij}}x^{\sum_{i<j}K_{ij}},
  \]
  division by it only shifts the powers of $x$. Thus the quotient is an $R$-valued field, independent of the ordering by \itemref{lem:field-products:a}. 
  It is also independent of the $K_{ij}$, since any other choice gives a product differing by a monomial in the differences $u_i-u_j$, which specializes to a unit.

  Finally, all four assertions extend to arbitrary states by multilinearity. 
  For states with finite Laurent-polynomial dependence on their own variables, expand termwise and factor out the scalar monomials. These finite sums and exponent shifts preserve the preceding bounds and summability. Hence, the assertions hold in this case as well.
\end{proof}

\section{From \texorpdfstring{$\mathscr{U}_{\sigma g}(V^{\otimes N})$}{sigma g-twisted} to \texorpdfstring{$\mathscr{U}_{g^N}(V)$}{g\textasciicircum N-twisted}}
\label{sec:varphi}

We first construct a $\mathscr{U}_h(V)$-valued field assignment on $V^{\otimes N}$, and then use it to define the required continuous homomorphism.

\begin{construction}\label{con:tensor-field}
  Let $R$ be a (graded) topological associative algebra, and let $F[-,u]\colon V\to\mathcal D_R(u^{\sfrac{1}{T}})$ be a (graded) field assignment satisfying the vacuum property, $h$-equivariance and the $h$-twisted Jacobi identity.
  For $a\in V$, set
  \begin{equation}\label{eq:transport-first-position}
    \mathcal{T}_N[\mathfrak{x}_a,z] = F[\Delta_N(z)a,z^{\sfrac{1}{N}}].
  \end{equation}
  We extend this assignment to a field assignment on $V^{\otimes N}$ as follows.
  Note that this field assignment should be $\tau$-equivariant. 
  Since $\tau=\sigma g$, applying $\tau^{i-1}$ moves the first tensor position to position $i$ and applies $g^{i-1}$ to its state. 
  Therefore, to obtain $a$ in position $i$, we start with $g^{1-i}a$ in the first position and then: $\mathfrak{x}_a^i=\tau^{i-1}\mathfrak{x}_{g^{1-i}a}$.

  Treating $u$ as an $N$-th root of $u^N$, \zcref{lem:delta-scaling} gives $\Delta_N(u^N) = u^{\VL{0}}\Delta_N(1)u^{-N\VL{0}}$.
  Define the field $B_i[a,u]$ to be 
  \begin{equation}\label{eq:B_i}
    \nota{B_i[a,u]}=F[\Delta_N(u^N)g^{1-i}a,u] = \sum_{n\in\sfrac{1}{T}\mathbb{Z}}F_n(\Delta_N(1)g^{1-i}a)u^{-n-N\wt a}.
  \end{equation}
  This will be the $i$-th factor in the product of fields for a simple tensor $A=a_1\otimes\cdots\otimes a_N$. 
  Namely, we will consider the product $B_1[a_1,u_1]\cdots B_N[a_N,u_N]$ and then specialize the $u_i$ to distinct branches of $z^{\sfrac{1}{N}}$.  
  For simplicity, put $\nota{\eta} = \zeta^T$ and take $x=z^{\sfrac{1}{N}}$, with $x^{\sfrac{1}{T}} = z^{\sfrac{1}{NT}}$.
  Apply \zcref{lem:field-products} to $F$ with $m=N$ and states $\Delta_N(u_i^N)g^{1-i}a_i$, and specialize to the branches
  \[
    u_i^{\sfrac{1}{T}} = \zeta^{1-i}x^{\sfrac{1}{T}}, \qquad u_i=\eta^{1-i}x.
  \]
  We thus can define the field $\mathcal{T}_N[A,z]$ as 
  \begin{equation}\label{eq:tensor-field}
    \mathcal{T}_N[A,z] = B_1[a_1,u_1]\cdots B_N[a_N,u_N]\Big|_{u_i^{\sfrac{1}{T}}=\zeta^{1-i}x^{\sfrac{1}{T}}}.
  \end{equation}
  By multilinearity, we obtain a linear map
  \[
    \nota{\mathcal{T}_N[-,z]} \colon V^{\otimes N} \longrightarrow \mathcal D_R(z^{\sfrac{1}{NT}}).
  \]
  Since $B_i[\vac,u]=1$, it has the prescribed restriction on $\mathfrak{x}_a$.
  When $R$ is graded, $\mathcal T_N$ is graded after dividing the degrees in $R$ by $N$.
\end{construction}

\begin{proposition}
  \label{prop:tensor-field-jacobi}
  The field assignment defined by \eqref{eq:tensor-field} has the vacuum property, $\tau$-equivariance, and the $\tau$-twisted Jacobi identity.
\end{proposition}
\begin{proof}
  Vacuum is immediate. 
  For a simple tensor $A$, the $i$-th factor in $\mathcal{T}_N[\tau A,\varsigma_{NT}z]$ equals the $(i-1)$-st factor in $\mathcal{T}_N[A,z]$ whenever $2\le i\le N$.
  For $i=1$, $h$-equivariance of $F$ gives
  \[
    B_1[ga_N,u_1]\Big|_{u_1^{\sfrac{1}{T}}=\zeta x^{\sfrac{1}{T}}}
    =B_N[a_N,u_N]\Big|_{u_N^{\sfrac{1}{T}}=\zeta^{1-N}x^{\sfrac{1}{T}}}.
  \]
  This proves $\tau$-equivariance of $\mathcal{T}_N[-,z]$.

  Next, we prove locality. Take $A=a_1\otimes\cdots\otimes a_N$ and $B=b_1\otimes\cdots\otimes b_N$ with homogeneous factors.
  Apply \zcref{lem:field-products:a} to $F$ with $m=2N$ and states $\Delta_N(u_i^N)g^{1-i}a_i$ and $\Delta_N(v_i^N)g^{1-i}b_i$ for $1\le i\le N$.
  Choose $K$ large enough that, for $1\le i,j\le N$,
  \[
    (u_i-v_j)^K B_i[a_i,u_i]B_j[b_j,v_j]
    =(u_i-v_j)^K B_j[b_j,v_j]B_i[a_i,u_i].
  \]
  By \zcref{lem:field-products:a}, after multiplying by suitable diagonal factors, we may move all $B_i[a_i,u_i]$ past all $B_j[b_j,v_j]$.
  Specializing $u_i=\eta^{1-i}z^{\sfrac{1}{N}}$ and $v_i=\eta^{1-i}w^{\sfrac{1}{N}}$, the factor $\prod_{i,j}(u_i-v_j)^K$ becomes $c(z-w)^{NK}$ with $c\ne0$. Hence, we obtain
  \[
    (z-w)^{NK}\mathcal{T}_N[A,z]\mathcal{T}_N[B,w] = (z-w)^{NK}\mathcal{T}_N[B,w]\mathcal{T}_N[A,z].
  \]
  This proves locality.

  Finally, we verify the diagonal identity \eqref{eq:jacobi-diagonal-expansion}. 
  For the diagonal expansions, set $u_i=\eta^{1-i}z^{\sfrac{1}{N}}$ and $v_i=\eta^{1-i}w^{\sfrac{1}{N}}$.
  Each $u_i-v_i$ is $z-w$ times a unit, while $u_i-u_j$, $v_i-v_j$, and $u_i-v_j$ are units for $i\ne j$.
  Thus \zcref{lem:field-products:c} and \itemref{lem:field-products:d} justify the simultaneous expansion and specialization.
  With these in mind, we compute the diagonal expansion as follows:
  \[
    \begin{aligned}
      \MoveEqLeft
      \iota_{w,z-w;0}\ab(\mathcal{T}_N[A,z]\mathcal{T}_N[B,w])\\
      \overset{\text{\eqref{eq:tensor-field}}}&=
        \iota_{w,z-w;0}\ab(
          \prod_{i=1}^{N}B_i[a_i,u_i]
          \prod_{j=1}^{N}B_j[b_j,v_j]
          \Bigg|_{\substack{
            u_i^{\sfrac{1}{T}}=\zeta^{1-i}z^{\sfrac{1}{NT}}\\
            v_j^{\sfrac{1}{T}}=\zeta^{1-j}w^{\sfrac{1}{NT}}}}
        )\\
      \overset{\text{\eqref{eq:B_i}}}&=
        \iota_{w,z-w;0}\ab(
          \prod_{i=1}^{N}F[\Delta_N(u_i^N)g^{1-i}a_i,u_i]
          \prod_{j=1}^{N}F[\Delta_N(v_j^N)g^{1-j}b_j,v_j]
          \Bigg|_{\substack{
            u_i^{\sfrac{1}{T}}=\zeta^{1-i}z^{\sfrac{1}{NT}}\\
            v_j^{\sfrac{1}{T}}=\zeta^{1-j}w^{\sfrac{1}{NT}}}}
        )\\
      \overset{\text{\zcref[noname]{lem:field-products:c}}}&=
        \prod_{i=1}^{N}F[
          Y(\Delta_N(u_i^N)g^{1-i}a_i,u_i-v_i)
          \Delta_N(v_i^N)g^{1-i}b_i,v_i]
        \Bigg|_{\substack{
          u_i^{\sfrac{1}{T}}=\zeta^{1-i}z^{\sfrac{1}{NT}}\\
          v_i^{\sfrac{1}{T}}=\zeta^{1-i}w^{\sfrac{1}{NT}}}}\\
      \overset{\text{\zcref[noname]{lem:delta-covariance}}}&=
        \prod_{i=1}^{N}F[
          \Delta_N(v_i^N)Y(g^{1-i}a_i,u_i^N-v_i^N)g^{1-i}b_i,v_i]
        \Bigg|_{\substack{
          u_i^{\sfrac{1}{T}}=\zeta^{1-i}z^{\sfrac{1}{NT}}\\
          v_i^{\sfrac{1}{T}}=\zeta^{1-i}w^{\sfrac{1}{NT}}}}\\
      &=\prod_{i=1}^{N}
        F[\Delta_N(v_i^N)g^{1-i}Y(a_i,z-w)b_i,v_i]
        \Big|_{v_i^{\sfrac{1}{T}}=\zeta^{1-i}w^{\sfrac{1}{NT}}}\\
      \overset{\text{\eqref{eq:B_i}}}&=
        \prod_{i=1}^{N}B_i[Y(a_i,z-w)b_i,v_i]
        \Big|_{v_i^{\sfrac{1}{T}}=\zeta^{1-i}w^{\sfrac{1}{NT}}}\\
      \overset{\text{\eqref{eq:tensor-field}}}&=
        \mathcal{T}_N[Y(a_1,z-w)b_1\otimes\cdots\otimes Y(a_N,z-w)b_N,w]\\
      &=\mathcal{T}_N[Y_{V^{\otimes N}}(A,z-w)B,w].
    \end{aligned}
  \]
  Together with locality and $\tau$-equivariance, this identity gives the $\tau$-twisted Jacobi identity by the argument following \eqref{eq:jacobi-diagonal-expansion} in the proof of \zcref{prop:jacobi-criterion}.
\end{proof}

For the rest of this section, take $R=\mathscr U_h(V)$ and
$F[-,u]=Y_h[-,u]$ from \zcref{conv:formal-vertex-operator}.

\begin{construction}\label{con:phi}
  For $n\in\sfrac{1}{T}\Z$ and homogeneous $A\in V^{\otimes N}$, define $\varphi_0(J_{\sfrac{n}{N}}(A))$ to be the coefficient of $z^{-\sfrac{n}{N}-\wt A}$ in $\mathcal{T}_N[A,z]$.
  This coefficient lies in $\mathscr{U}_h(V)_{-n}$.
  By linearity and freeness, this defines a unique unital algebra homomorphism
  \[
    \nota{\varphi_0}\colon\mathbb{U}_{\sfrac{1}{NT}}(V^{\otimes N})
    \longrightarrow\mathscr{U}_h(V).
  \]
  For very simple tensors,
  \[
    \varphi_0(J_{\sfrac{n}{N}}(\mathfrak{x}_a^i))
    =\zeta^{nT(i-1)}J_n(\Delta_N(1)g^{1-i}a).
  \]
\end{construction}

\begin{remark}
  When $T = 1$, namely $g = \id_V$, the definition of $\varphi_0$ reduces to:
  \[
    \varphi_0(J_{\sfrac{n}{N}}(\mathfrak{x}_{a}^{i})) = \zeta^{nT(i-1)}J_{n}(\Delta_N(1)a)
    \txforall a \in V, n \in \sfrac{1}{T}\Z.
  \]
\end{remark}

\begin{lemma}\label{lem:continuous-phi}
  The assignment $\varphi_0$ is continuous and extends to the degreewise completion.
\end{lemma}
\begin{proof}
  Since $\varphi_0$ multiplies degrees by $N$, \zcref[noname]{A2} give
  \[
    \varphi_0(\cNL[\sfrac{s}{N}]
      \mathbb{U}_{\sfrac{1}{NT}}(V^{\otimes N}))
    \subset \NL[s]\mathscr{U}_h(V)
    \txand
    \varphi_0(\cNR[\sfrac{s}{N}]
      \mathbb{U}_{\sfrac{1}{NT}}(V^{\otimes N}))
    \subset \NR[s]\mathscr{U}_h(V).
  \]
  Hence $\varphi_0$ is continuous. The target is separated and degreewise complete, so $\varphi_0$ extends uniquely to the degreewise completion. 
\end{proof}

\begin{proposition}\label{prop:phi-homomorphism}
  The assignment $\varphi_0$ induces a continuous graded homomorphism
  \[
    \nota{\varphi}\colon\mathscr{U}_{\sigma g}(V^{\otimes N})
    \longrightarrow\mathfrak{d}_N^{\ast}\mathscr{U}_{g^N}(V).
  \]
\end{proposition}
\begin{proof}
  The continuous extension of $\varphi_0$ given by \zcref{lem:continuous-phi} commutes with the fixed-degree sums in the defining Jacobi relations.
  By \zcref{prop:tensor-field-jacobi}, the vacuum, equivariance and twisted Jacobi relations therefore map to zero.
  Since $\mathscr{U}_h(V)$ is separated, the kernel of this extension is a closed two-sided ideal containing $\widehat{\mathbb{I}}_{\tau}$.
  Thus this extension descends continuously to $\mathscr{U}_{\tau}(V^{\otimes N})$.
\end{proof}

\section{From \texorpdfstring{$\mathscr{U}_{g^N}(V)$}{g\textasciicircum N-twisted} to \texorpdfstring{$\mathscr{U}_{\sigma g}(V^{\otimes N})$}{sigma g-twisted}}
\label{sec:psi}

We use $\Delta_N^{-1}$ to construct a $\mathscr{U}_\tau(V^{\otimes N})$-valued field assignment on $V$.

\begin{proposition}\label{prop:inverse-field}
  Let $R$ be a (graded) topological associative algebra, and let $F[-,z]\colon V^{\otimes N}\to\mathcal D_R(z^{\sfrac1{NT}})$ be a (graded) field assignment satisfying the vacuum property, $\tau$-equivariance and the $\tau$-twisted Jacobi identity.
  Treat $x$ and $x^{\sfrac{1}{T}}$ as an $N$-th root and an $NT$-th root of $x^N$, respectively. Define a field assignment $\nota{\mathcal S_N[-,x]}\colon V\to\mathcal D_R(x^{\sfrac{1}{T}})$ by
  \begin{equation}\label{eq:inverse-transport-field}
    \mathcal{S}_N[a,x]=F[\mathfrak{x}_{\Delta_N(x^N)^{-1}a},x^N].
  \end{equation}
  This assignment has the vacuum property, $h$-equivariance and the $h$-twisted Jacobi identity.
  When $R$ is graded, $\mathcal S_N$ is graded after multiplying the degrees
  in $R$ by $N$.
\end{proposition}
\begin{proof}
  By \zcref{lem:delta-scaling}, we have $\Delta_N(x^N)^{-1}=x^{N\VL{0}}\Delta_N(1)^{-1}x^{-\VL{0}}$.
  For homogeneous $a$, this gives
  \[
    \mathcal{S}_N[a,x] = 
    \sum_{n\in\sfrac{1}{T}\Z}
    F_{\sfrac{n}{N}}(\mathfrak{x}_{\Delta_N(1)^{-1}a})x^{-n-\wt a}.
  \]
  These coefficients tend to zero as $n\to+\infty$, since $\Delta_N(1)^{-1}a$ is a finite sum of homogeneous states and $F$ is a field assignment. 
  This proves well-definedness of $\mathcal{S}_N[-,x]$.
  The vacuum property follows from \zcref{lem:delta-vacuum}. Since $\tau^N$ acts diagonally as $h$ and $\Delta_N^{-1}$ commutes with $h$, the $\tau$-equivariance of $F$ gives $h$-equivariance of $\mathcal{S}_N$.

  We next prove locality. For homogeneous $a,b$, let $\alpha$ and $\beta$ range over the coefficients of the finite Laurent polynomials $\Delta_N(x^N)^{-1}a$ and $\Delta_N(y^N)^{-1}b$, respectively. 
  Applying locality of $F$ to the finitely many pairs $F[\mathfrak{x}_\alpha,z]$, $F[\mathfrak{x}_\beta,w]$, taking a common $K\in\N$, and substituting $z=x^N$, $w=y^N$ gives
  \begin{equation}
    \label{eq:S-locality}\tag{\textasteriskcentered}
    (x^N-y^N)^K\mathcal{S}_N[a,x]\mathcal{S}_N[b,y]
    =(x^N-y^N)^K\mathcal{S}_N[b,y]\mathcal{S}_N[a,x].
  \end{equation}
  Modulo any open left ideal $W\subset R$, \zcref{lem:field-products:a} gives lower bounds for the exponents of $x^{\sfrac{1}{T}}$ and $y^{\sfrac{1}{T}}$ in this common distribution.
  For each fixed sum of exponents, it therefore gives a finite Laurent polynomial. 
  Multiplying this polynomial by the two ordered expansions of $(x^N-y^N)^{-K}$ recovers the corresponding parts of $\mathcal S_N[a,x]\mathcal S_N[b,y]$ and $\mathcal S_N[b,y]\mathcal S_N[a,x]$ modulo $W$.

  Now, we show that the factor $(x^N-y^N)^K$ in \eqref{eq:S-locality} can be replaced by $(x-y)^K$. 
  For this, we show that $\mathcal{S}_N[a,x]\mathcal{S}_N[b,y]$ has no poles at the branches $x^{\sfrac{1}{T}}=\zeta^{-r}y^{\sfrac{1}{T}}$ if $r\not\equiv0\pmod N$.  
  In fact, for such a $r$ there is no tensor position where both the states $\tau^r\mathfrak{x}_\alpha$ and $\mathfrak{x}_\beta$ have nontrivial factors, so $Y(\tau^r\mathfrak{x}_\alpha,x^N-y^N)\mathfrak{x}_\beta$ has no negative powers of $x^N-y^N$. 
  On the other hand, the diagonal expansions of $\Delta_N(x^N)^{-1}a$ and $\Delta_N(y^N)^{-1}b$ at any branch have no negative powers of $x^N-y^N$.
  Hence applying \zcref{lem:field-products:b} to $F$, with the substitutions $z=x^N$ and $w=y^N$, gives a local expansion of $\mathcal{S}_N[a,x]\mathcal{S}_N[b,y]$ with no negative powers of $x^N-y^N$ at the branch $x^{\sfrac{1}{T}}=\zeta^{-r}y^{\sfrac{1}{T}}$. 
  The same holds for $\mathcal S_N[b,y]\mathcal S_N[a,x]$.
  Binomial expansion in $x^{\sfrac{1}{T}}-\zeta^{-r}y^{\sfrac{1}{T}}$ preserves total degree in the root variables. 
  Thus each finite Laurent polynomial above is divisible by $(x^{\sfrac{1}{T}}-\zeta^{-r}y^{\sfrac{1}{T}})^K$, and this factor cancels from its two ordered quotients whenever $r\not\equiv0\pmod N$.
  The remaining factor is thus
  \[
    \prod_{q=0}^{T-1}(x^{\sfrac{1}{T}}-\zeta^{-Nq}y^{\sfrac{1}{T}})^K=(x-y)^K.
  \]
  This gives the desired locality relation modulo every $W$, hence in $R$ by separatedness:
  \[
    (x-y)^K\mathcal{S}_N[a,x]\mathcal{S}_N[b,y]
    =(x-y)^K\mathcal{S}_N[b,y]\mathcal{S}_N[a,x].
  \]

  It remains to verify the $h$-twisted Jacobi identity.
  On the branch $x^{\sfrac{1}{T}}=y^{\sfrac{1}{T}}$, the expansion of $(x^N-y^N)/(x-y)$ is a unit in $\k\dparen{y^{1/T}}\dbrack{x-y}$. So we can rewrite the local expansion from \zcref{lem:field-products:b} in $x-y$ rather than in $x^N-y^N$. 
  Using such an expansion, we obtain
  \[
    \begin{aligned}
      \iota_{y,x-y;0}\ab(\mathcal{S}_N[a,x]\mathcal{S}_N[b,y])
      \overset{\text{\eqref{eq:inverse-transport-field}}}&=
        \iota_{y,x-y;0}\ab(
          F[\mathfrak{x}_{\Delta_N(x^N)^{-1}a},x^N]
          F[\mathfrak{x}_{\Delta_N(y^N)^{-1}b},y^N])\\
      \overset{\text{\zcref[noname]{lem:field-products:b}}}&=
        F[Y_{V^{\otimes N}}(
          \mathfrak{x}_{\Delta_N(x^N)^{-1}a},x^N-y^N)
          \mathfrak{x}_{\Delta_N(y^N)^{-1}b},y^N]\\
      &=F[\mathfrak{x}_{
          Y(\Delta_N(x^N)^{-1}a,x^N-y^N)\Delta_N(y^N)^{-1}b},y^N]\\
      \overset{\text{\zcref[noname]{lem:delta-covariance}}}&=
        F[\mathfrak{x}_{\Delta_N(y^N)^{-1}Y(a,x-y)b},y^N]\\
      &=\mathcal{S}_N[Y(a,x-y)b,y].
    \end{aligned}
  \]
  Together with locality and $h$-equivariance, this identity gives the $h$-twisted Jacobi identity by the argument following \eqref{eq:jacobi-diagonal-expansion} in the proof of \zcref{prop:jacobi-criterion}.
\end{proof}

For the rest of this section, take $R=\mathscr U_\tau(V^{\otimes N})$ and
$F[-,z]=Y_\tau[-,z]$ from \zcref{conv:formal-vertex-operator}.

\begin{construction}\label{con:psi}
  For $n\in\sfrac{1}{T}\Z$ and homogeneous $a\in V$, define $\psi_0(J_n(a))$ to be the coefficient of $x^{-n-\wt a}$ in $\mathcal{S}_N[a,x]$. 
  This coefficient lies in $\mathscr{U}_{\tau}(V^{\otimes N})_{-\sfrac{n}{N}}$. 
  By linearity and freeness, this defines a unique unital algebra homomorphism
  \[
    \nota{\psi_0}\colon\mathbb{U}_{\sfrac{1}{T}}(V)
    \longrightarrow\mathscr{U}_{\tau}(V^{\otimes N}).
  \]
  Note that, for general $a\in V$, this homomorphism is explicitly given by
  \[
    \psi_0(J_n(a))=J_{\sfrac{n}{N}}(\mathfrak{x}_{\Delta_N(1)^{-1}a}).
  \]
\end{construction}

\begin{lemma}\label{lem:continuous-psi}
  The assignment $\psi_0$ is continuous and extends to the degreewise completion.
\end{lemma}
\begin{proof}
  Since $\psi_0$ divides degrees by $N$, multiplicativity gives
  \[
    \psi_0(\cNL[s]\mathbb{U}_{\sfrac{1}{T}}(V))
    \subset\cNL[\sfrac{s}{N}]\mathscr{U}_{\tau}(V^{\otimes N})
    \txand
    \psi_0(\cNR[s]\mathbb{U}_{\sfrac{1}{T}}(V))
    \subset\cNR[\sfrac{s}{N}]\mathscr{U}_{\tau}(V^{\otimes N}).
  \]
  Hence $\psi_0$ is continuous by \zcref[noname]{A2}. The target is separated and degreewise complete, so $\psi_0$ extends uniquely to the degreewise completion.
\end{proof}

\begin{proposition}\label{prop:psi-homomorphism}
  The assignment $\psi_0$ induces a continuous graded homomorphism
  \[
    \nota{\psi}\colon\mathfrak{d}_N^{\ast}\mathscr{U}_{g^N}(V)
    \longrightarrow\mathscr{U}_{\sigma g}(V^{\otimes N}).
  \]
\end{proposition}
\begin{proof}
  The continuous extension of $\psi_0$ given by \zcref{lem:continuous-psi} commutes with the fixed-degree sums in the defining Jacobi relations.
  By \zcref{prop:inverse-field}, the vacuum, equivariance and twisted Jacobi relations therefore map to zero.
  Since $\mathscr{U}_{\tau}(V^{\otimes N})$ is separated, the kernel of this extension is a closed two-sided ideal containing $\widehat{\mathbb{I}}_{h}$. Thus this extension descends continuously to $\mathscr{U}_{h}(V)$.
\end{proof}

\section{Main results and applications}
\label{sec:conclusion}

We show that the homomorphisms constructed in the preceding two sections are inverse to each other and extend the results to arbitrary permutations. 
We then derive consequences for module categories, level algebras and mode transition algebras.

\subsection{Main theorems}

We begin with the cyclic case.

\begin{theorem}\label{thm:main-iso}
  Let $\sigma$ be a cyclic permutation of length $N$, and let $g$ be a finite-order automorphism of $V$. 
  Then $\varphi$ and $\psi$ give mutually inverse isomorphisms of almost-canonically seminormed rings
  \[
    \begin{tikzcd}[column sep=large]
      \mathscr{U}_{\sigma g}(V^{\otimes N})
        \arrow[r, shift left, "\varphi"] &
      \mathfrak{d}_N^{\ast}\mathscr{U}_{g^N}(V).
        \arrow[l, shift left, "\psi"]
    \end{tikzcd}
  \]
\end{theorem}

\begin{proof}
  The formulas in \zcref[noname]{con:phi,con:psi} give $\varphi\psi(J_n(a))=J_n(a)$ and $\psi\varphi(J_{\sfrac{n}{N}}(\mathfrak{x}_a))=J_{\sfrac{n}{N}}(\mathfrak{x}_a)$. 
  Then equivariance gives $\psi\varphi(J_{\sfrac{n}{N}}(\mathfrak{x}^i_a))=J_{\sfrac{n}{N}}(\mathfrak{x}^i_a)$ for $1\le i\le N$. 
  By its construction, $\mathscr{U}_{g^N}(V)$ is topologically generated by the modes $J_n(a)$, and \zcref{lem:tensor-generation} shows that $\mathscr{U}_{\sigma g}(V^{\otimes N})$ is topologically generated by the modes of very simple tensors.
  Continuity therefore makes both composites the identity of the algebras.
\end{proof}

A general permutation $\sigma\in\mathfrak{S}_N$ can be written as
\begin{equation}\label{eq:cycle-decomposition}
  \sigma = \pi ( \sigma_1 \oplus \cdots \oplus \sigma_k ) \pi^{-1},
\end{equation}
where $\pi\in\mathfrak{S}_N$, each $\sigma_i\in\mathfrak{S}_{N_i}$ is a cycle of length $N_i$, and $N_1+\cdots+N_k=N$.
Here we use the \emph{direct sum of permutations}
\[
  \mathfrak{S}_{N_1}\times\cdots\times\mathfrak{S}_{N_k}
  \longrightarrow\mathfrak{S}_N
  \colon
  (\sigma_1,\cdots,\sigma_k)
  \longmapsto \sigma_1\oplus\cdots\oplus\sigma_k,
\]
where $\sigma_1\oplus\cdots\oplus\sigma_k$ acts as $\sigma_i$ on the $i$-th consecutive block of $N_i$ positions.
By \zcref{lem:conjugation-invariance}, we may assume $\sigma_i = (1 2 \cdots N_i)$.

\begin{theorem}\label{thm:main-iso-general}
  Let $\sigma$ be as in \eqref{eq:cycle-decomposition}, and let $g$ be a finite-order automorphism of $V$.
  Then there is an isomorphism of almost-canonically seminormed rings
  \[
    \mathscr{U}_{\sigma g}(V^{\otimes N})
    \cong
    \mathfrak{d}_{N_1}^{\ast}\mathscr{U}_{g^{N_1}}(V)
    \ctensor\cdots\ctensor
    \mathfrak{d}_{N_k}^{\ast}\mathscr{U}_{g^{N_k}}(V).
  \]
\end{theorem}
\begin{proof}
  Choose a period $T$ of $g$ and put $N_{\sigma}=\operatorname{lcm}(N_1,\ldots,N_k)$.
  By \zcref{rem:enveloping-choice-of-T}, we may use the period $N_{\sigma}T$ for $\sigma g$ and the period $N_{\sigma}T/N_i$ for $g$ in the application of \zcref{thm:main-iso} to the $i$-th cycle, with compatible choices of roots of unity.
  These cyclic isomorphisms preserve degrees and neighborhoods, so they induce an isomorphism of degreewise completed tensor products by \zcref{recol:acsr-tensor-products}.
  Since the diagonal action of $g$ commutes with $\pi$, we obtain
  \begin{align*}
    \mathscr{U}_{\sigma g}(V^{\otimes N})
    \overset{\text{\zcref{lem:conjugation-invariance}}}&{\cong}
    \mathscr{U}_{(\sigma_1\oplus\cdots\oplus\sigma_k)g}(V^{\otimes N}) \\
    \overset{\text{\zcref{prop:enveloping-tensor-product}}}&{\cong}
    \mathscr{U}_{\sigma_1g}(V^{\otimes N_1}) \ctensor \cdots \ctensor \mathscr{U}_{\sigma_kg}(V^{\otimes N_k}) \\
    \overset{\text{\zcref{thm:main-iso}}}&{\cong}
    \mathfrak{d}_{N_1}^{\ast}\mathscr{U}_{g^{N_1}}(V) \ctensor \cdots \ctensor \mathfrak{d}_{N_k}^{\ast}\mathscr{U}_{g^{N_k}}(V).
    \qedhere
  \end{align*}
\end{proof}

\begin{remark}
  We cannot in general apply \zcref{prop:enveloping-tensor-product} to identify the right-hand side of \zcref{thm:main-iso-general} with $\mathscr{U}_{g^{N_1}\otimes\cdots\otimes g^{N_k}}(V^{\otimes k})$.
  The proposition uses the total grading, whereas the theorem rescales the factors before taking the total grading and completing degreewise.
  For unequal cycle lengths, these completions need not agree. See \zcref{rem:rescaled-tensor-products}.
\end{remark}

\subsection{Level algebras and mode transition algebras}
We first consider the cyclic case, then use the tensor-product decomposition to handle arbitrary permutations.

The following corollary generalizes the result in \cite{DXY22} about (twisted) Zhu algebras.
\begin{corollary}\label{cor:higher-zhu-isomorphism}
  Let $\sigma$ be a cyclic permutation of length $N$, let $g$ be a finite-order automorphism of $V$, and let $n\in\sfrac{1}{T}\N$. Then, $A_{\sigma g, \sfrac{n}{N}}(V^{\otimes N})\cong A_{g^N, n}(V)$.
\end{corollary}
\begin{proof}
  By \zcref{rem:zhu-comparison}, the assertion follows from the quotient isomorphism induced by $\varphi$:
  \[
    \frac{\mathscr{U}_{\sigma g}(V^{\otimes N})_0}{\NL[\sfrac{n}{N}+\sfrac{1}{NT}]\mathscr{U}_{\sigma g}(V^{\otimes N})_0}
    \overset{\varphi}{\cong}
    \frac{\mathscr{U}_{g^N}(V)_0}{\NL[n+\sfrac{1}{T}]\mathscr{U}_{g^N}(V)_0}.\qedhere
  \]
\end{proof}

\begin{lemma}\label{lem:tensor-filtration-graded}
  Use the setting of \zcref{thm:main-iso-general} and put
  $N_{\sigma}=\operatorname{lcm}(N_1,\ldots,N_k)$. For
  $n\in\sfrac1{N_{\sigma}T}\N$ and $p\in\sfrac1{N_{\sigma}T}\Z$, there are canonical isomorphisms
  \begin{align*}
    \mathfrak{L}^{\gr n}_p(\mathscr{U}_{\sigma g}(V^{\otimes N}))
    &\cong\bigoplus_{\substack{n_i\in\sfrac{1}{T}\N,\ p_i\in\sfrac{1}{T}\Z\\
                 \sum_i \sfrac{n_i}{N_i}=n,\ \sum_i \sfrac{p_i}{N_i}=p}}
       \bigotimes_i\mathfrak{L}^{\gr n_i}_{p_i}(\mathscr{U}_{g^{N_i}}(V)),\\
    \mathfrak{R}^{\gr n}_p(\mathscr{U}_{\sigma g}(V^{\otimes N}))
    &\cong\bigoplus_{\substack{n_i\in\sfrac{1}{T}\N,\ p_i\in\sfrac{1}{T}\Z\\
                 \sum_i \sfrac{n_i}{N_i}=n,\ \sum_i \sfrac{p_i}{N_i}=p}}
       \bigotimes_i\mathfrak{R}^{\gr n_i}_{p_i}(\mathscr{U}_{g^{N_i}}(V)).
  \end{align*}
  Only finitely many nonzero summands occur.
\end{lemma}
\begin{proof}
  On the common lattice $\sfrac{1}{N_{\sigma}T}\Z$, \zcref{thm:main-iso-general,recol:acsr-tensor-products} identify the filtrations with the convolution filtrations determined by the factors $\mathfrak{d}_{N_i}^{\ast}\mathscr{U}_{g^{N_i}}(V)$.
  The associated graded pieces with negative filtration index vanish because $\NL[0]U=\NR[0]U=U$.
  Hence, taking associated graded spaces of these convolution filtrations gives the displayed isomorphisms.

  There are only finitely many nonzero summands. Indeed, note that
  \[
    \mathfrak{L}^{\gr n_i}_{p_i}=0\quad\text{if }p_i<-n_i,
    \qquad
    \mathfrak{R}^{\gr n_i}_{p_i}=0\quad\text{if }p_i>n_i.
  \]
  For fixed $n$ and $p$, there are only finitely many choices for the nonnegative $n_i$ and, by these bounds, for the $p_i$.
\end{proof}

\begin{corollary}\label{cor:mode-transition-tensor}
  In the setting of \zcref{thm:main-iso-general}, there is an isomorphism of bigraded algebras
  \[
    \mathfrak{A}(\mathscr{U}_{\sigma g}(V^{\otimes N}))
    \cong\bigotimes_{i=1}^k
       \mathfrak{A}(\mathfrak{d}_{N_i}^{\ast}\mathscr{U}_{g^{N_i}}(V)).
  \]
\end{corollary}
\begin{proof}
  Applying \cite[Theorem D]{SIC} to the tensor product in \zcref{thm:main-iso-general} gives the algebra isomorphism.
  The level-zero case of \zcref{lem:tensor-filtration-graded} identifies its bigraded components.
\end{proof}

\begin{corollary}\label{cor:zhu-tensor}
  With the same setting as in \zcref{thm:main-iso-general}, there is an isomorphism of twisted Zhu algebras
  \[
    A_{\sigma g}(V^{\otimes N}) \cong A_{g^{N_1}}(V) \otimes \cdots \otimes A_{g^{N_k}}(V).
  \]
\end{corollary}
\begin{proof}
  In \zcref{cor:mode-transition-tensor}, each factor on the right
  is supported in bidegrees $(p,-q)$ with $p,q\ge0$.
  Taking the $(0,0)$-component therefore gives the assertion by
  \zcref{def:mode-transition-algebra,rem:zhu-comparison}.
\end{proof}

\begin{corollary}\label{cor:sic-general}
  In the setting of \zcref{thm:main-iso-general}, put
  $h=g^{N_1}\otimes\cdots\otimes g^{N_k}$ on $V^{\otimes k}$.
  Then $\mathscr{U}_{\sigma g}(V^{\otimes N})$ satisfies the strong identity condition if and only if $\mathscr{U}_h(V^{\otimes k})$ does.
  These conditions hold if every $\mathscr{U}_{g^{N_i}}(V)$ satisfies the strong identity condition.
  If all these factor algebras are nonzero, the converse holds as well.
\end{corollary}

  When $\sigma$ is an $N$-cycle, the condition is equivalent to the strong identity condition for $\mathscr{U}_{g^N}(V)$.

\begin{proof}
  By \zcref{thm:main-iso-general,prop:enveloping-tensor-product}, $\mathscr{U}_{\sigma g}(V^{\otimes N})$ and $\mathscr{U}_h(V^{\otimes k})$ are the degreewise completed tensor products of the factors $\mathfrak{d}_{N_i}^{\ast}\mathscr{U}_{g^{N_i}}(V)$ and $\mathscr{U}_{g^{N_i}}(V)$, respectively.
  If any factor $\mathscr{U}_{g^{N_i}}(V)$ is zero, both $\mathscr{U}_{\sigma g}(V^{\otimes N})$ and $\mathscr{U}_h(V^{\otimes k})$ are zero and trivially satisfy the strong identity condition.
  Otherwise, uniform rescaling preserves the strong identity condition of each factor, so \zcref{recol:acsr-tensor-products} identifies the condition for either completed tensor product with that for every factor algebra $\mathscr{U}_{g^{N_i}}(V)$.
\end{proof}

If $\one_p^{(i)}$ are strong identity elements of $\mathscr{U}_{g^{N_i}}(V)$, then under \zcref{cor:mode-transition-tensor} the $\sigma$-twisted enveloping algebra $\mathscr{U}_{\sigma g}(V^{\otimes N})$ has strong identity elements
  \[
    \one_s=
    \sum_{\substack{p_i\in\sfrac{1}{T}\N\\ \sum_i p_i/N_i=s}}
      \bigotimes_i\one_{p_i}^{(i)}
    \qquad\ab(s\in\sfrac1{N_{\sigma}T}\N),
  \]
  where $N_{\sigma}=\operatorname{lcm}(N_1,\ldots,N_k)$.

\subsection{Module categories and rationality}

We now apply the algebra isomorphisms and the corresponding field constructions to study the module categories.

In the setting of \zcref{thm:main-iso-general}, put
\[
  \nota{\tau}=\sigma g,
  \qquad \nota{h}=g^{N_1}\otimes\cdots\otimes g^{N_k},
  \qquad \nota{N_{\sigma}}=\operatorname{lcm}(N_1,\ldots,N_k).
\]
Write $\nota{J_n^{(i)}(a)}=1\otimes\cdots\otimes J_n(a)\otimes\cdots\otimes1$,
where $J_n(a)\in\mathscr{U}_{g^{N_i}}(V)$ occurs in position $i$,
and set $\nota{\VL{0}^{(i)}}=J_0^{(i)}(\cfv)$.
By \zcref{prop:module-tensor-products}, a weak $h$-twisted
$V^{\otimes k}$-module is equivalently a vector space carrying
pairwise commuting weak $g^{N_i}$-twisted $V$-actions.
We use the same symbols for the operators of these factor actions.

\begin{lemma}\label{lem:virasoro-zero-mode}
  Let $\cch$ be the central charge of $V$.
  Under the isomorphism in \zcref{thm:main-iso-general}, the Virasoro operator
  $\VL{0}^{\tau} = J_0(\sum_{j=1}^{N}\mathfrak{x}_{\cfv}^{j}) \in \mathscr{U}_{\tau}(V^{\otimes N})$ maps to
  \[
    \sum_{i=1}^{k}\frac{\VL{0}^{(i)}}{N_i}
    +\frac{\cch}{24}\sum_{i=1}^{k}
      \ab(N_i-\frac{1}{N_i})\,
      1_{\mathscr{U}_{g^{N_1}}(V)} \otimes\cdots\otimes 1_{\mathscr{U}_{g^{N_k}}(V)}.
  \]
\end{lemma}
\begin{proof}
  The conjugation in the proof of \zcref{thm:main-iso-general} fixes the tensor-product conformal vector.
  Under the block decomposition there, $\VL{0}^{\tau}$ is the sum of the Virasoro operators of the blocks.
  For the block of length $N_i$, the definition of $\Delta_{N_i}$ and the Virasoro relations give
  \[
    \Delta_{N_i}(1)\cfv
    =\frac{1}{N_i^2}\ab(\cfv+\frac{\cch(N_i^2-1)}{24}\vac).
  \]
  Apply \zcref{con:phi} with $n=0$ to its $N_i$ tensor positions.
  Since $g\cfv=\cfv$, the Virasoro operator of this block maps to
  \[
    N_i\,J_0(\Delta_{N_i}(1)\cfv)
    =\frac{1}{N_i}J_0(\cfv)
      +\frac{\cch}{24}\ab(N_i-\frac{1}{N_i})\,
        1_{\mathscr{U}_{g^{N_i}}(V)}.
  \]
  Summing these images in their respective tensor factors gives the stated formula.
\end{proof}

The following recovers the weak, admissible and ordinary module equivalences in \cite{BDM02}.

\begin{theorem}\label{thm:module-equivalences-general}
  The cycle constructions in the preceding two sections give the following:
  \begin{enumerate}[beginpenalty=10000]
    \item\label{thm:module-equivalences-general:weak} 
        There are equivalences
          \[
            \cat{WMod}_{\tau}(V^{\otimes N})
            \simeq\bigboxtimes_{i=1}^k\cat{WMod}_{g^{N_i}}(V)
            \simeq\cat{WMod}_{h}(V^{\otimes k}).
          \]
    \item\label{thm:module-equivalences-general:exhaustive} 
        These equivalences restrict to
          \[
            \cat{ExMod}_{\tau}(V^{\otimes N})
            \simeq\bigboxtimes_{i=1}^k\cat{ExMod}_{g^{N_i}}(V)
            \simeq\cat{ExMod}_{h}(V^{\otimes k}).
          \]
    \item\label{thm:module-equivalences-general:admissible} 
        Under the equivalence, an exhaustive $h$-twisted module $M$ corresponds to an admissible $\tau$-twisted module if and only if $M$ admits a $\sfrac1{N_{\sigma}T}\N$-grading such that
          \[
            J_n^{(i)}(a)M(p)\subset M\ab(p-\tfrac{n}{N_i})
            \quad\ab(n\in\sfrac{1}{T}\Z).
          \]
    \item\label{thm:module-equivalences-general:ordinary} 
        An exhaustive $h$-twisted module $M$ corresponds to an ordinary $\tau$-twisted module if and only if the operator
          \[
            \tfrac{\VL{0}^{(1)}}{N_1}+\cdots+\tfrac{\VL{0}^{(k)}}{N_k}\Big|_M
          \]
          is diagonalizable over $\k$, with finite-dimensional eigenspaces $M_\lambda$ satisfying $M_{\lambda+n}=0$ for each $\lambda\in\k$ and all sufficiently negative $n\in\sfrac{1}{N_{\sigma}T}\Z$.
    \item\label{thm:module-equivalences-general:geometrically-ordinary} 
        An exhaustive $h$-twisted module $M$ corresponds to a geometrically ordinary $\tau$-twisted module if and only if the condition in \itemref{thm:module-equivalences-general:ordinary} holds after a base change of $\k$.
  \end{enumerate}
  If all cycle lengths are equal, the conditions in \itemref{thm:module-equivalences-general:admissible,thm:module-equivalences-general:ordinary,thm:module-equivalences-general:geometrically-ordinary} become the usual admissible, ordinary and geometrically ordinary conditions after uniform rescaling.
\end{theorem}
\begin{proof}
  We first prove \itemref{thm:module-equivalences-general:weak,thm:module-equivalences-general:exhaustive}, and start with the case where $k=1$. Let $M$ be a weak module on either side. By \zcref{rem:weak-module-distributions}, the constructions in \zcref{con:tensor-field,prop:inverse-field} apply with $R=\End(M)$ and give a weak module on the other side.
  For $g^N$-twisted $M$, the composite recovers all the original fields on $V$, since $\Delta_N$ and its inverse cancel.
  For $\tau$-twisted $M$, the same cancellation recovers the original fields on the first tensor position of $V^{\otimes N}$.
  The original field assignment and the one obtained after the two constructions are $\tau$-equivariant, so agreement extends to all very simple tensors and then to all states by \zcref{lem:tensor-generation}.
  On exhaustive modules, these constructions coincide with pullback along $\varphi$ and $\psi$, since their mode actions agree.
  As for the general case, \itemref{thm:module-equivalences-general:weak} follows from
  \[
    \cat{WMod}_{\tau}(V^{\otimes N})
    \simeq \bigboxtimes_{i=1}^k\cat{WMod}_{\sigma_i g}(V^{\otimes N_i})
    \simeq \bigboxtimes_{i=1}^k\cat{WMod}_{g^{N_i}}(V)
    \simeq \cat{WMod}_h(V^{\otimes k}),
  \]
  where the first and last equivalences are given by \zcref{prop:module-tensor-products} and the middle one is the cyclic case.
  By \zcref{rem:rescaled-tensor-products}, the same factorwise
  constructions give \itemref{thm:module-equivalences-general:exhaustive}.

  For \itemref{thm:module-equivalences-general:admissible}, it is due to that $J_n^{(i)}(a)$ has degree $-n/N_i$ in the $i$-th rescaled factor of \zcref{thm:main-iso-general}.
  Assertion \itemref{thm:module-equivalences-general:ordinary} follows from \zcref{lem:virasoro-zero-mode}. Indeed, the displayed operator and the Virasoro operator on the corresponding $\tau$-twisted module have the same eigenspaces, with all eigenvalues shifted by $\tfrac{\cch}{24}\sum_{i=1}^{k}\ab(N_i-\tfrac1{N_i})$, which preserves ordinarity.
  Since the constructions commute with base changes, \itemref{thm:module-equivalences-general:geometrically-ordinary} follows.
\end{proof}

\begin{corollary}\label{cor:rationality-general}
  In the setting of \zcref{thm:module-equivalences-general}, $V^{\otimes N}$ is $\tau$-rational if and only if $V^{\otimes k}$ is $h$-rational.
  These conditions hold if $V$ is $g^{N_i}$-rational for every $i$.
  If every $\mathscr{U}_{g^{N_i}}(V)$ is nonzero, the converse holds as well.
\end{corollary}
\begin{proof}
  By \zcref{thm:main-iso-general,prop:enveloping-tensor-product,rem:exhaustive-module-correspondence}, the two rationality conditions are $\k$-rationality of the degreewise completed tensor products with factors $\mathfrak{d}_{N_i}^{\ast}\mathscr{U}_{g^{N_i}}(V)$ and $\mathscr{U}_{g^{N_i}}(V)$, respectively.
  If $\mathscr{U}_{g^{N_i}}(V)=0$ for some $i$, both tensor products vanish and both rationality conditions hold trivially.
  Otherwise, since uniform rescaling preserves $\k$-rationality, iterating \zcref{lem:tensor-k-rationality} makes each condition equivalent to $\k$-rationality of every $\mathscr{U}_{g^{N_i}}(V)$, hence to $g^{N_i}$-rationality of $V$ for every $i$.
\end{proof}

\bibliographystyle{amsalpha}
\bibliography{bib}

\end{document}